\documentclass{elsarticle_mod}
\usepackage[english]{babel}

\usepackage{amsmath}
\usepackage{amssymb,amsfonts}
\usepackage{amsthm}
\usepackage{rotating}
\usepackage{graphicx}
\usepackage{float}
\usepackage{multirow}
\usepackage{algorithm,algpseudocode} 
\usepackage{color}
\usepackage{booktabs} 
\usepackage{verbatim}
\usepackage{textcomp}
\usepackage{subfigure}

\newtheorem{proposition}{Proposition}

\newcommand{\argmin}{\operatorname*{arg\, \min}}
\newcommand{\R}{\mathbb{R}}

\newcommand{\ve}[1]{\mbox{\boldmath $ #1$}}
\newcommand{\bC}{\boldsymbol{C}}
\newcommand{\bM}{\boldsymbol{M}}
\newcommand{\be}{\boldsymbol{e}}
\newcommand{\bk}{\boldsymbol{k}}

\begin{document}
\begin{frontmatter}

\title{Computational approaches for parametric imaging of dynamic PET data}

		\author[address1,address4]{Serena Crisci}
\ead{serena.crisci@unife.it}
\address[address1]{Dipartimento di Matematica e Informatica, Universit\`{a} di Ferrara, Ferrara, Italy}
\author[address2,address4]{Michele Piana}
\ead{piana@dima.unige.it}
\address[address2]{ Dipartimento di Matematica, Universit\`{a} di Genova, and CNR-SPIN Genova,Italy}
\author[address1,address4]{Valeria Ruggiero}
\ead{valeria.ruggiero@unife.it}
\author[address3]{Mara Scussolini}
\ead{scussolini@dima.unige.it}
\address[address3]{Dipartimento di Matematica, Universit\`{a} di Genova, Genova,Italy}
\address[address4]{Member of the INdAM Research group GNCS}

\vspace{10pt}

	\begin{abstract}
Parametric imaging of nuclear medicine data exploits dynamic functional images in order to reconstruct maps of kinetic parameters related to the metabolism of a specific tracer injected in the biological tissue. From a computational viewpoint, the realization of parametric images requires the pixel-wise numerical solution of compartmental inverse problems that are typically ill-posed and nonlinear. In the present paper we introduce a fast numerical optimization scheme for parametric imaging relying on a regularized version of the standard affine-scaling Trust Region method. The validation of this approach is realized in a simulation framework for brain imaging and comparison of performances is made with respect to a regularized Gauss-Newton scheme and a standard nonlinear least-squares algorithm.
	\end{abstract}

\end{frontmatter}
\vspace{2pc}

	\section{Introduction}
	
	Positron Emission Tomography (PET) \cite{baetal05} utilizes an isotope produced in a cyclotron to provide dynamical images of the metabolism-based isotope accumulation in the biological tissue. PET dynamic images of the tracer distribution are obtained by applying a reconstruction algorithm to the measured radioactivity and provide a reliable estimate of the tracer-related metabolism in the tissue \cite{nawu01,shva82}.
	
	From a technical viewpoint, compartmental analysis \cite{guetal01,sctu02,waetal06} allows processing these dynamic PET data in order to estimate a set of physiological kinetic parameters that explain such metabolism in a quantitative manner (specifically, these parameters express the effectiveness of the tracer in changing its functional status within the tissue). 
	Compartmental analysis requires, first, the formulation of a forward model for the tracer concentration represented by a Cauchy problem, in which the kinetic parameters are the coefficients of the differential equations; then, the numerical solution of the corresponding nonlinear inverse problem, in which the kinetic parameters are the unknowns and the tracer concentrations in the tissue are the input data.
	
	Relying on compartmental analysis, parametric imaging \cite{guetal97,kaetal13,reve14} allows the pixel-wise determination of the kinetic parameters by means of two possible alternative approaches. On the one hand, direct parametric imaging \cite{kaetal05,waqi13} utilizes as input the PET raw sinograms and solves the inverse problem that relates them to the parameters; on the other hand, indirect parametric imaging \cite{guetal97,kaetal13,caetal17,Scussolini} is applied to the reconstructed PET images and solves pixel-wise the compartmental analysis problem. Direct approaches do not need the application of image reconstruction methods but have typically to deal with the intertwining of spatial and temporal correlations, which makes the optimization process more complex; this same optimization is more straightforward in indirect approaches but requires a higher computational burden, due to the need of solving a large number of nonlinear inverse problems.\newline
	
	\noindent The present paper deals with indirect parametric imaging and introduces a regularized optimization method for the solution of the nonlinear ill-posed inverse problem of compartmental analysis. The idea of the method is to introduce a regularizing strategy \cite{Wang} in the standard affine-scaling Trust Region method \cite{Coleman,Bellavia_2006}, which allows a better reduction of the numerical instabilities induced by the presence of the experimental noise in the measured data. 
	
	From a formal viewpoint, we prove a convergence result for the regularized algorithm, which enables a generalization to the non-negatively constrained case of the
convergence analysis developed in \cite{Wang} for the unconstrained problem. The numerical validation of the method is performed against synthetic data realized from an 'ad hoc' modification of the Hoffman Brain Phantom often used in PET and CT imaging (http://depts.washington.edu/petctdro/DROhoffman\_main.html). Specifically, we mimicked a two-compartment experiment for the kinetics of [18F]fluoro-2-deoxy-D-glucose (FDG), which is the mostly utilized tracer in PET diagnostic and prognostic activities \cite{boetal10,fletal01,kaetal16,maetal18,maetal13}. Using this simulation we could compare the computational effectiveness and reconstruction accuracy of the method with respect to the performances provided by two frequently used indirect parametric imaging methods. 
	
	The structure of the paper is as follows. Section 2 sets up the two-compartment problem for FDG kinetics. Section 3 describes in detail the nonlinear optimization method for the solution of this problem. Section 4 illustrates the validation experiment and its results. Our conclusions are offered in Section 5.
	
	\section{Compartmental analysis of dynamic PET data} \label{sec:comp-analysis}
	
\begin{figure}
		\centering
		\includegraphics[width=0.5\textwidth]{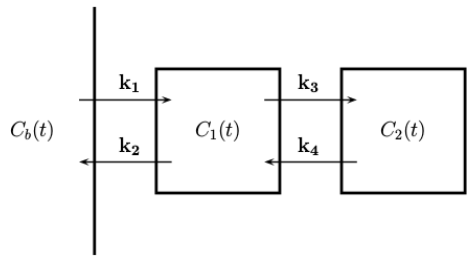}
		\caption{Compartmental model for FDG kinetics in a living tissue: external arterial blood compartment of concentration $C_b$, free tracer compartment of concentration $C_f$, metabolized tracer compartment of concentration $C_m$, and four kinetic parameters $k_1,k_2,k_3,k_4$.}
		\label{fig:model-2C}
	\end{figure}

	Compartmental analysis of nuclear medicine data is the mathematical framework for the quantitative assessment of tracer kinetics in the biological tissue \cite{Carson,Cherry,Morris,Schmidt,gaetal15,gaetal14}. The compartmental model of a specific organ comprises compartments representing the functional states of the tracer radioactive molecules
	(e.g.: physical location as intravascular space, extracellular space, intracellular space, etc., or chemical state as metabolic form, binding state, etc.),	and kinetic parameters, which are the input/output tracer rates for each compartment. Figure \ref{fig:model-2C} illustrates the standard two-compartmental model describing the FDG metabolism in the organ under consideration \cite{Sokoloff}. This model reproduces the main steps of the FDG path in a PET experiment. First, the tracer is injected into the blood with a concentration mathematically modelled by the Input Function (IF), here assumed as known and represented by the tracer concentration $C_b$ in the arterial blood compartment. Then, the FDG metabolism within the tissue is characterized by two functional states: the free compartment with concentration $C_f$, associated to the tracer molecules outside the tissue cells, and the metabolized compartment with concentration $C_m$, associated to FDG molecules within the cytoplasm. Finally, the FDG kinetics is described by four rate constants connecting the model compartments: $k_1$ and $k_2$ describe the exchange rates between the input and free pools, and $k_3$ and $k_4$ describe the exchange rates at the basis of the phosphorylation/dephosphorylation process.


	The system of Ordinary Differential Equations (ODEs) for the two-compartment model is
	\begin{equation}\label{eq:ODE_2C}
	\frac{d\bC}{dt}(t) = \dot{\bC}(t) = \bM \bC(t) + k_1 C_b(t) \be_1 ,
	\end{equation}
	where
	\begin{equation}\label{eq:ODE_2C_matrix}
	\bC = \begin{pmatrix} C_f \\ C_m \end{pmatrix} , \quad
	\bM = \begin{pmatrix} -(k_2+k_3) & k_4 \\ k_3 & -k_4 \end{pmatrix} , \quad
	\be_1 = \begin{pmatrix} 1 \\ 0 \end{pmatrix} ,
	\end{equation}
	where $t$ 
		is the time variable, and $\pm k_i$, $i =1,2,3,4$, represent incoming and outgoing fluxes.
	In standard applications of compartmental analysis, the initial conditions are $C_f(0) = C_m(0) = 0$, meaning that the PET experiment starts at time $t = 0$ when there is no available tracer into the biological system.
	The analytical solution of \eqref{eq:ODE_2C} represents the forward model equation of determining the compartment concentrations given
	the kinetic parameters, and takes the form
	\begin{equation}\label{eq:sol_2C}
	\bC(t;\bk) = k_1 \int_0^t \exp(\bM(t-u)) C_b(u) \be_1 \, du ,
	\end{equation}
	where the entries of the vector $\bk = (k_1,k_2,k_3,k_4)^T \in \R^4$ have to be non-negative real values.
	
	The compartmental input function $C_b(t)$ can be obtained experimentally either from serial sampling of the arterial blood or reconstructed dynamic images \cite{suetal15}, when a large arterial pool such as the left ventricle is in the field of view for many frames, or by using reference tissue methods \cite{saetal15,scetal18}. However, PET images cannot offer enough resolution power to provide information on $\bC(t;\bk)$. Therefore the measurement equation 
	\begin{equation}\label{eq:inv_2C}
	\tilde{C}_{(p,q)}(t;\ve k) = \ve \alpha^T \ve C_{(p,q)}(t;\bk) + V_{(p,q)} C_b(t) \ , \quad \ve \alpha = \begin{pmatrix} 1-V_{(p,q)} \\ 1-V_{(p,q)} \end{pmatrix} ,
	\end{equation}
	should be added to equation (\ref{eq:sol_2C}) to connect the compartment model to the PET data. In this equation $(p,q)$ represents a specific image pixel, $\tilde{C}_{(p,q)}$ denotes the measured tracer concentration at pixel $(p,q)$ of the organ image, $\ve C_{(p,q)}$ is the formal analytic solution of \eqref{eq:sol_2C}, and $V_{(p,q)}$ is the fraction of tissue volume occupied by the blood. In general, the blood volume fraction depends on the pixel position, but within a homogeneous tissue it can be assumed as a known constant.

	In equation (\ref{eq:inv_2C}) the unknown kinetic parameters are functions of $(p,q)$ and therefore the inverse problem represented by this equation should be solved numerically and pixel-wise. Rather coarse approximations allow a linearization of this equation \cite{paetal83,lo00}. However, the pixel-wise solution of the exact nonlinear equation requires the availability of an effective optimization scheme for the regularization of the ill-posed nonlinear compartmental inverse problem represented by the equation and eventually for the reconstruction of the four parametric images associated to $k_1$, $k_2$, $k_3$, and $k_4$.

	\section{Computational approaches for nonlinear ill-posed problems} \label{sec:TR}
	The compartmental inverse  problem described in the previous section is a special case of the following more general formulation.
	Given a set of measurements $\ve y^0$ 
	of tracer concentration provided by PET images, corresponding to a finite sample of $N$ time points
	$t_1,...,t_N$, we have to determine the kinetic parameters $\ve k \in\R^n$, $n\leq N$, by solving the
	non-negatively constrained nonlinear system
	\begin{equation}\label{nonlinearsystem}
	{\ve F}(\ve k)= {\ve{y}^0},\quad \quad \mbox{subject to }{\ve k}\geq 0.
	\end{equation}
	Here $\ve{y}^0 = (\tilde{C}_{(p,q)}(t_1,\ve{k}),...,\tilde{C}_{(p,q)}(t_N,\ve{k}))^T$, and $\ve {F}\colon\R^n\rightarrow \R^N$ 
	is the continuously differentiable function at the right hand side of (\ref{eq:inv_2C}). In real experiments a noisy version $\ve {y}^{\delta}$ of $\ve{y}^0$ is at disposal, where $\delta$ is a known bound on the measurement error, with $\delta\leq \|\ve{y}^0\|$.
	A standard approach to address equation (\ref{nonlinearsystem}) \cite{Engl,Kalten} 
	consists in approximating a solution $\ve k^\dag$ of this nonlinear system by solving the following non-negatively
	nonlinear least squares problem via an iterative regularization technique with semiconvergent behaviour:
	\begin{equation}\label{problem}
	\min_{\ve{k}\geq  0}\  \Phi(\ve k) \equiv \frac{1}{2}\|\ve {y^\delta} - {\ve F}(\ve k)\|^2.
	\end{equation}
	In view of the discrepancy principle \cite{Engl},  the iterative method is stopped at the iteration ${\bar{j}(\delta)}$ satisfying the following condition
	\begin{equation}\label{stopcrit}
	\|\ve {y}^\delta - \ve F(\ve k^{\bar{j}(\delta)})\|\leq \tau \delta <\|\ve {y}^\delta - \ve F(\ve k^{{j}})\|\quad 0\leq j\leq {\bar{j}(\delta)},
	\end{equation}
	for a suitable $\tau>1$.
	
	In this section, we describe a method for computing a regularized solution of problem~\eqref{problem}; in particular, we combine the regularizing approach developed in~\cite{Wang} for unconstrained ill-posed problems with the affine scaling trust-region (TR) schemes for a box-constrained minimization problem~\cite{Coleman,Macconi}. The key point to link these methods is represented by the following Proposition~\ref{mono_proj}, which shows that possible projection
	steps do not prevent the convergence of the iterative scheme. Therefore, the main contribution of this section is to show that the theoretical framework developed for the unconstrained problem \cite{Wang} still holds also in the non-negatively constrained case.
	\subsection{A regularizing affine scaling trust-region method for non-negatively nonlinear least-squares problems} \label{subsec:reg-TR-nonneg}
	For unconstrained nonlinear ill-posed least-squares problems, the state-of-the-art approaches are the regularized Levenberg-Marquadt (LM) method, proposed by Hanke~\cite{Hanke}, and its reformulation within a Trust-Region~(TR) framework, proposed by Wang et al.~\cite{Wang} and, more recently, by Bellavia et al. \cite{Bellavia}.
	As in the standard TR algorithm, the regularizing TR iteration requires to compute,
		at each iteration, a trial step $\ve p^j$, 
		by minimizing the quadratic model $m_j(\ve p)$ within a region around the current iterate $\ve{k}^j$:
	\begin{eqnarray}\label{eqn:tr}
	&&\min_{\ve p}\ m_j(\ve{p})\equiv \frac{1}{2} \ve p^T \ve B^j \ve p+ \ve p^T\ve g^j \nonumber \\
	&& \mbox{s.t. }\ \|\ve p\|\leq \Delta_j
	\end{eqnarray}
	where $\ve B^j\equiv \ve J(\ve{k}^j)^T \ve J(\ve{k}^j)$ is the Gauss-Newton approximation of the Hessian of $\ve F$,
	$\ve g^{j}\equiv \nabla \Phi(\ve{k}^j)=\ve J(\ve{k}^j)^T(\ve F(\ve{k}^j)-\ve{y}^\delta)$  and $\Delta_j$ denotes the TR radius; this can be expanded or reduced depending on whether a sufficient reduction of the model is achieved or not, i.e. if the ratio $\rho_j = \displaystyle \frac{\Phi(\ve k^j + \ve { p}^j)-\Phi(\ve k^j) } {m_j(\ve{ p}^j)}$ between the actual reduction in the objective functional and the predicted reduction in the quadratic model is lower  than some
	positive threshold $\beta\in(0, 1)$.
	The regularizing property is accomplished by requiring that
	the TR constraint is active at the solution,
	i.e., the solution $\ve p^j$ of \eqref{eqn:tr} must be such that $\|\ve p^j\|=\Delta_j$ so that the associated Lagrange multiplier $\alpha_j$ plays the role of a penalization parameter in a Tikhonov-like regularization. Indeed, given $\ve{k}^j$, the new iterate can be viewed as the solution of the penalized subproblem arising at the iteration of LM method:
		\begin{equation}\label{lmHanke}
 \ve k^{j+1} =	\ve{k}^j+ \ve{p}^j =	\argmin_{\ve k} \{ \|\ve {y}^{\delta} - \ve F(\ve k^j) - \ve J(\ve k^j) (\ve k - \ve k^j) \|^2 + \alpha_j \| \ve k -\ve k^j\|^2 \}.
		\end{equation}

This regularization technique for an unconstrained problem can be combined with the TR methods for box-constrained nonlinear least-squares problems. To this aim, we introduce a regularizing technique in the affine-scaling TR method~\cite{Bellavia_2006,Coleman,Macconi}
	requiring that the TR constraint in the subproblem \eqref{eqn:tr} is active at the solution.  In particular,
		given $\ve k^j>0$ and $\ve g^j\neq 0$,
		we find the solution $\alpha_j >0$ of the nonlinear equation $\Delta_j-\|\ve p(\alpha)\|=0$, 
		where $\displaystyle \ve p(\alpha)=(\ve J(\ve k^j)^T \ve J(\ve k^j) +\alpha \ve I_n)^{-1} \ve J(\ve k^j)^T (\ve y^\delta -\ve F(\ve k^j))$. 
By setting $\ve p^j=\ve p(\alpha_j)$, in order to ensure the strict feasibility of a new iterate, 
	the $i$-th entry of $\ve k^{j+1} =	\ve{k}^j+\ve{\bar{p}}^j$ is computed in accordance with the following rule:
		\begin{equation} \label{proj_step}
		\ve{\bar{p}}_i^j= \left\{\begin{array}{cc}
		\ve p_i^j & \quad \mbox{if } (\ve{k}^j+\ve{p}^j)_i> 0 \quad \\
		t(\Pi(\ve{k}^j+\ve{p}^j)-\ve k^j)_i & \mbox{if } (\ve k^j+\ve p^j)_i\leq 0
		\end{array}\right.
		\end{equation}
		where $\Pi(\cdot)$ denotes the Euclidean projection onto the non-negative orthant of $\R^n$ and  $t\in (0,1)$.
		Clearly, in view of the properties of the projection operator, $\|\ve{\bar{p}}^j\|\leq \|\ve p^j\|$.
	
As emphasized in~\cite{Coleman}, a key point to assure 
the convergence 
of the affine-scaling TR method is that the new iterate $\ve{k}^j+\ve{\bar{p}}^j$ 
must be able to achieve at least as much reduction in the quadratic model  as
the one achieved by the generalized Cauchy point $\ve{p}_C^j = - \lambda_C^j  \ve D(\ve k^j) \ve g^j$, 
where  $\ve D(\ve k)$ is a diagonal matrix such that 
	\begin{equation}\label{scaling}
	\ve D(\ve k)_{i,i} = \left\{\begin{array}{ll}
	|\ve k_i| &  \mbox{if }\nabla\Phi(\ve k)_i\geq 0  \\
	1&  \mbox{otherwise}\end{array}\right.
	\end{equation}

and  $\lambda_C^j$  is defined as follows

{\footnotesize \begin{equation} \label{proj_cauchy}
	\hskip -2cm\lambda^j_C= \left\{\begin{array}{ll}
	\min\left\{ \displaystyle \frac{\Delta_j}{\|\ve D(\ve k^j) \ve g^j\|},
	\displaystyle \frac{\|\ve D^{1/2}(\ve k^j) \ve g^j\|^2}{\|\ve J(\ve k^j) \ve D(\ve k^j) \ve g^j \|^2 }\right\} &  \mbox{if } \left(\ve{k}^j-\hat{\lambda}^j_C  \ve D(\ve k^j) \ve g^j\right)_i> 0 \quad \\
	t \min_i \left\{\displaystyle \frac{\ve k_i^j}{(\ve D(\ve k^j) \ve g^j)_i} \colon  (\ve D(\ve k^j) \ve g^j)_i> 0 \right\}  & \mbox{if } \left(\ve{k}^j-\hat{\lambda}^j_C  \ve D(\ve k^j) \ve g^j\right)_i\leq 0
	\end{array}\right.
	\end{equation}}
\hskip -0.2cm with $t\in (0,1)$.
	If $\rho_j^C=\displaystyle \frac{m_j(\ve{\bar{p}}^j)}{m_j(\ve{p}^j_C)}>\beta_C\in (0,1)$ and $\rho_j>\beta\in [0.25,1)$, 
	the current trial step is accepted and the next iterate is updated as  $\ve k^{j+1}=\ve k^{j}+\ve{\bar{p}}^j$, otherwise the TR radius is reduced. In particular,
	if $\rho_j^C \leq \beta_C$, %
the unsatisfactory reduction of the quadratic model at $\ve{\bar{p}}^j$ with respect to the reduction obtained with the generalized Cauchy step highlights that
	we have to increase the effect of the regularization term by reducing
	the TR radius and computing a new reduced step; 
	this vector tend to line up with $\ve g^j$ and the new generalized Cauchy step,
	so that the sufficient reduction of the quadratic model is obtained.  
	Furthermore, when $\|\ve g ^j\| \neq 0$,
 after a successful iteration of the method, the TR radius can be further adjusted  by increasing or reducing it within a prefixed range, accordingly to a strategy proposed in~\cite{Bellavia} (see Eq.~(5.5)-(5.6)),
	as follows:\newline
	\begin{equation}\label{eqn:mu}
	\Delta_{j+1} = \max\left(\mu_{j+1}\|\ve F(\ve k^{j+1})- \ve{ y}^{\delta}\|,1.2 \frac{(1-q)\|\ve g^{j+1}\|}{\|\ve B^{j+1}\|}\right) ,
	\end{equation}
	where
	\begin{equation}
	\mu_{j+1}=\left\{ \begin{array}{cc}
	\theta \mu_j & \quad \mbox{if } q_j<q \quad \\
	\frac{\mu_j}{\eta}            & \ \ \qquad \mbox{if } q_j>1.1q\ \quad \\
	\mu_j             & \quad \mbox{otherwise}
	\end{array}  \right.
	\end{equation}
	with 
	$q\in(0, 1)$, $q_j = \displaystyle \frac{\|\ve{y}^\delta-\ve F(\ve{k}^j)-\ve{J}(\ve{k}^j)\ve{\bar{p}}^j\|}{ \|\ve{y}^\delta-\ve F(\ve{k}^j)\|}$ and $\theta, \eta\in (0,1)$.

	The regularizing affine-scaling TR method, called in the following reg-AS-TR, is summarized in Algorithm \ref{Alg:TR};	for data affected by noise, the  stopping criterion is based on the discrepancy principle~\eqref{stopcrit}. 
	
The convergence analysis of reg-AS-TR requires to prove Proposition~\ref{mono_proj}, which is analogous to Theorem 2.1 of \cite{Wang}, i.e., we need to prove that the distance between $\ve k^j$ and the exact solution $\ve k^\dag$ decreases for $j\leq \bar{j}(\delta)$. 
	To this aim, we give two essential assumptions on the local properties of the nonlinear system \eqref{nonlinearsystem}, very similar to the ones used in \cite{Hanke,Wang} to handle ill-posed problems.
	
	\begin{itemize}
		\item[{A1}.] Given an initial guess $\ve k^0 > 0$, there exist $\nu, c>0$ such that $\ve k^{\dag}\in B_{\nu}(\ve k^0)=\{\ve k\geq 0:
		\|\ve k-\ve k^0\|\leq \nu\}$
		and for all $\ve{\bar{k}},\ve{k} \in B_{2\nu}(\ve k^0)=\{\ve k\geq 0:
		\|\ve k-\ve k^0\|\leq 2\nu\}$ the following condition holds:
		\begin{equation}
		\|\ve F(\ve{\bar{k}})- \ve F(\ve k) -\ve J(\ve k)(\ve{\bar{k}}-\ve k)\| \leq c \| \ve{\bar{k}}-\ve k\|  \| \ve F(\ve{\bar{k}})- \ve F(\ve k) \|
		\end{equation}
		\item[{A2}.] $\|\ve k^0-\ve k^{\dag}\|< \min(\frac{q}{c},\nu)$ for noisy-free data ($\delta=0$) and $\|\ve k^0-\ve k^{\dag}\|< \min(\frac{q\tau-1}{c(1+\tau)},\nu)$ for noisy data ($\delta>0$) with $\tau>1/q$.
	\end{itemize}

	We highlight that, when at the first steps of the algorithm these assumptions are not verified, the initial iterations
	can enable to restrict the domain so that they hold from a certain~$j$.
	Now, we are able to state the following key proposition (for the proof see the Appendix).

	\begin{proposition}\label{mono_proj}
		Let  assume that 
		$\ve J(\ve k^j)^T \ve J(\ve k^j)+\alpha_j \ve I_n$ is positive definite, $\ve g^j\neq 0$ and 
		\begin{equation}\label{qcond}
		\|\ve{y}^\delta-\ve F(\ve{k}^j)-\ve{J}(\ve{k}^j)\ve {\ve {\bar p}}^j \|\geq q \|\ve{y}^\delta-\ve F(\ve{k}^j)\|
		\end{equation}
		for a suitable $q\in(0,1)$,  
		with $j\geq 0$ and with $j\leq \bar{j}(\delta)$ when $\delta>0$.
		Moreover, let assume that, for a suitable $\gamma_{\delta}>1$, the following condition holds for $\ve k^j>0$:
		\begin{equation}\label{gamma_cond1}
		\|\ve y^{\delta}- \ve F(\ve k^j) -\ve J(\ve k^j)(\ve{k}^\dag-\ve k^j)\| \leq \frac{q}{\gamma_\delta} \| \ve y^{\delta}- \ve F(\ve k^j) \|.
		\end{equation}
		Thus  we have
		\begin{equation}
		\|\ve{k}^\dag-\ve k^j\|^2-\|\ve{k}^\dag-\ve{k}^{j+1}\|^2> \frac{2t(\gamma_{\delta}-1)q}{\gamma_{\delta}}  \|\ve{y}^\delta-\ve F(\ve{k}^j)\| \|\ve{v}^j\| \label{primawang}
		\end{equation}
		with $\ve{v}^j= (\ve J(\ve{k}^j)\ve J(\ve{k}^j)^T+\alpha_j \ve I_N )^{-1}(\ve{y}^\delta-\ve F(\ve{k}^j))$.
	\end{proposition}

	We remark that condition~\eqref{gamma_cond1} with $j=0$ follows directly from the assumptions A1-A2 with
	$\gamma_0\geq\displaystyle \frac{q}{c\|\ve k^{\dag}-\ve k^0\|}>1$ for noise-free data. For $\delta>0$, condition~\eqref{gamma_cond1} with $j=0$
	is obtained with $\gamma_{\delta}\geq \displaystyle \frac{q\tau}{1+c\|\ve k^{\dag}-\ve k^0\|(1+\tau)}>1$,
	combining the assumptions A1-A2 with the inequality $\displaystyle \frac{\|\ve y^\delta- \ve F(\ve k^j)\|}{\delta}>\tau$ which is satisfied for $j\leq \bar{j}(\delta)$ (see \eqref{stopcrit}).
	As a consequence of Proposition~\ref{mono_proj}, $\ve k^1$ belongs to $B_{2\nu}(\ve k^0)$
	and to $B_{\nu}(\ve k^\dag)$. Therefore, for the same argument above, condition~\eqref{gamma_cond1} holds by induction for $j\geq 0$ and
	for $j\leq \bar{j}(\delta)$ when $\delta>0$; as a consequence, the sequence $\|\ve k^j -\ve k^\dag\|$ is decreasing.
		
	Based on the above proposition and the convergence results of the affine-scaling TR methods,
	the same properties of the regularizing TR method for an unconstrained nonlinear least-squares problem can be easily extended to the non-negatively constrained case.
	Under Assumptions A1-A2 on the exact solution $\ve k^\dag$, reg-AS-TR
	terminates after $\bar{j}(\delta)<\infty$ iterations, where $\delta$ is the noise level on the data,
	whereas for $\delta=0$ or $\delta\rightarrow 0$ the sequence $\{\ve k^j\}$ generated by Algorithm \ref{Alg:TR} converges to a solution of the
	original problem. 
	
	As a final remark, we point out that the ill-posedness and nonlinearity of the method, together with the local properties of reg-AS-TR imply that the effectiveness of our numerical scheme may be significantly influenced by the accuracy of both the initialization and the noise estimate. The reliability with which these two aspects are addressed is an essential requirement for the accuracy of the reconstruction results.

	\begin{algorithm}[H]
		\caption{Regularizing affine-scaling Trust-Region (reg-AS-TR) method \label{Alg:TR}} 
		\begin{algorithmic}[1]
			{\scriptsize
				\begin{tabular}{p{1.0\textwidth}}
					
					{\bf Initialize}: \ Choose $\ve{k}^{0}>0, \beta\in [0.25, 1), \gamma,\beta_C \in(0,1)$,\\
					$0<\Delta_{min}<\Delta_{max}, \ q\in(0,1),\ \mu_0=0.001$\\[.1cm]
					\ $j=0$; \\
					\While {the stopping rule is not satisfied}\vspace{.1cm} 
					\State  Evaluate $\ve B_j = \ve J(\ve k^j)^T J(\ve k^j)$ and $\ve g^j =\ve{J}^T(\ve{k}^j)(\ve F(\ve{k}^j) - {\ve{y}^\delta})$ \vspace{.1cm}
					\State  $\Delta_j =  \max\left(\mu_j\|\ve F(\ve{k}^j) - {\ve{y}^{\delta}}\|, \ 1.2\frac{(1-q)\|\ve g^j \| }{\|\ve B_j \|} \right) \in \left[ \Delta_{min}, \Delta_{max}\right] $  \vspace{.1cm}
					\Repeat \vspace{.1cm}
					\State Compute a feasible solution $\ve {\bar p}^j$ to the trust-region problem \eqref{eqn:tr}  \vspace{.1cm}
					\State Compute the Cauchy point $\ve{ p}_C^j$ 
					\State Compute $\rho_j^C = \frac{m_j(\ve {\bar p}^j)}{m_j(\ve{ p}_C^j)}$ and  $\rho^j = \frac{\Phi(\ve k^j + \ve {\bar p}^j)-\Phi(\ve k^j) } {m_j(\ve{\bar p}^j)}$  \vspace{.1cm}
					\State If $\rho_j^C\leq\beta_C $ or $\rho_j\leq\beta$ then set $\Delta_j = \gamma\Delta_j$   \vspace{.1cm}
					\Until $\rho_j^C>\beta_C $ and $\rho_j>\beta$  \vspace{.1cm}
					\State $\ve k^{j+1}=\ve k^j +\ve{\bar{p}}^j$ \vspace{.1cm}
					\State  j = j+1 \vspace{.1cm}
					\State Update $\mu_{j+1}$ as specified in \eqref{eqn:mu}\vspace{.1cm}
					\EndWhile
				\end{tabular}
			}
		\end{algorithmic}
	\end{algorithm}	
	\section{Numerical experiments} \label{sec:num-exp}
	
	The numerical validation of reg-AS-TR is performed using synthetic PET data generated by means of a digital phantom of the human brain. All simulations were realized on a workstation equipped with an Intel Xeon QuadCore E5620 processor at 2,40 GHz and 18 Gb of RAM, by implementing the method in the Matlab$^\circledR$ R2019a environment.
	
	\subsection{Simulation setting} \label{subsec:simulation}

	The starting point was the 3D Hoffman Digital Reference Object, a digital representation of the Hoffman Brain Phantom used in PET and CT imaging studies, freely available from the Imaging Research Laboratory of the Department of Radiology at the Medical Center of the University of Washington (http://depts.washington.edu/petctdro/DROhoffman\_main.html).
		
\begin{figure}[h!]
		\centering		
		\subfigure[Brain regions \label{subfig:regions}]
		{\includegraphics[width=0.26\textwidth]{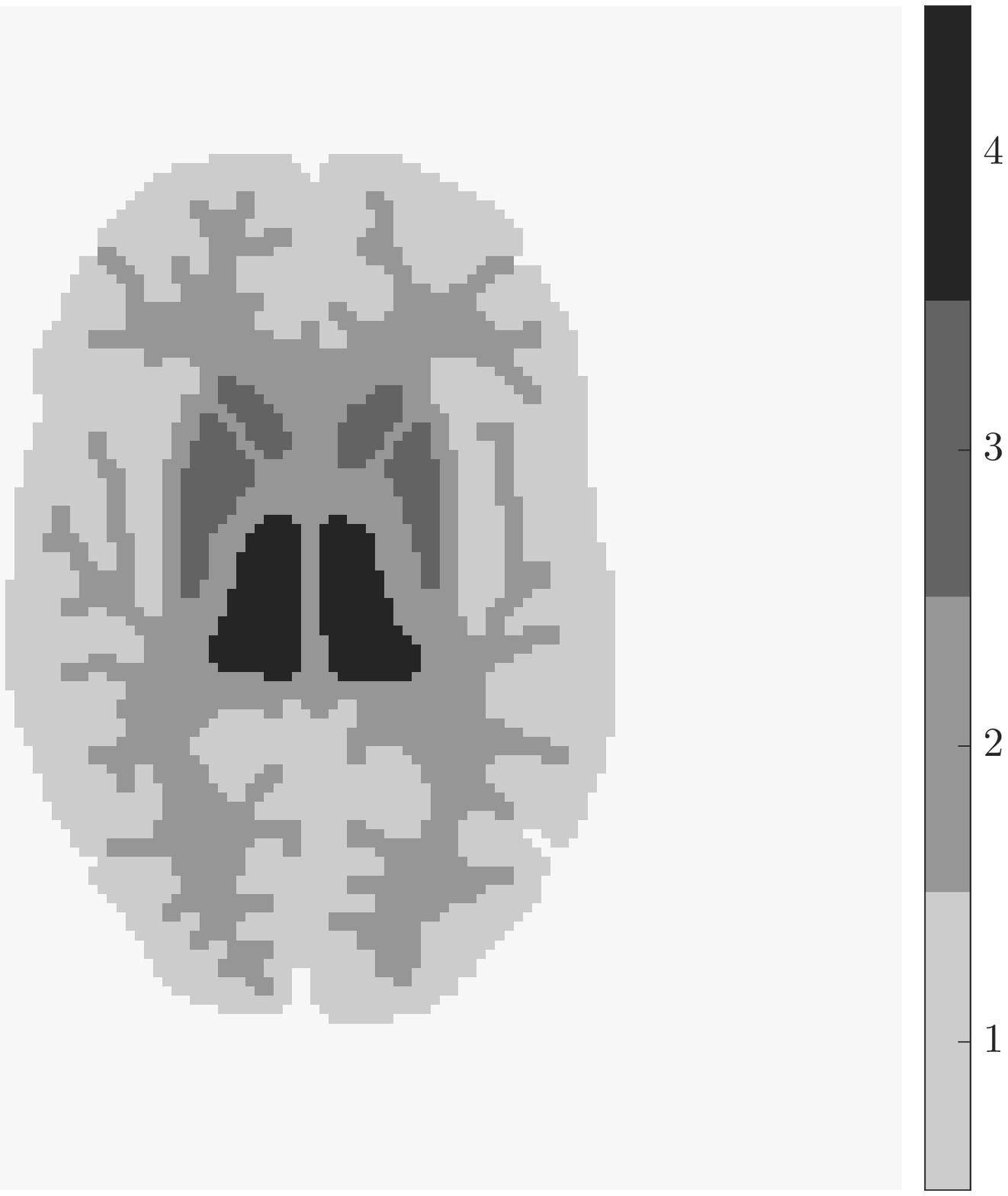}} \quad \hspace{1.4cm}
		\subfigure[Last frame PET data \label{subfig:last-frame}]
		{\includegraphics[width=0.26\textwidth]{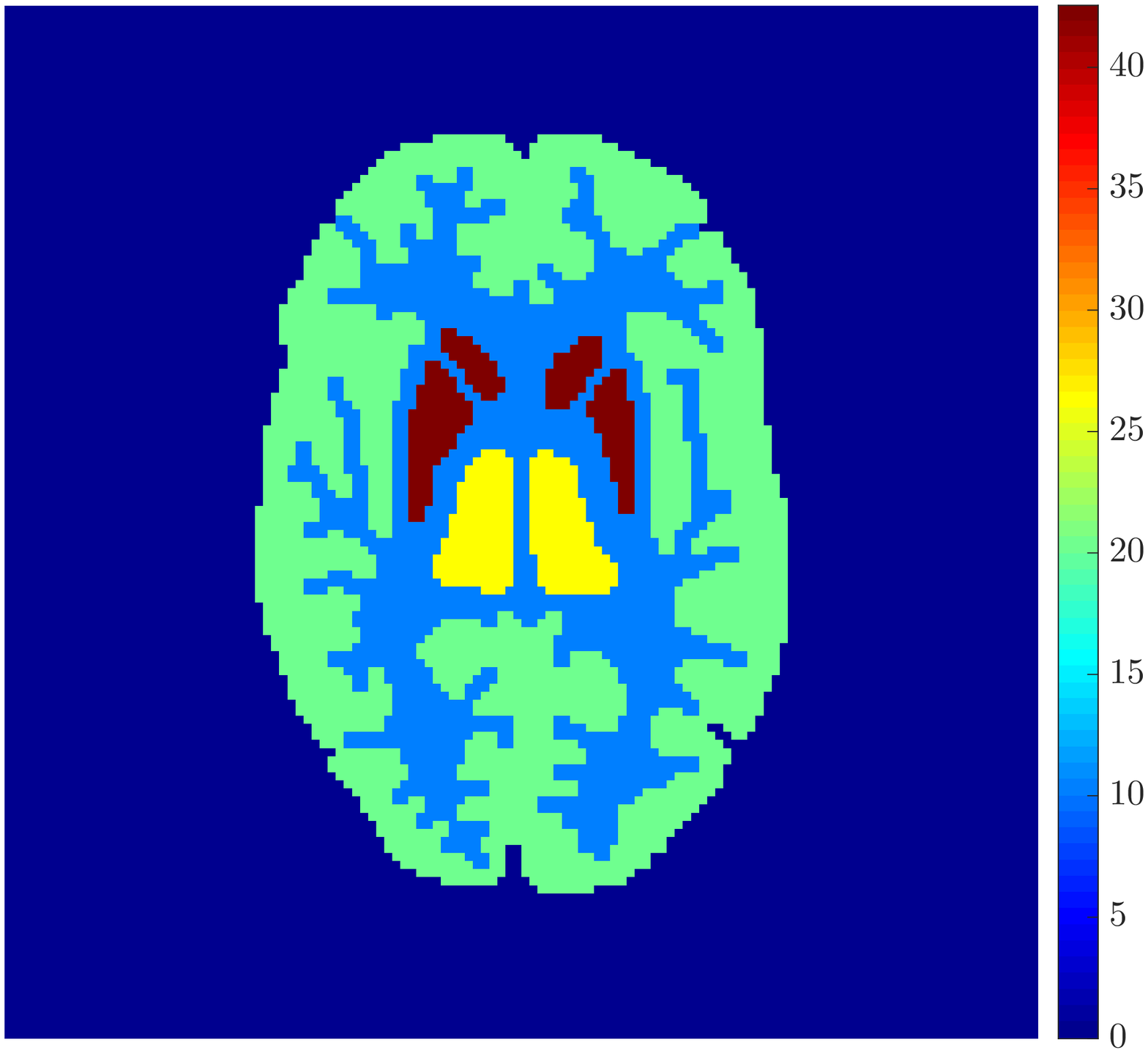}}  \\		
		\subfigure[Time concentration curves \label{subfig:TACs}] 
		{\includegraphics[width=0.28\textwidth]{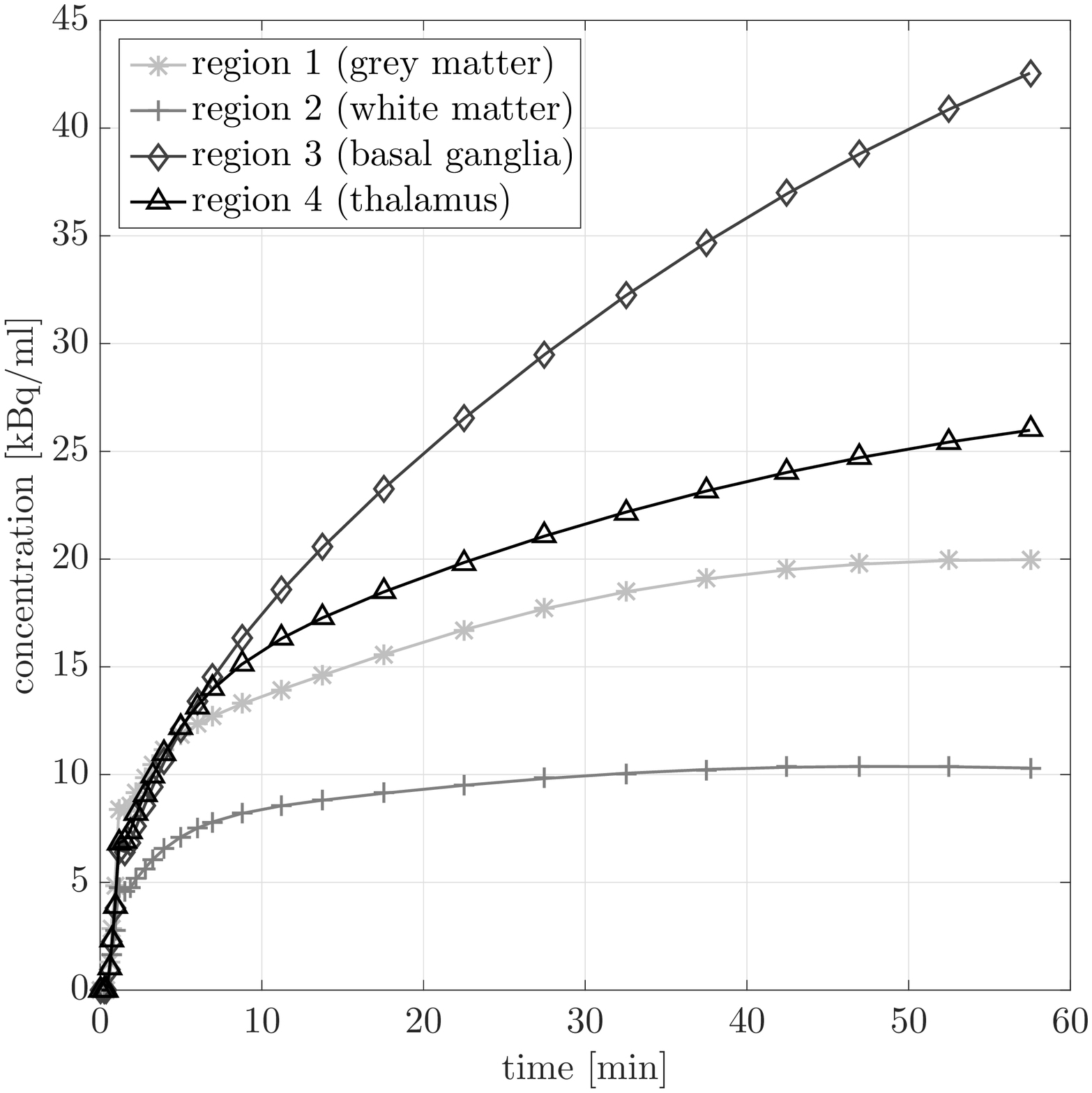}} \quad\hspace{1.3cm}
		\subfigure[IF \label{subfig:IF}]
		{\includegraphics[width=0.28\textwidth]{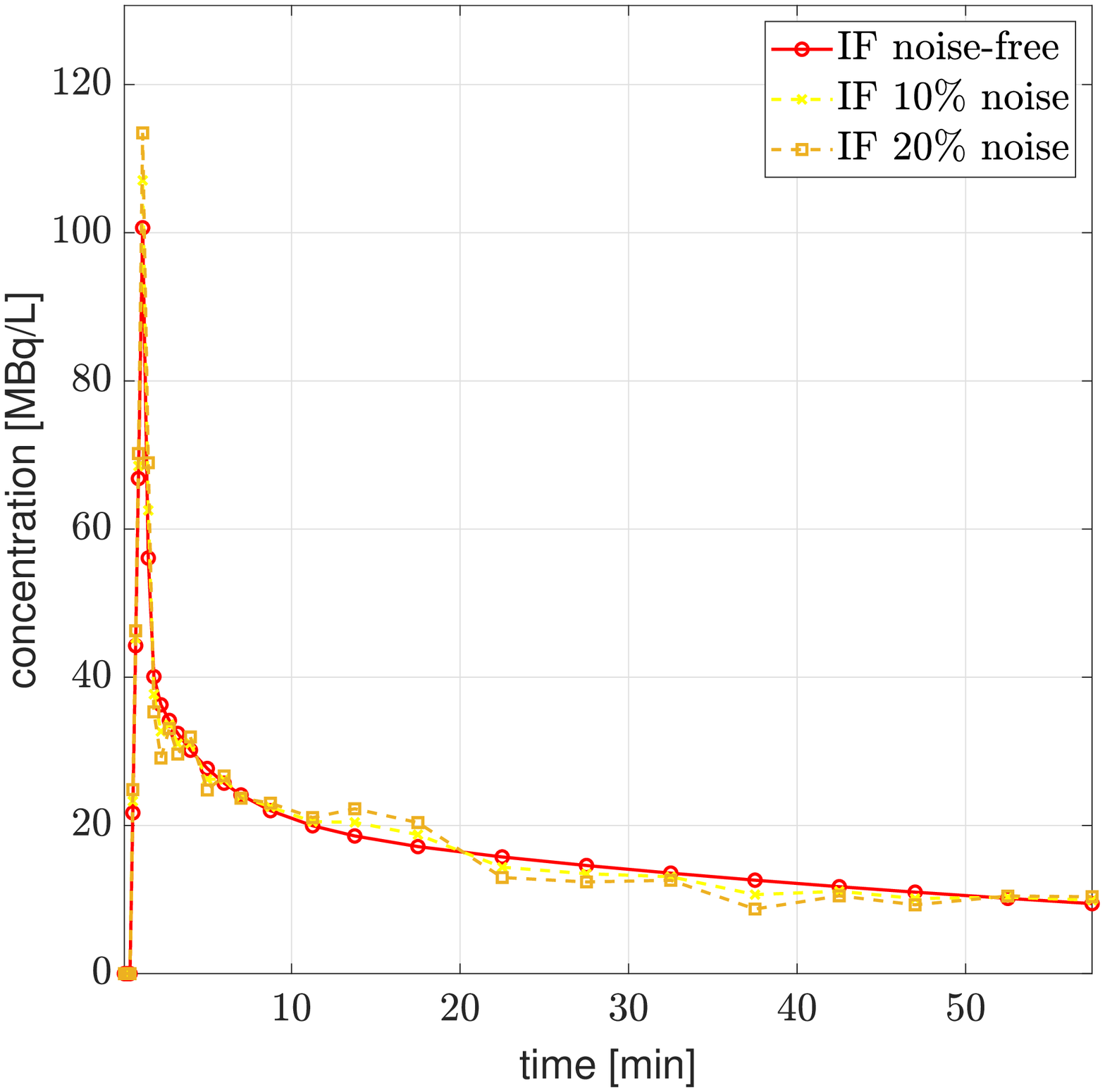}}
		\caption{Simulation layout: (a) $128 \times 128$ Hoffman brain image subdivided into four homogeneous regions; (b) last frame of the simulated PET dynamic data ; (c) time-dependent concentration curves of all brain regions; (d) simulated IF for AA of 350 MBq as noise-free, and with $10\%$, and $20\%$ of noise.}
		\label{fig:phantom}		
	\end{figure}

\begin{figure}
		\centering
		\subfigure[$k_1$ \label{subfig:k1-gt}]
		{\includegraphics[width=0.2\textwidth]{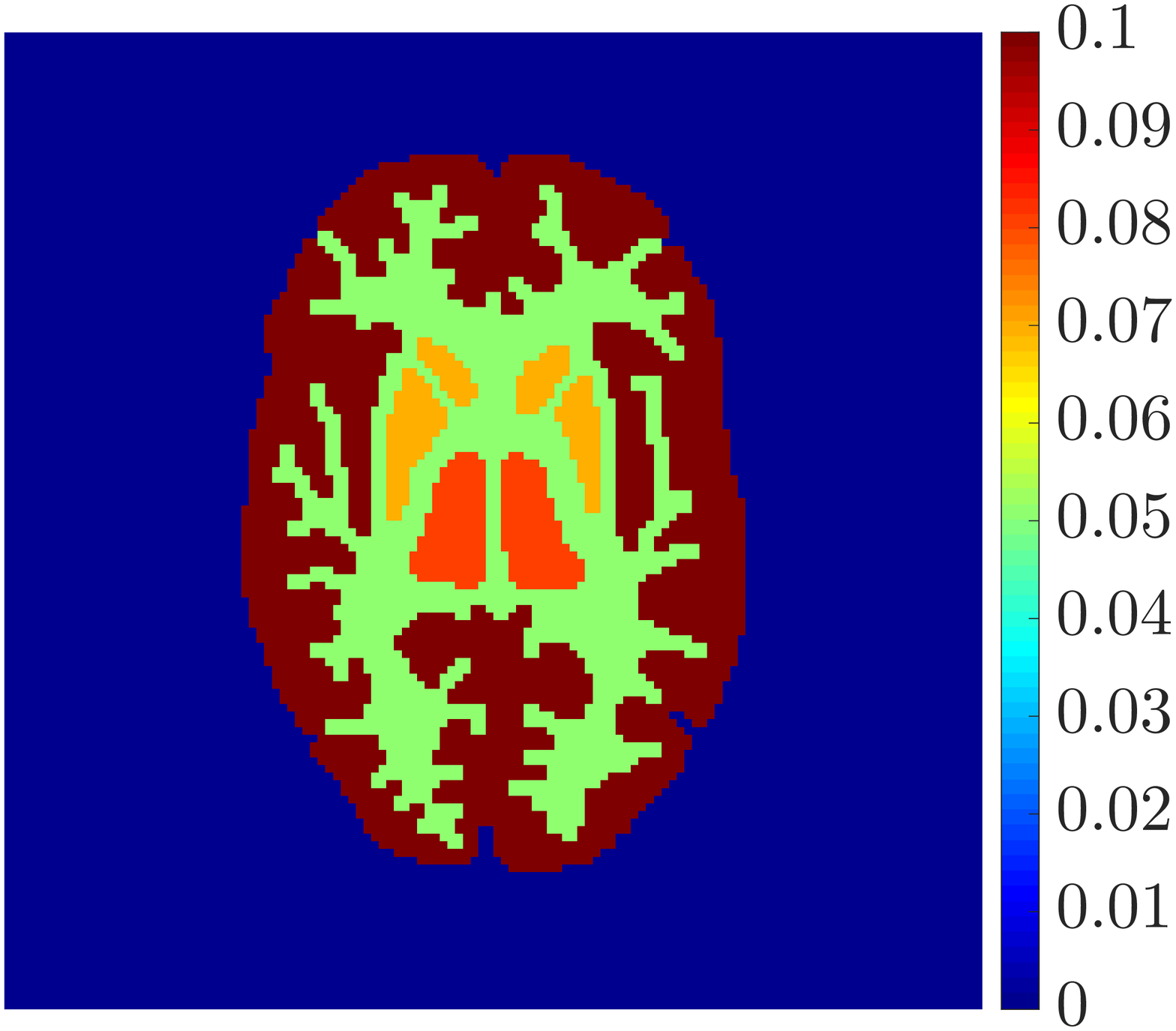}} \quad  \hspace{0.8cm}
		\subfigure[$k_2$ \label{subfig:k2-gt}]
		{\includegraphics[width=0.205\textwidth]{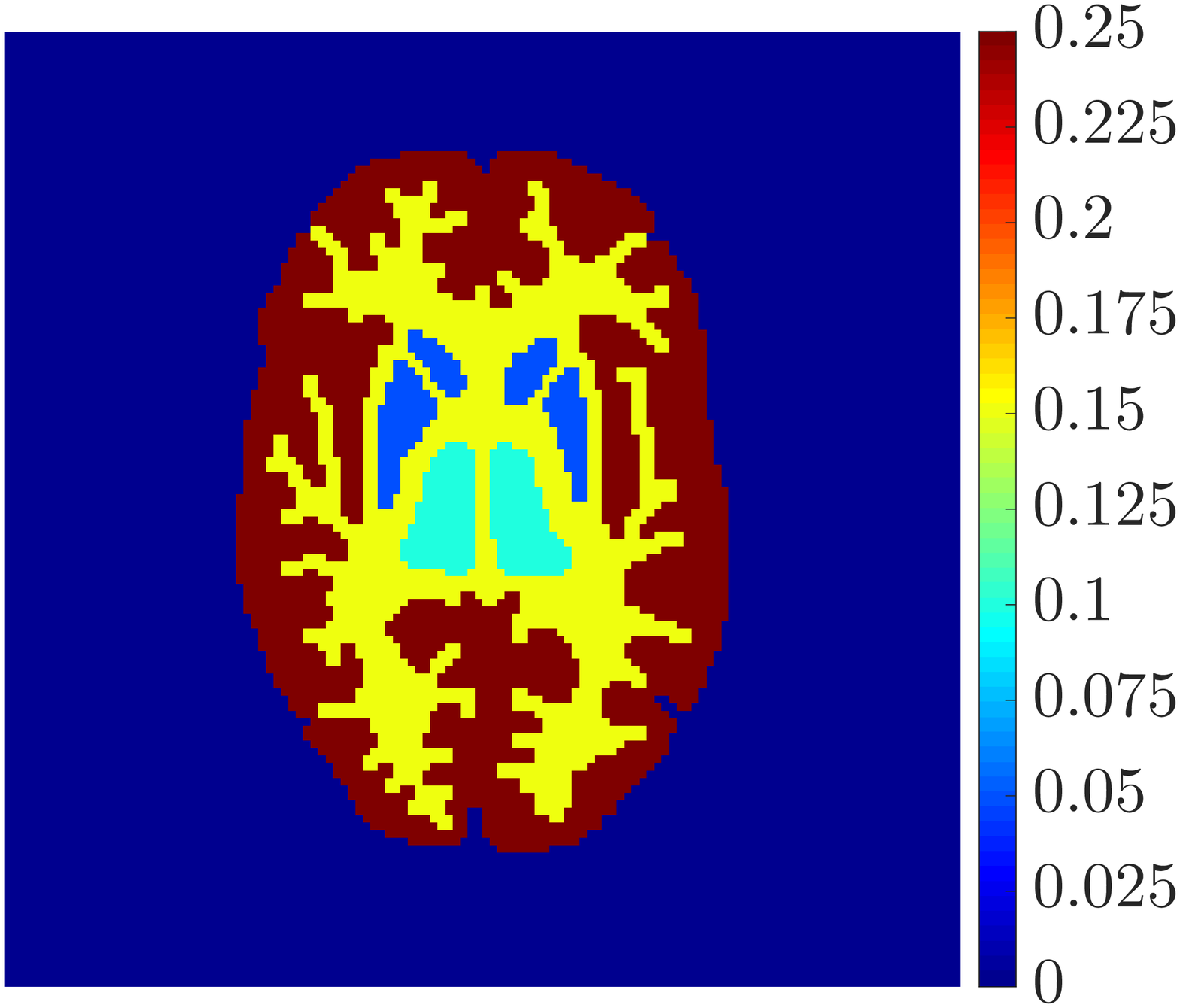}} \quad \\
		\subfigure[$k_3$ \label{subfig:k3-gt}]
		{\includegraphics[width=0.2\textwidth]{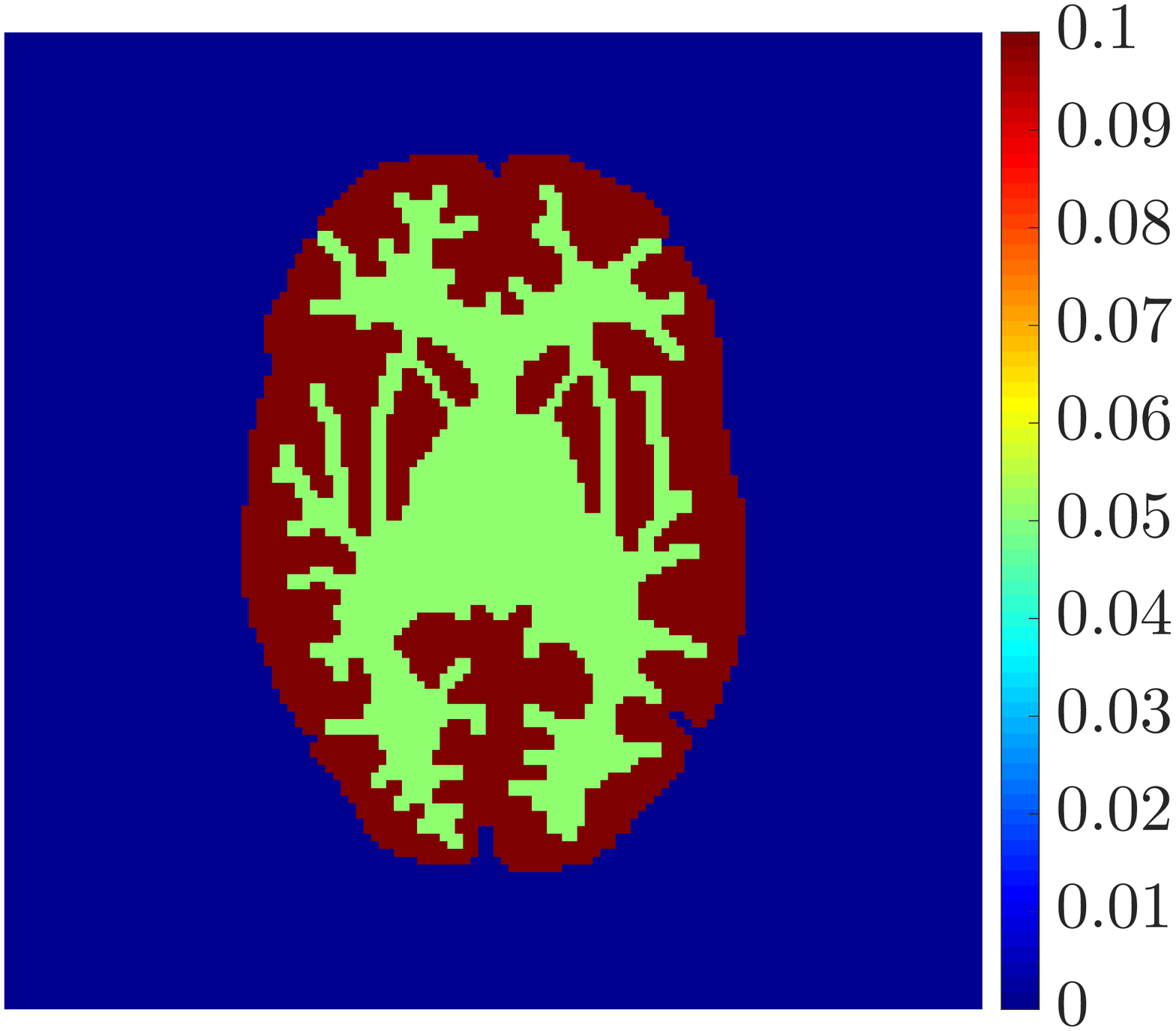}} \quad  \hspace{0.8cm}
		\subfigure[$k_4$ \label{subfig:k4-gt}]
		{\includegraphics[width=0.205\textwidth]{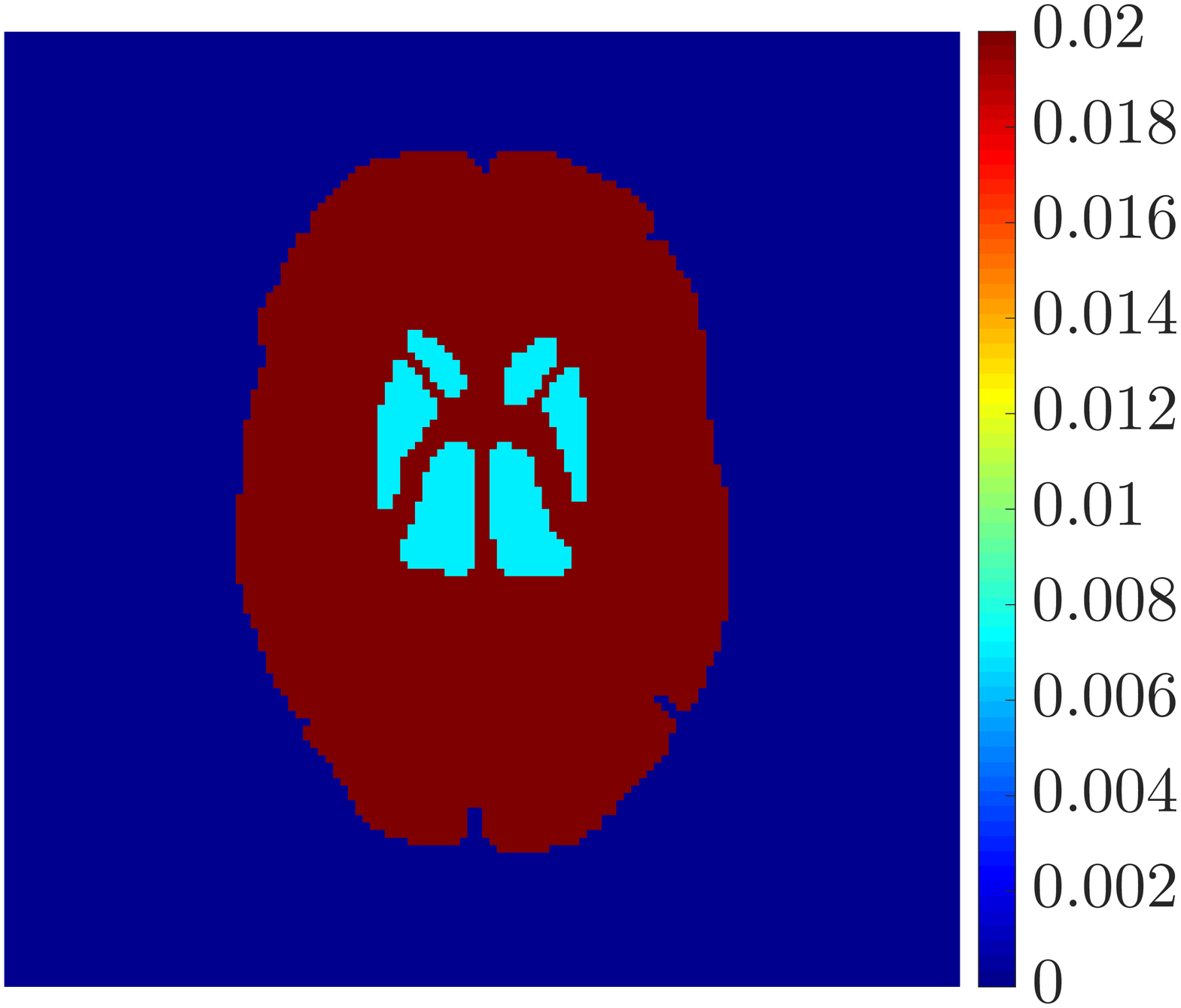}}
		\caption{Ground truth parametric images.}
		\label{fig:k-images-gt}
	\end{figure}

	The 3D Hoffman brain phantom is composed of 250 slices, covering the entire head, consisting in black/white images of size $256\times 256$. We reduced the image size to $128 \times 128$ to resemble typical PET acquisitions, preserving the shape and features of the original phantom.
	 For sake of simplicity, we selected a middle slice including eight anatomical structures that can be subdivided into the four homogeneous functional regions in Figure \ref{subfig:regions}: grey matter (region 1), white matter (region 2), basal ganglia (region 3), and thalamus (region 4).
	Then, for each region, we assigned a ground-truth set of rate constants of the two-compartment model for FDG kinetics (described in Section 2) and a specific blood volume fraction $V$. The numerical values of such parameters, as reported in Table~\ref{table:1}, have been chosen in order to reproduce a realistic framework for the FDG uptake of a human brain \cite{Guo,Huang,Wang-Qi2,Wang-Qi1}. The ground-truth parametric images are shown in \figurename~\ref{fig:k-images-gt}.
	
	In order to model the IF we implemented the following procedure \cite{Vriens}. We considered a mathematical function (see Eq. (2) in \cite{Vriens}) consisting of an increasing linear component followed by a tri-exponential decay; we fitted the free parameters of this function against measurements for $80$ subjects; we selected the median estimated parameters computed over all $80$ subjects (see Table 2 in \cite{Vriens}), a median initial distribution volume ($12.7$ L corresponding to $0.1683$ L/kg body weight), and an Administered Activity (AA) of $350$ MBq (typical of human PET acquisitions).
	The resulting simulated IF is shown in \figurename~\ref{subfig:IF}.
		
	The dynamic PET data were generated by solving the compartmental forward problem for each pixel of the processed Hoffman brain image. In particular, the two-compartment concentrations were evaluated by means of the integral equation \eqref{eq:sol_2C} with the ground truth values of the compartmental parameters and the simulated IF, at 28 time frames (6 $\times$ 10 sec, 3 $\times$ 20 sec, 3 $\times$ 30 sec, 4 $\times$ 60 sec, 3 $\times$ 150 sec, 9 $\times$ 300 sec) with a time sampling typical of standard PET experiments, for a total time interval of 60 minutes.
	Then, the measurement equation \eqref{eq:inv_2C} was computed to create the time concentration curves characteristic for each brain region (\figurename~\ref{subfig:TACs}). The last frame of the obtained dynamic PET images is reported in \figurename~\ref{subfig:last-frame}.
	
	Once the noise-free dynamic PET images were obtained, we projected the images into the sinogram space by means of the Radon transform, and we added Poisson noise to the projected data through the Matlab function \emph{poissrnd}. In this way, we created ten independent identically-distributed noisy data.
	In addition to the noise-free IF case, we considered two  
further instances where the IF was perturbed by two Gaussian noise levels:
$C_b^{c} = C_b(t_i)(1 + c \cdot r)$, for time points $t_i$, $i=1,\dots, N$, where $r$ is randomly generated from a standard normal distribution of mean~$0$ and standard deviation~$1$, and $c =  0.10, \ 0.20$ (\figurename~\ref{subfig:IF}). 

	\begin{table}	
		\begin{center}
			\caption{Ground truth values of the kinetic parameters $k_1$, $k_2$, $k_3$, $k_4$, (min$^{-1}$) and the blood volume fraction $V$, for each one of the four homogeneous regions.}\label{table:1}
			\begin{tabular}{cccccr}
				\hline
				& $k_1$ & $k_2$ & $k_3$ & $k_4$ & $V$ \\
				\hline
				region 1 &  0.100  &  0.250  &  0.100  &  0.020  &  0.050	\\
				region 2 &  0.050  &  0.150  &  0.050  &  0.020  &  0.030 	\\
				region 3 &  0.070  &  0.050  &  0.100  &  0.007  &  0.040	\\
				region 4 &  0.080  &  0.100  &  0.050  &  0.007  &  0.050	\\
				\hline
			\end{tabular}
		\end{center}
	\end{table}

	\subsection{Setup of the algorithms} \label{subsec:results}
	
	The parametric reconstruction by means of reg-AS-TR was performed as follows.
	
			In order to remove blurring artifacts from the images of each dataset, we applied a well-known deblurring technique based on the minimization of the Kullback-Leibler divergence with a smooth total-variation regularization term referred to as \emph{hypersurface potential}; this minimization is performed by means of the Scaled Gradient Projection (SGP) method proposed in \cite{Bonettini-etal-2009} (see also \cite{Crisci_etal2019}),
		starting from the inverse Radon transform of the noisy sinogram data. The deblurring procedure exploits the parallel toolbox of Matlab enabling the use of \texttt{GPUarray} and it requires about $7$ minutes overall.\\
	 
The stopping criterion of reg-AS-TR is the following:
		\begin{equation}\label{stop}
		   \epsilon_{j} < \tau_1 \quad \text{or} \quad
		    \left( \epsilon_{j}   < \tau_2 \quad \text{and} \quad  \left|1 - \frac{\epsilon_{j-1}}{\epsilon_j}\right|  <10^{-2} \right)
			\end{equation}
			where $\epsilon_{j}$ = $\|\ve y^{\delta} - \ve F(\ve k^j) \|$, 
			$\tau_1$ is the sample standard deviation computed at the current pixel and~$\tau_2$ is a multiple of~$\tau_1$, which changes accordingly when the procedure switches between boundary ($\tau_2 =10 \tau_1$ ) and inner pixels ($\tau_2 =3 \tau_1$ ) of a region.
	   In addition, if condition~(\ref{stop}) is not satisfied, the execution terminates when stagnation  
	   or the maximum number of iterations are reached.
	
	 The stopping rule implemented allows to diversify the initialization procedure of reg-AS-TR. 
	In general, the initial vector is randomly chosen in an interval determined by a priori knowledges on the physiology. 
	However, when the current pixel is strictly inside a functional region and some neighboring pixels have been already successfully processed, the initialization value is the mean over the values obtained on these neighboring pixels.

The reconstruction accuracy of reg-AS-TR has been assessed by comparison with both the ground truth and the parametric images provided by a recently introduced regularized Gauss-Newton method (reg-GN) \cite{Scussolini}. For sake of comparison, the setup of reg-GN is coherent with what is done in that paper, i.e.:
	
	\begin{itemize}
		\item \emph{Deblurring.} The noise on the PET datasets was reduced by applying a Gaussian smoothing filter (mean~$0$, standard deviation~$1$, window 3$\times$3) directly to the noisy PET images.
		\item \emph{Initialization.} The starting point of the kinetic parameters was chosen randomly in intervals determined by knowledge on the physiology.
		\item \emph{Stopping criterion.} The iterative scheme is stopped when the relative error between the experimental dynamic concentration and the model-predicted one is less than an appropriate threshold, or the maximum number of iterations is reached.
	\end{itemize}
	
\subsection{Results}

\begin{figure}
		\centering
		\includegraphics[width=0.16\textwidth]{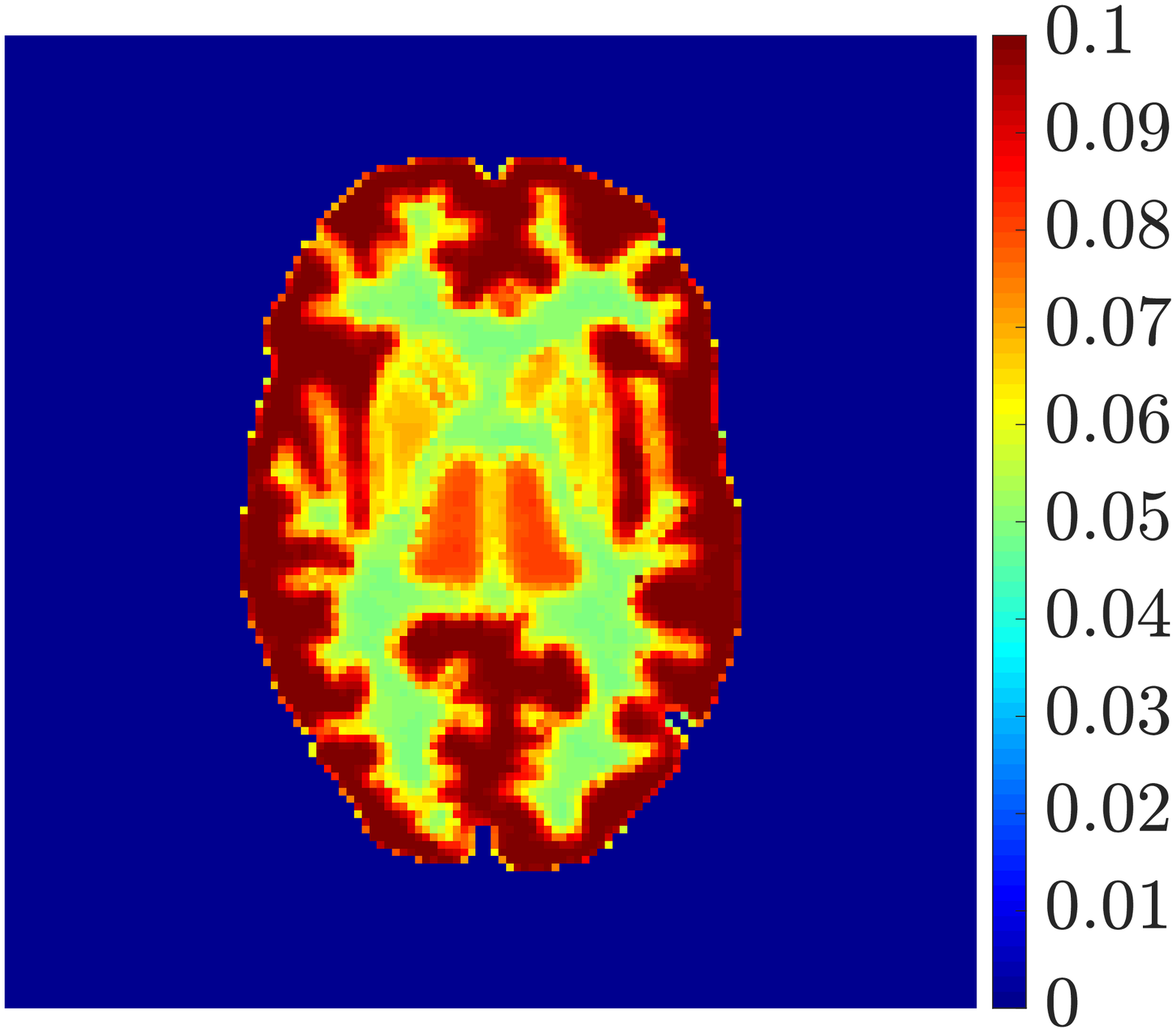}\quad \hspace{0.55cm}
		\includegraphics[width=0.16\textwidth]{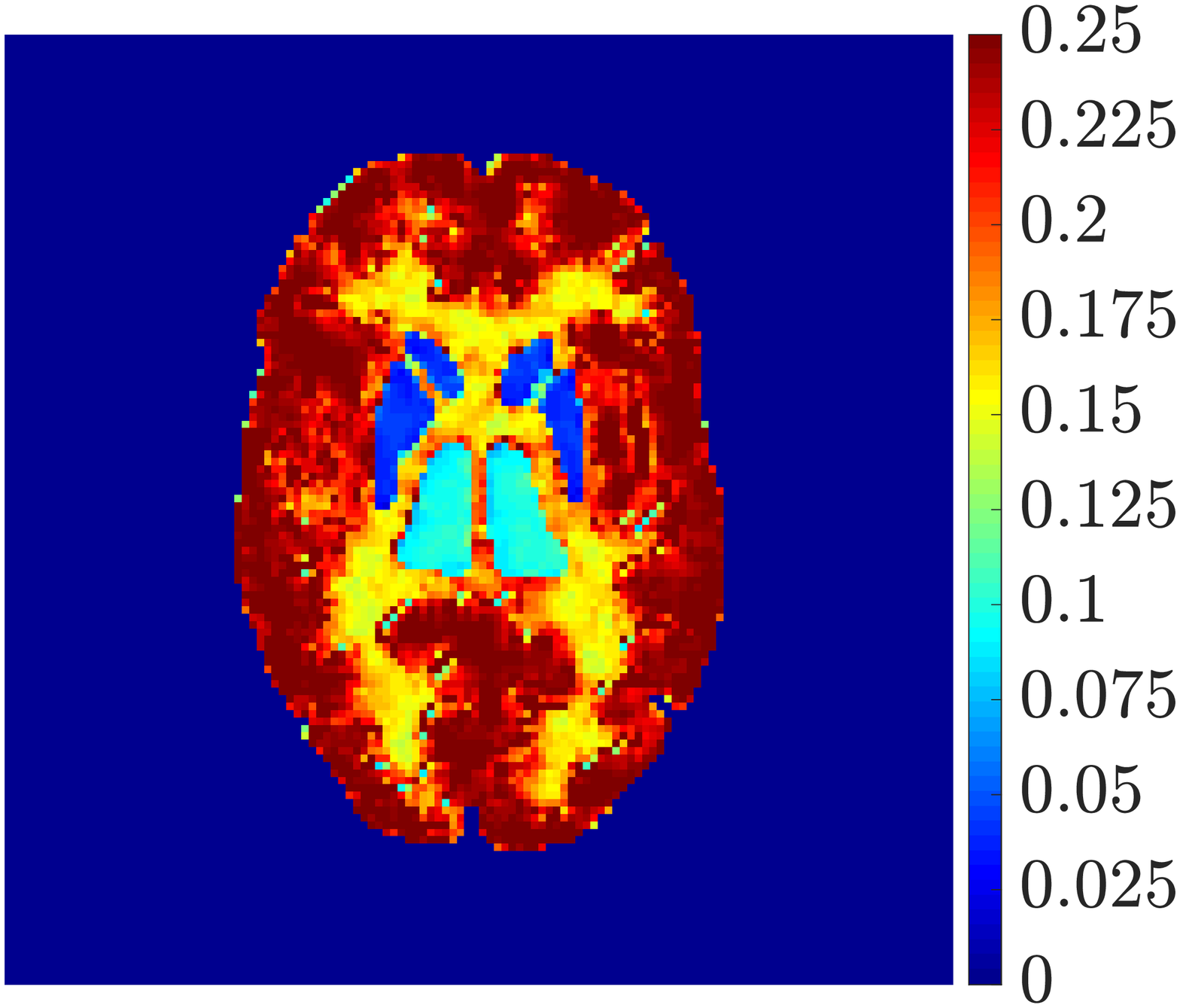}\quad  \hspace{0.55cm}
		\includegraphics[width=0.16\textwidth]{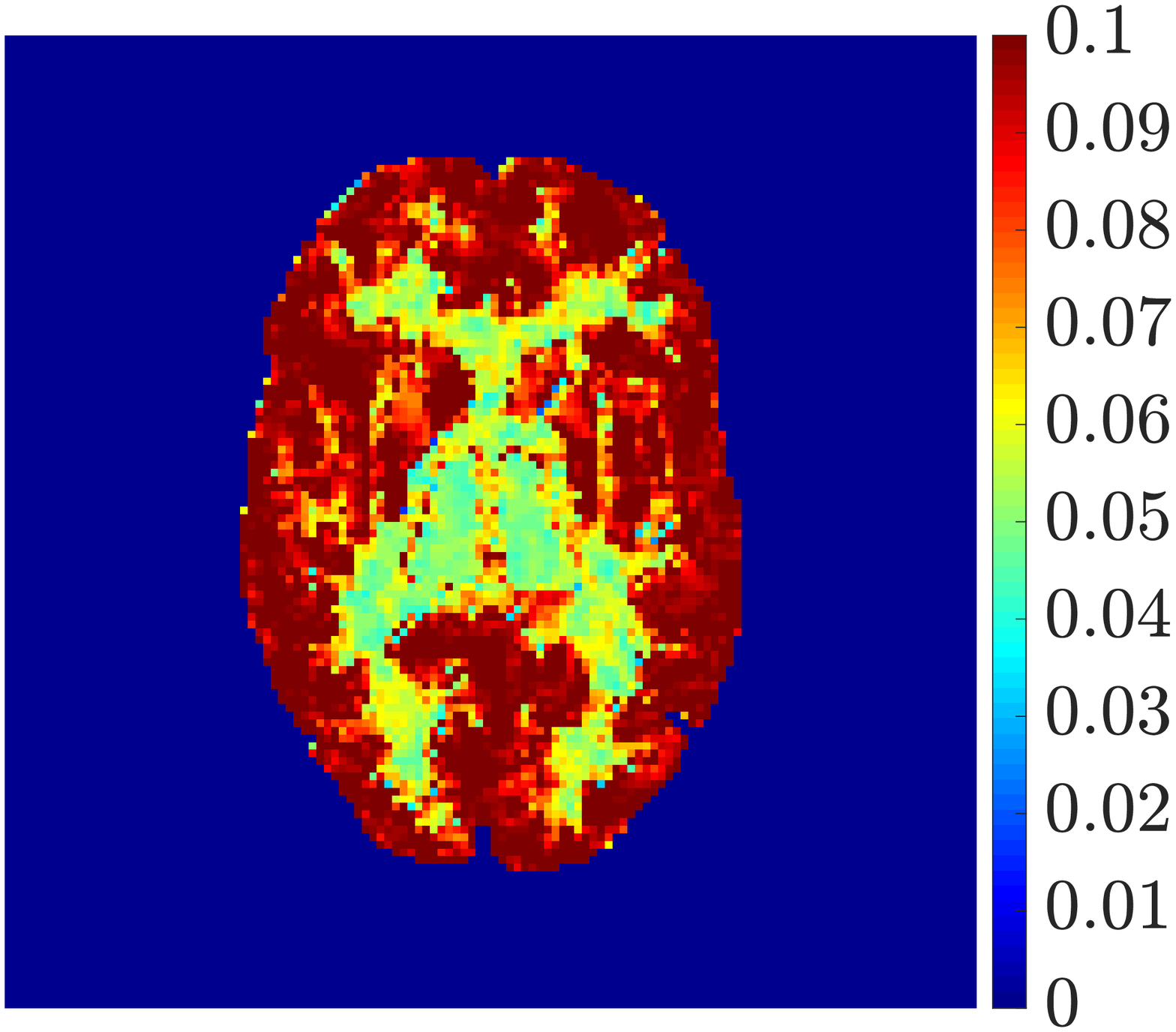}\quad  \hspace{0.55cm}
		\includegraphics[width=0.16\textwidth]{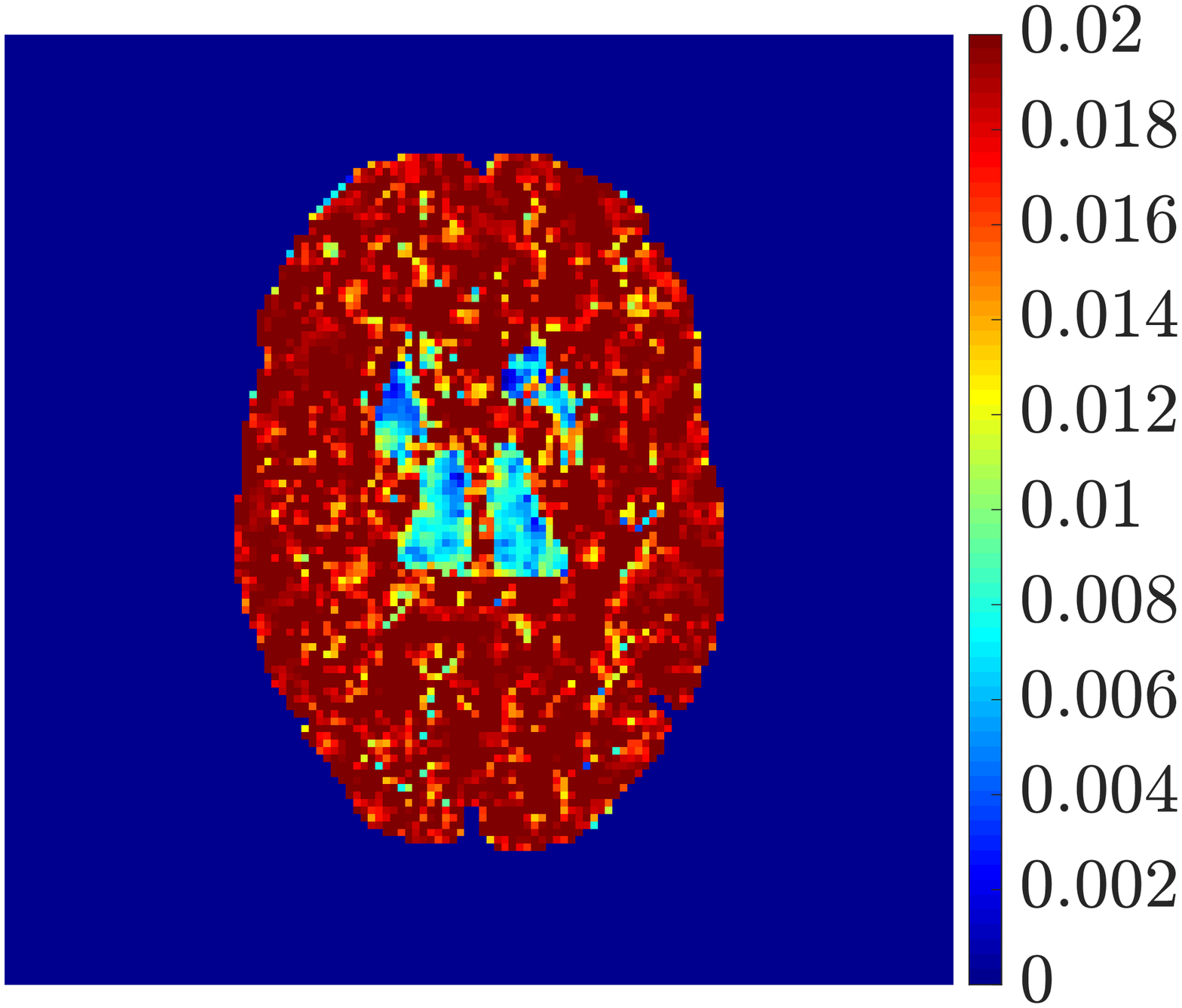}\\

		\includegraphics[width=0.16\textwidth]{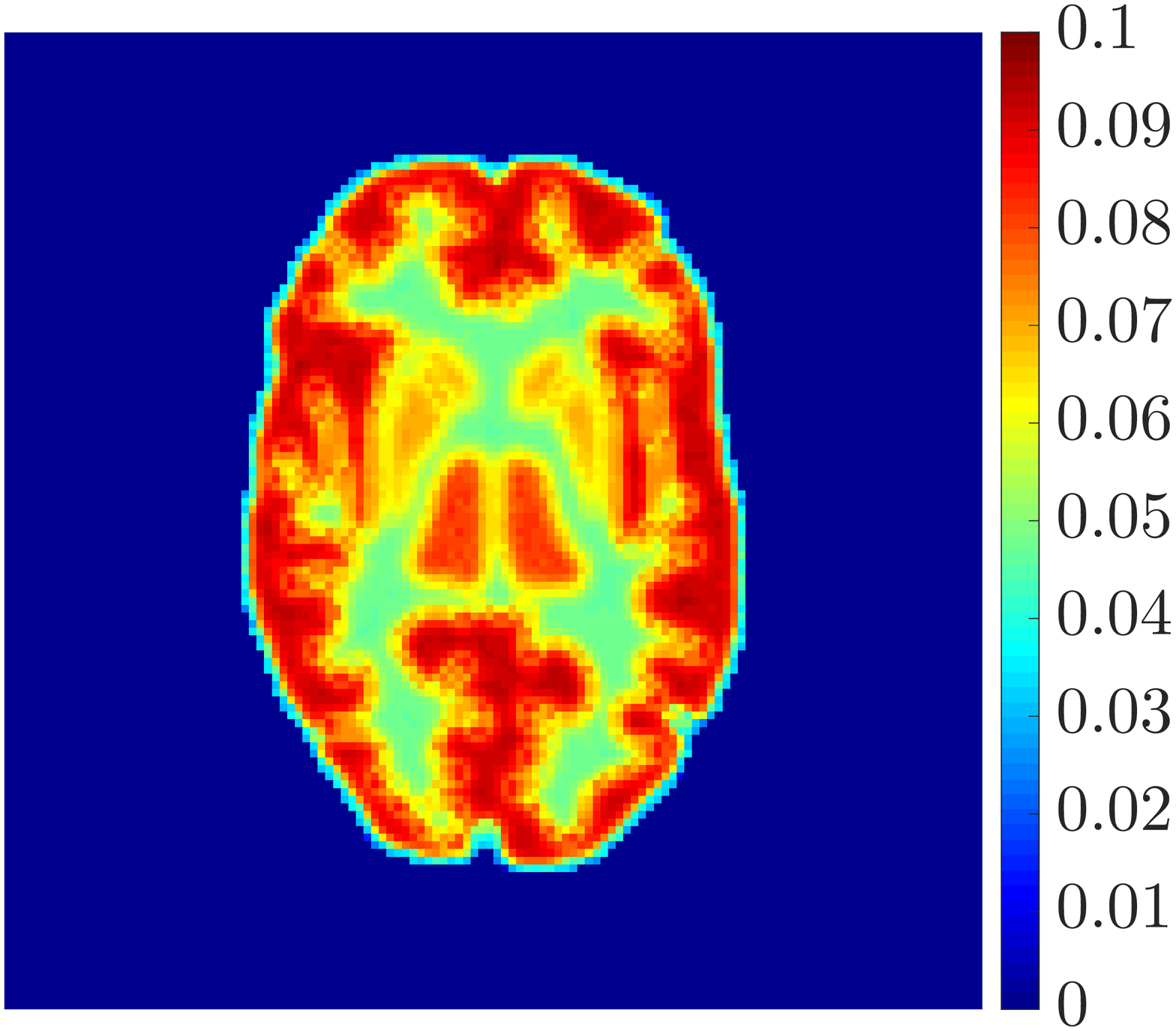}\quad  \hspace{0.55cm}
		\includegraphics[width=0.16\textwidth]{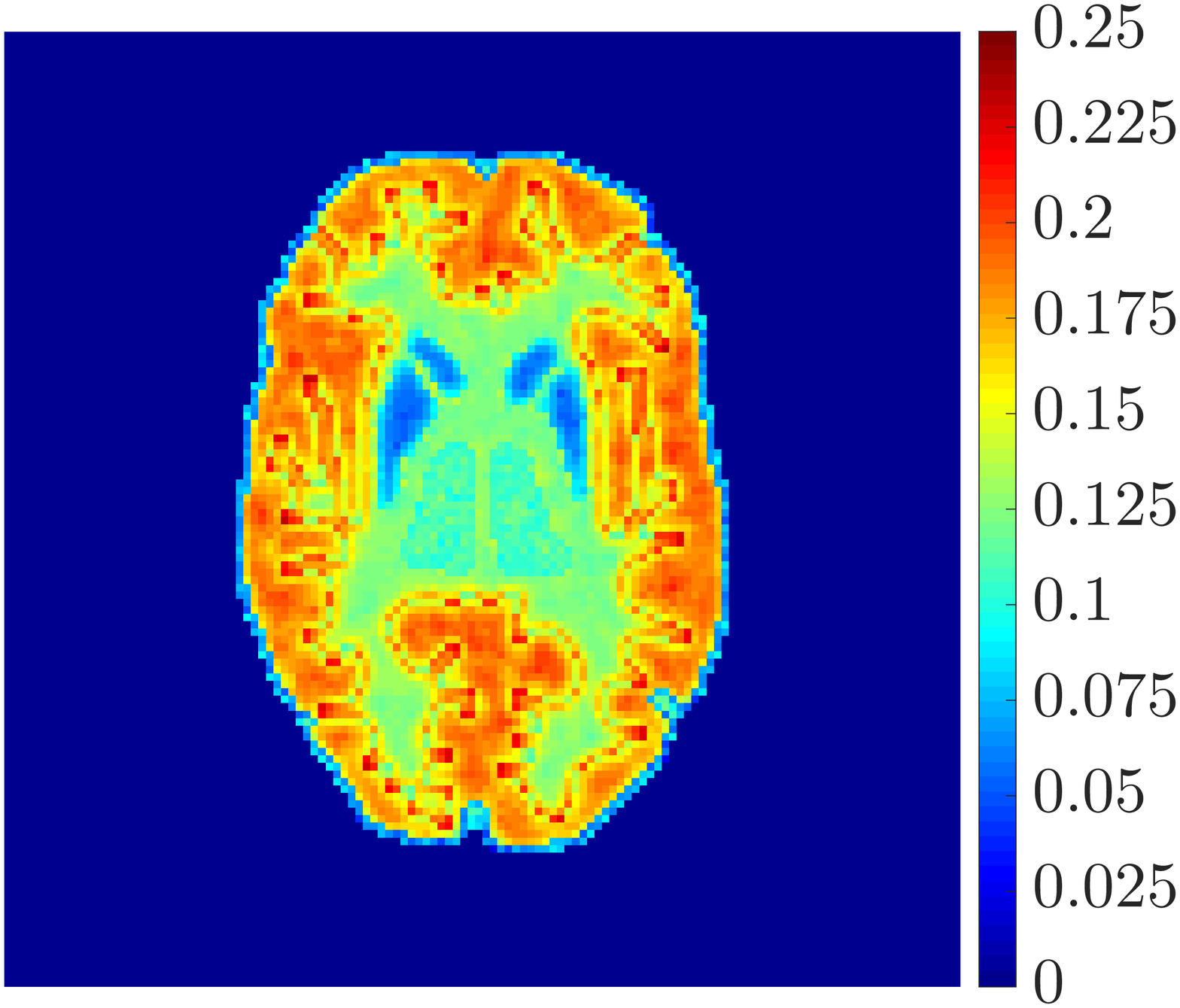}\quad  \hspace{0.55cm}
		\includegraphics[width=0.16\textwidth]{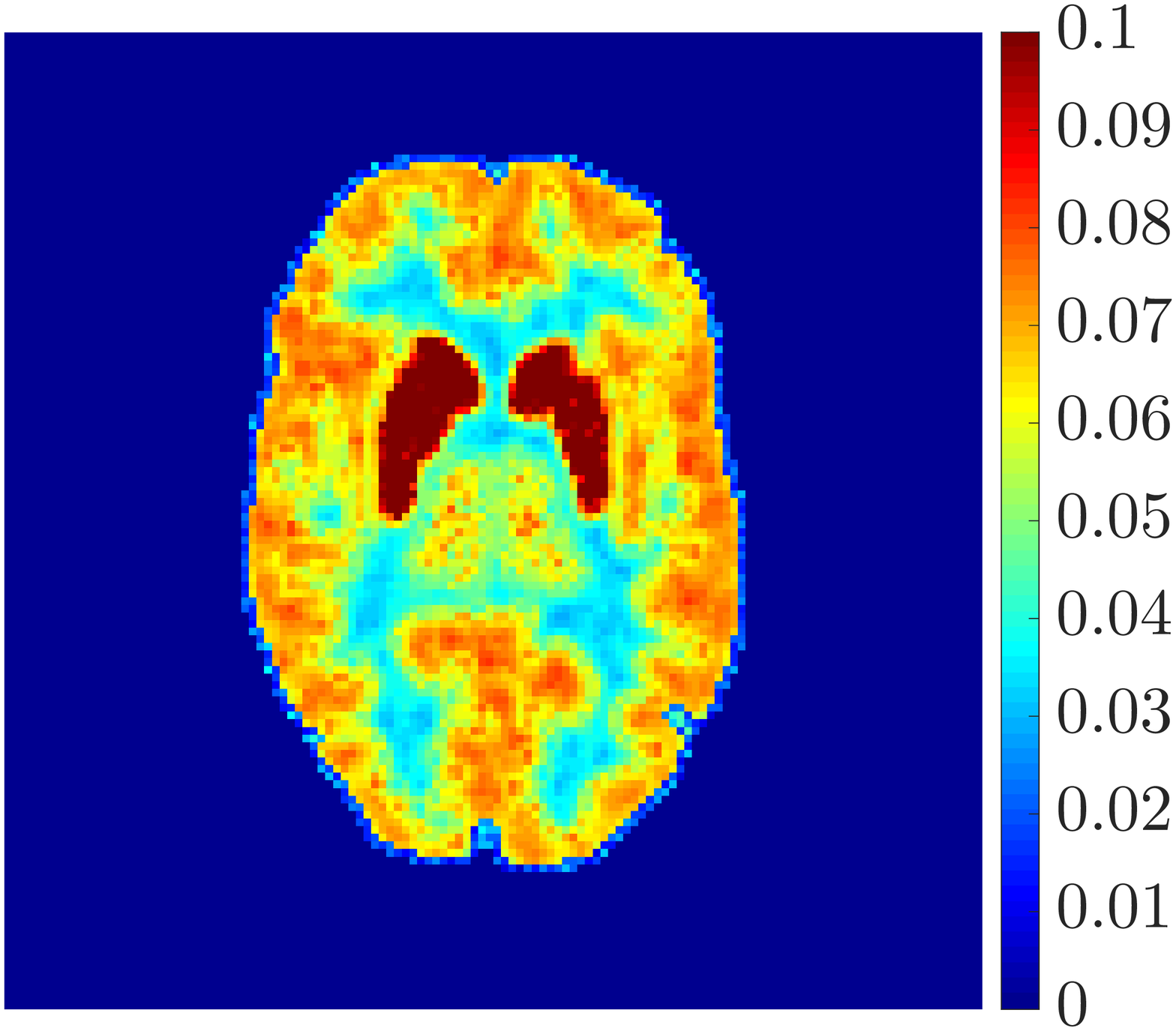}\quad  \hspace{0.55cm}
		\includegraphics[width=0.16\textwidth]{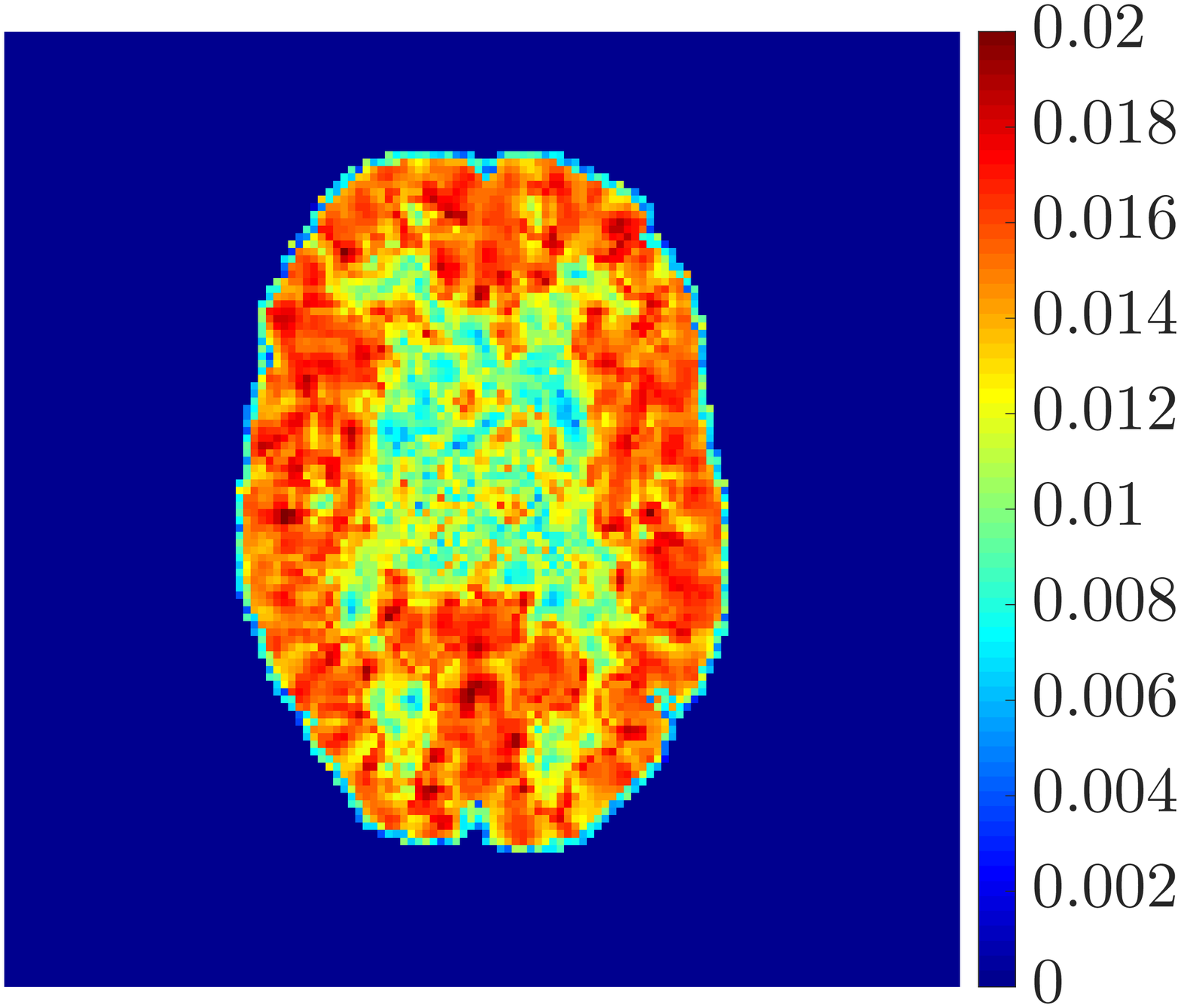}\\

		\includegraphics[width=0.16\textwidth]{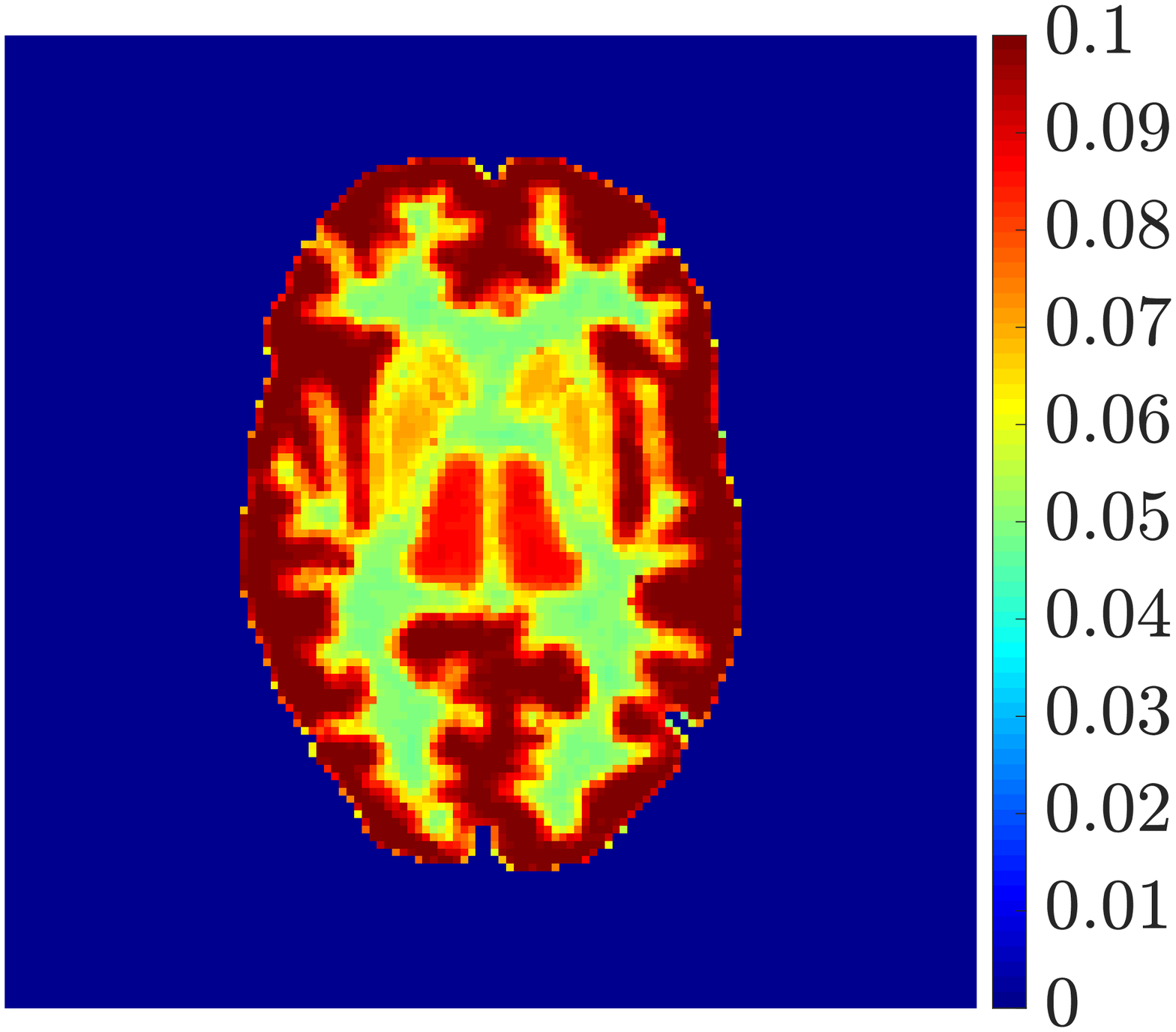}\quad  \hspace{0.55cm}
		\includegraphics[width=0.16\textwidth]{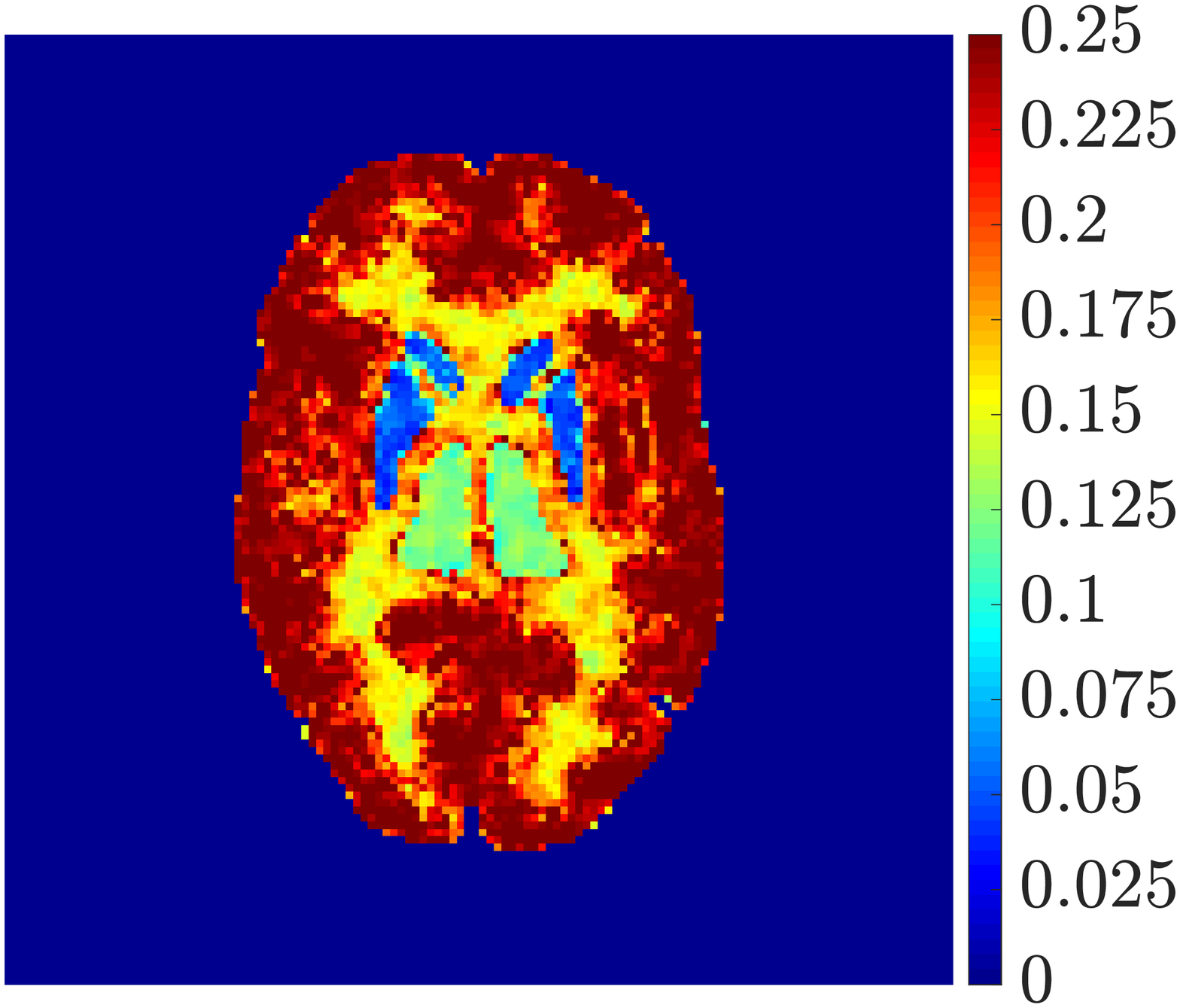}\quad  \hspace{0.55cm}
		\includegraphics[width=0.16\textwidth]{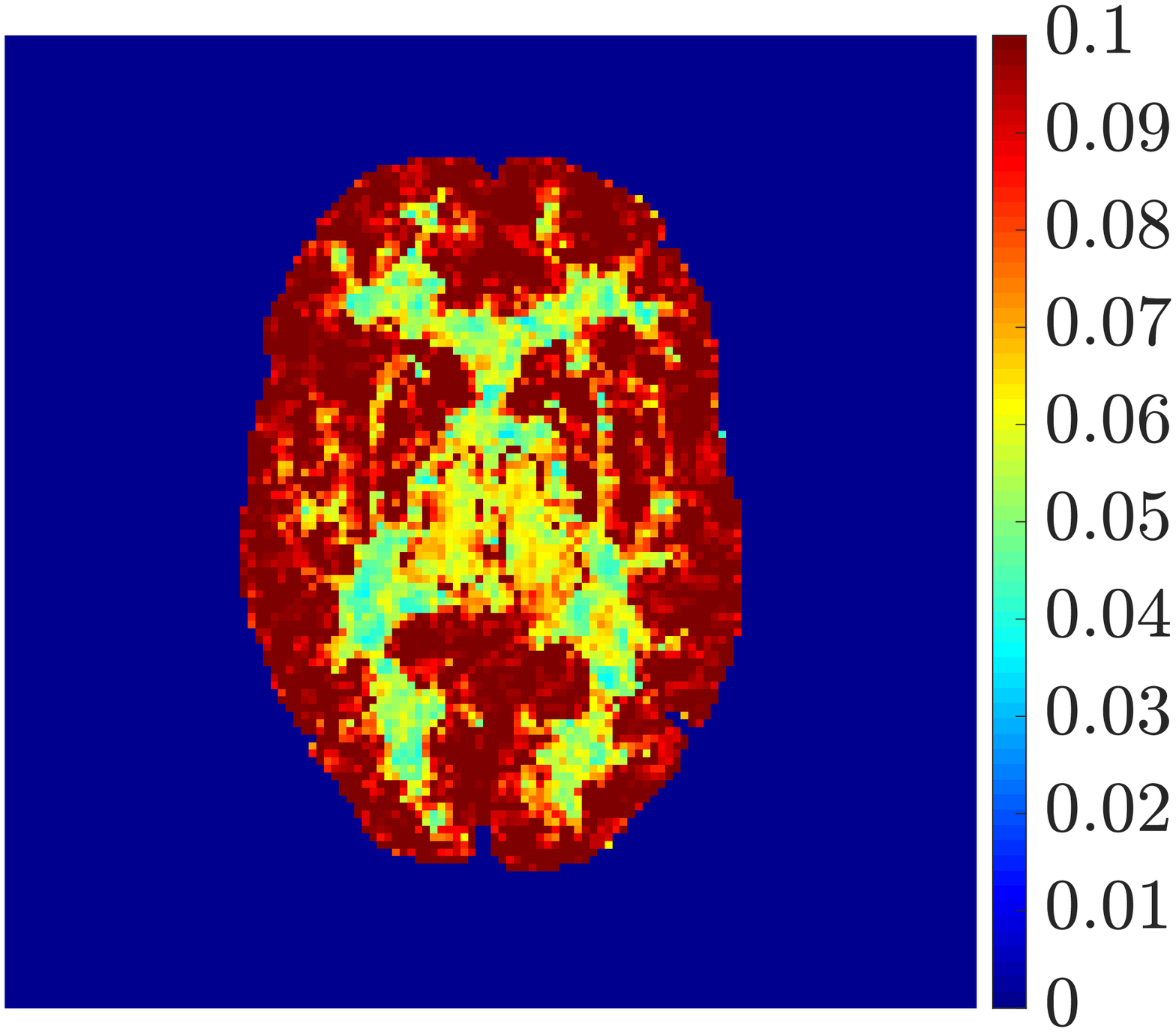}\quad   \hspace{0.55cm}
		\includegraphics[width=0.16\textwidth]{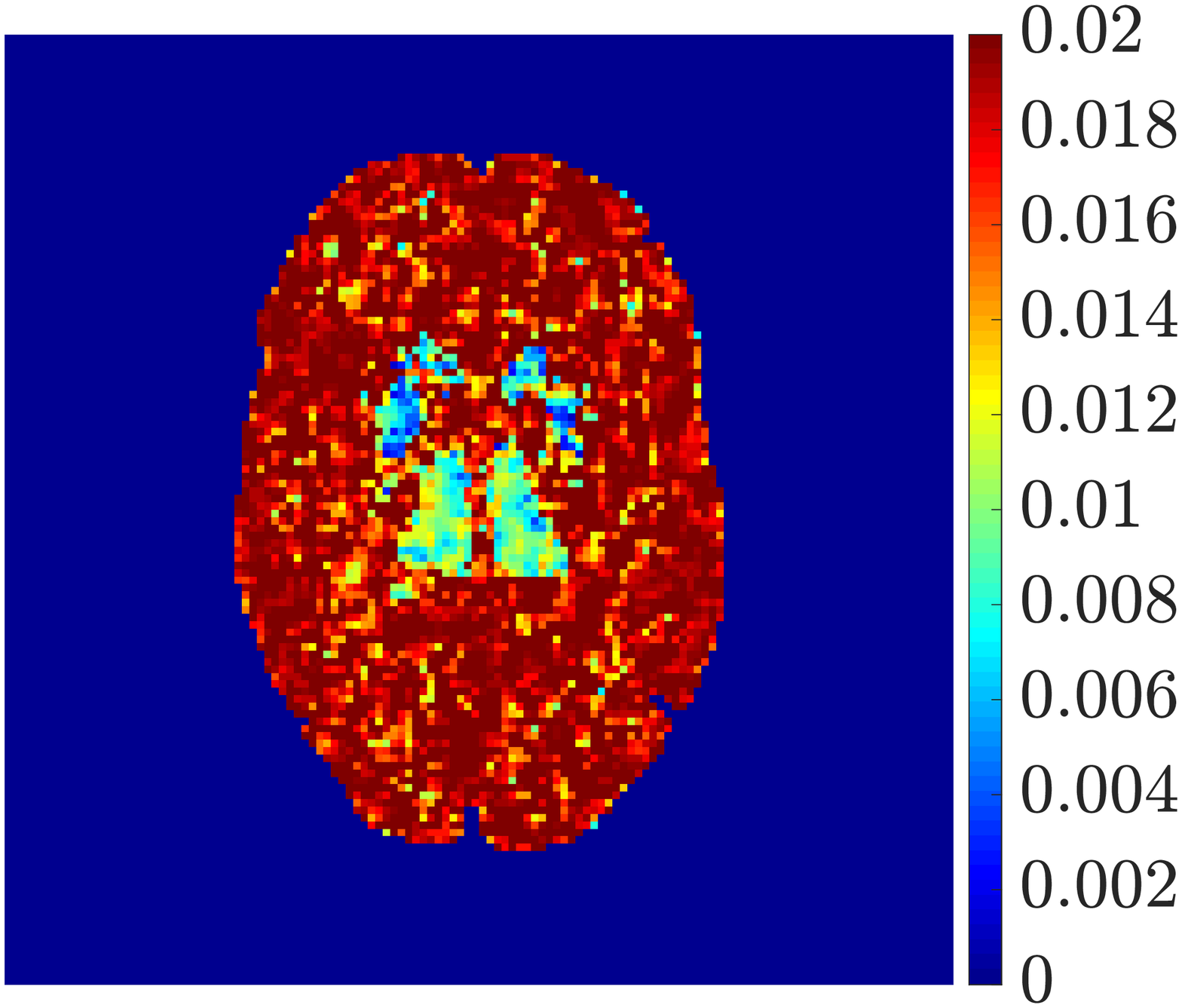}
		\caption{From left to right: mean parametric images corresponding to $k_1$, $k_2$, $k_3$, $k_4$, obtained by using reg-AS-TR (first row), reg-GN (second row), \emph{lsqcurvefit} (third row). Case noise-free IF.}
		\label{fig:k-rec-noise-free}
	\end{figure}

\begin{figure}
		\centering
	\includegraphics[width=0.16\textwidth]{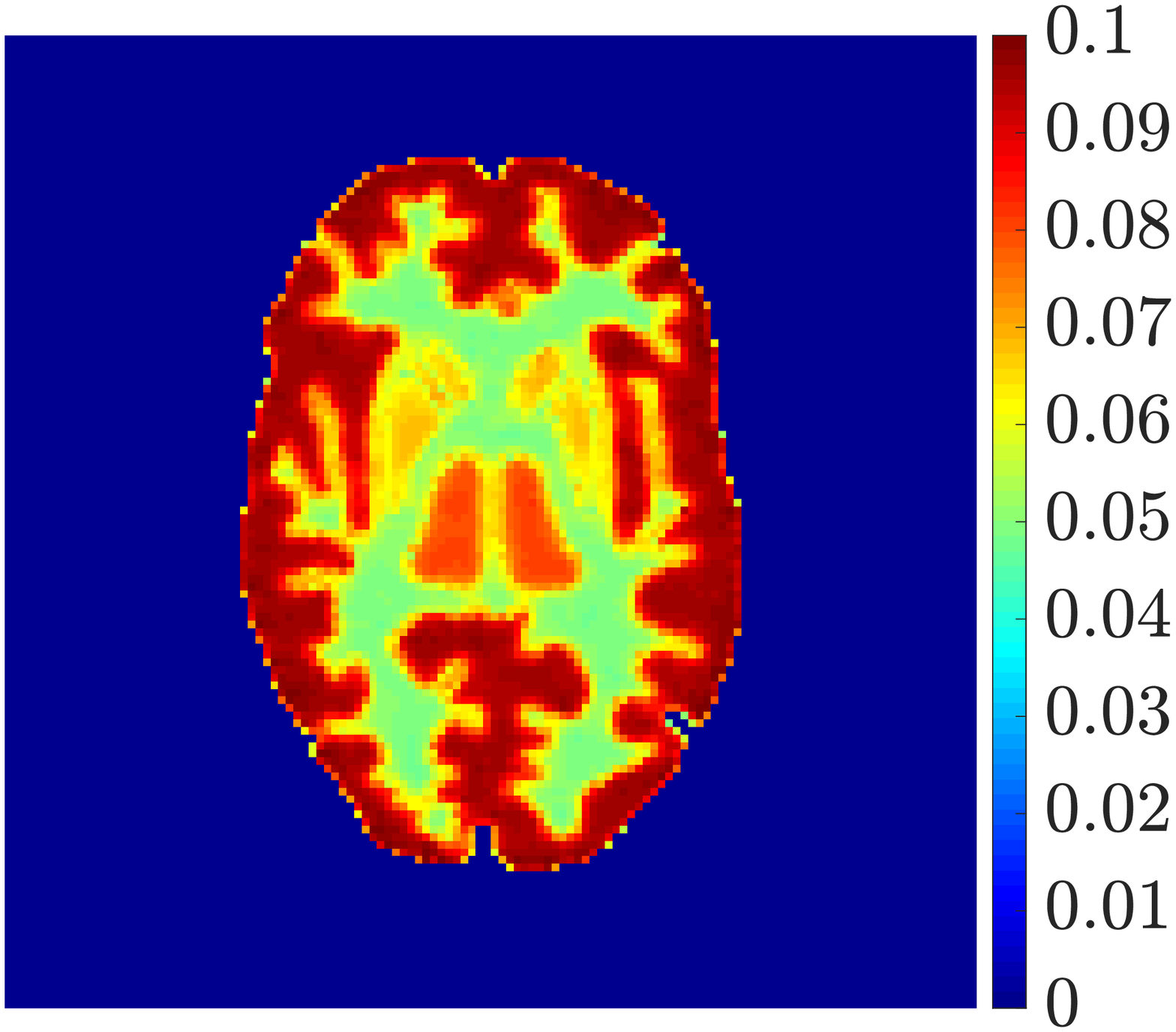}\quad  \hspace{0.55cm}
	\includegraphics[width=0.16\textwidth]{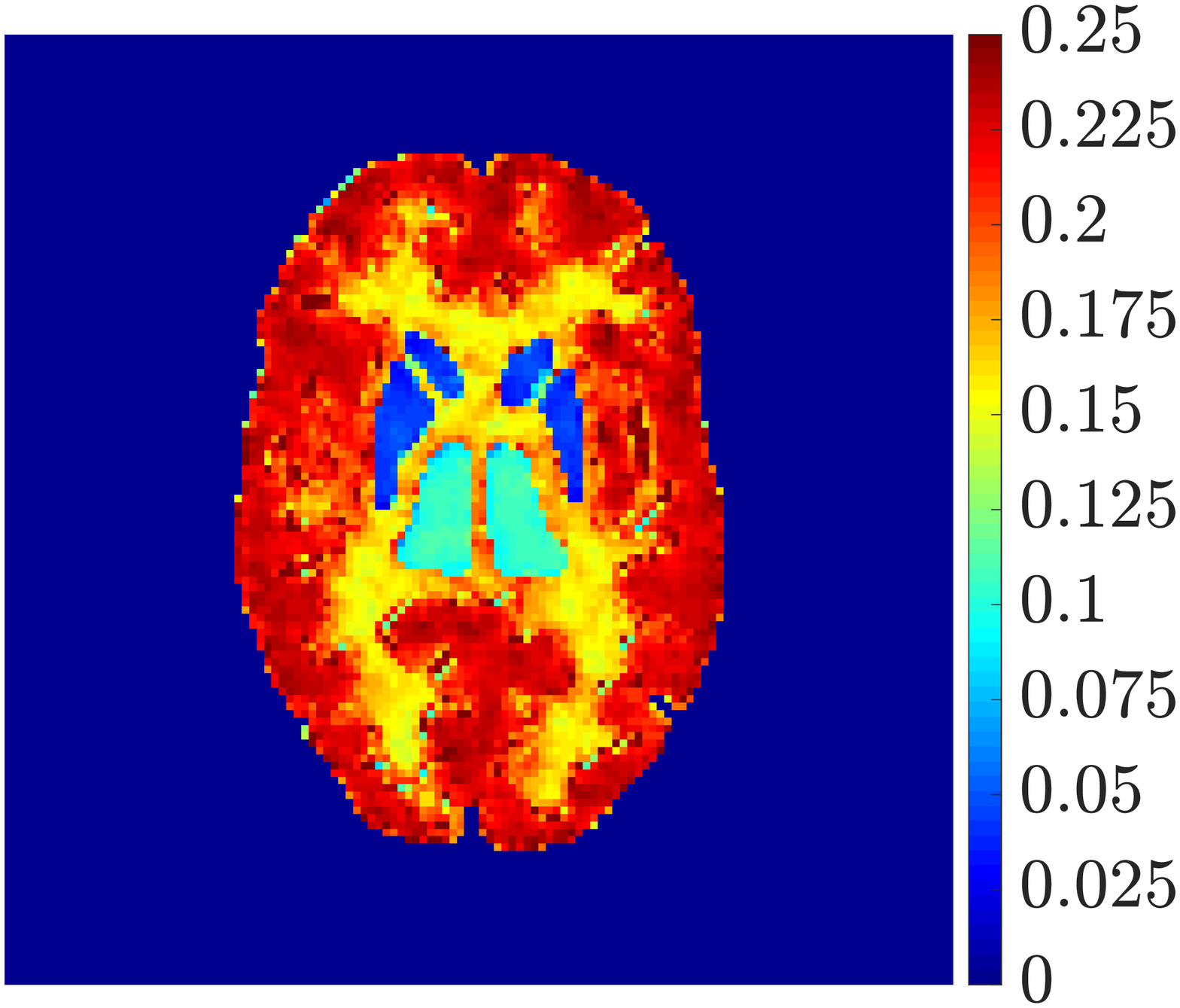}\quad   \hspace{0.55cm}
	\includegraphics[width=0.16\textwidth]{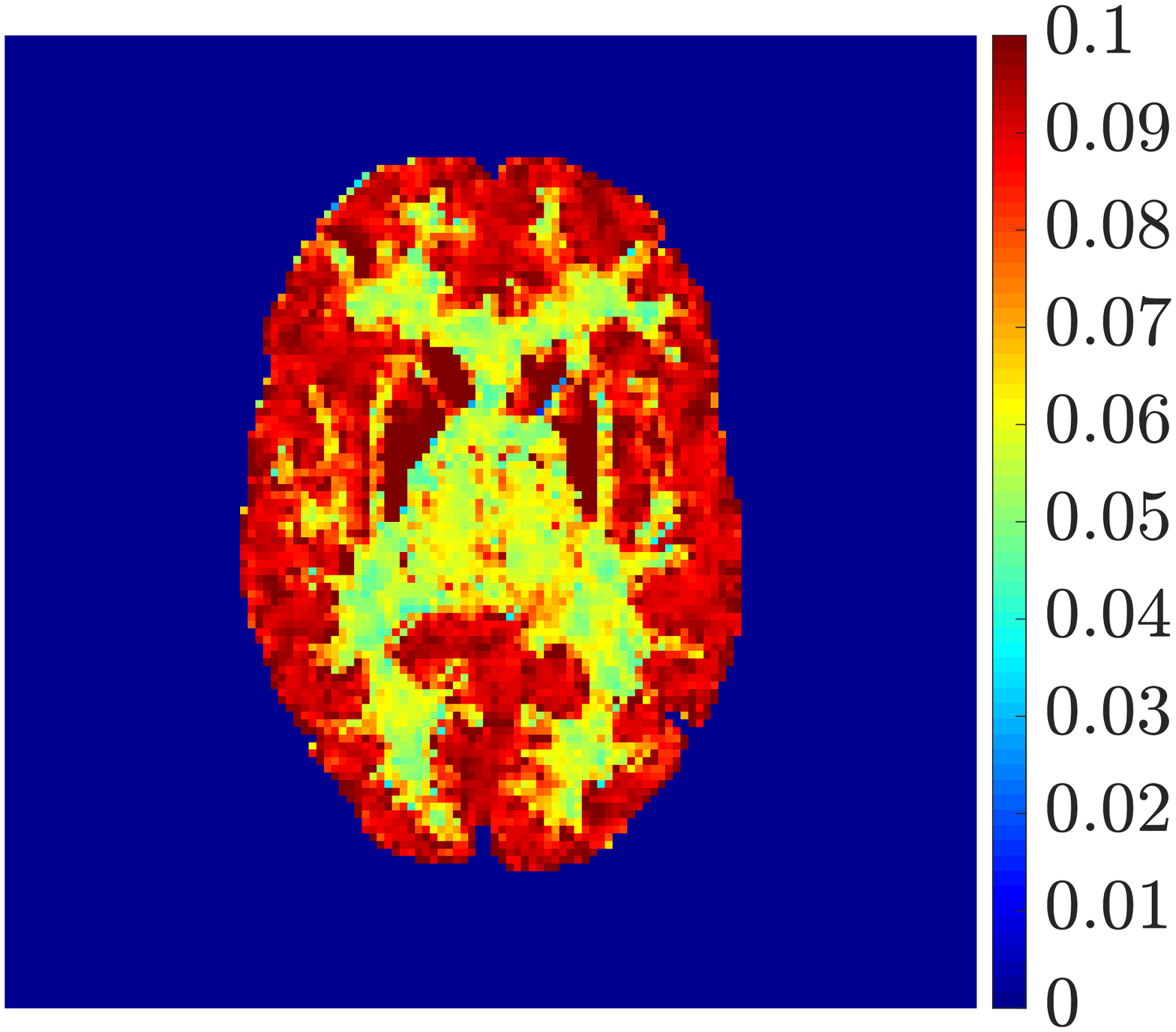}\quad   \hspace{0.55cm}
	\includegraphics[width=0.16\textwidth]{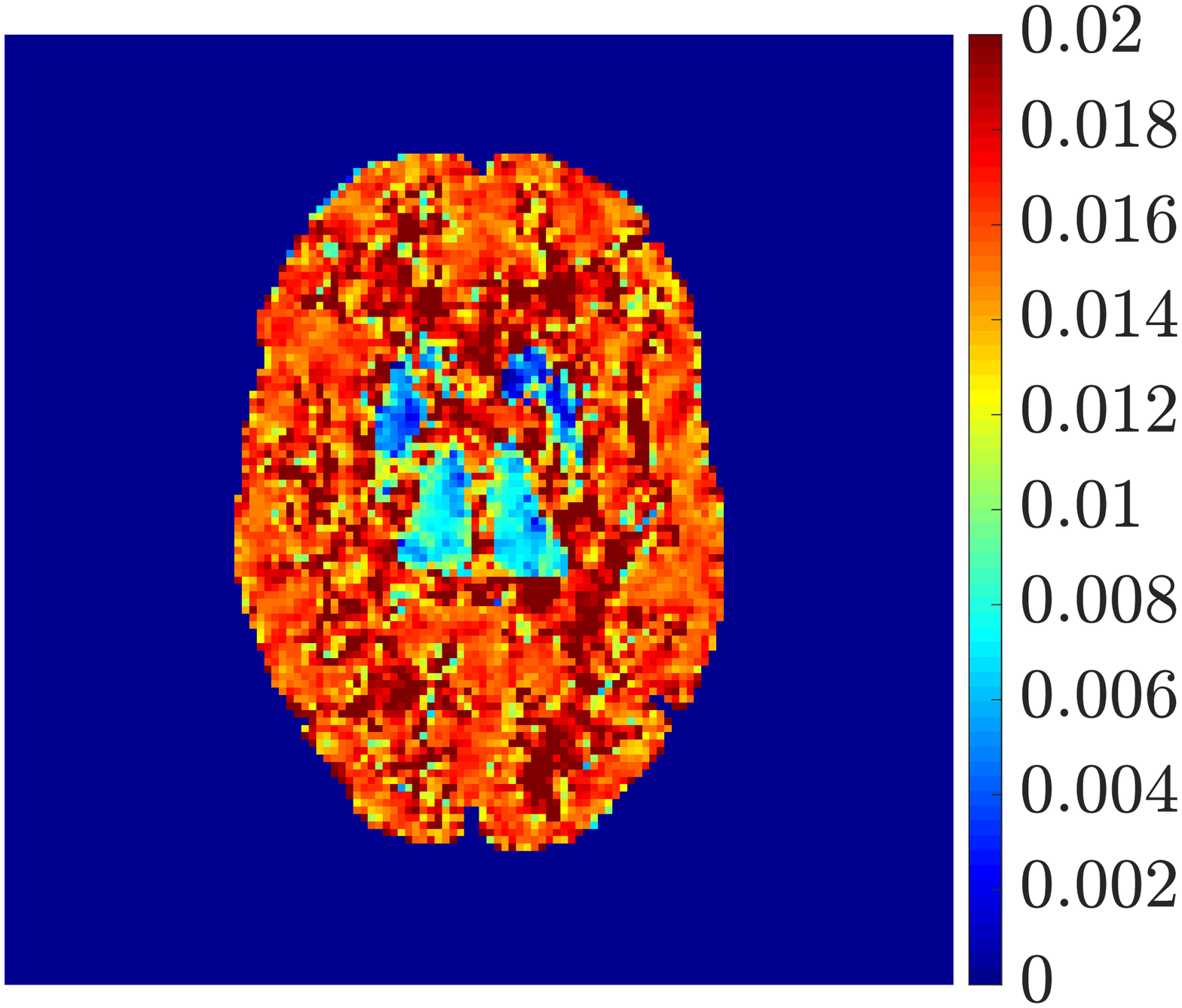}\\
	
	\includegraphics[width=0.16\textwidth]{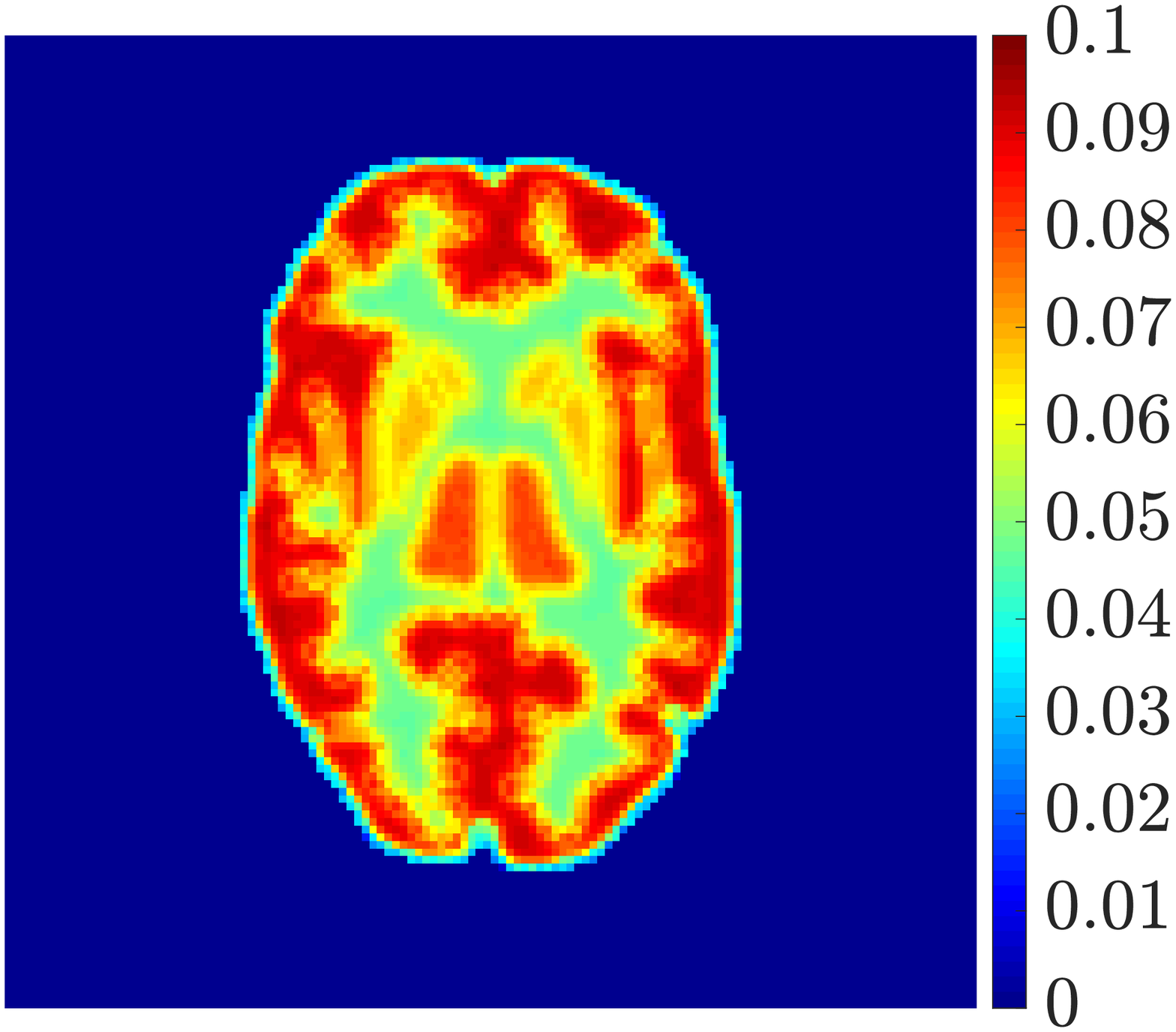}\quad   \hspace{0.55cm}
	\includegraphics[width=0.16\textwidth]{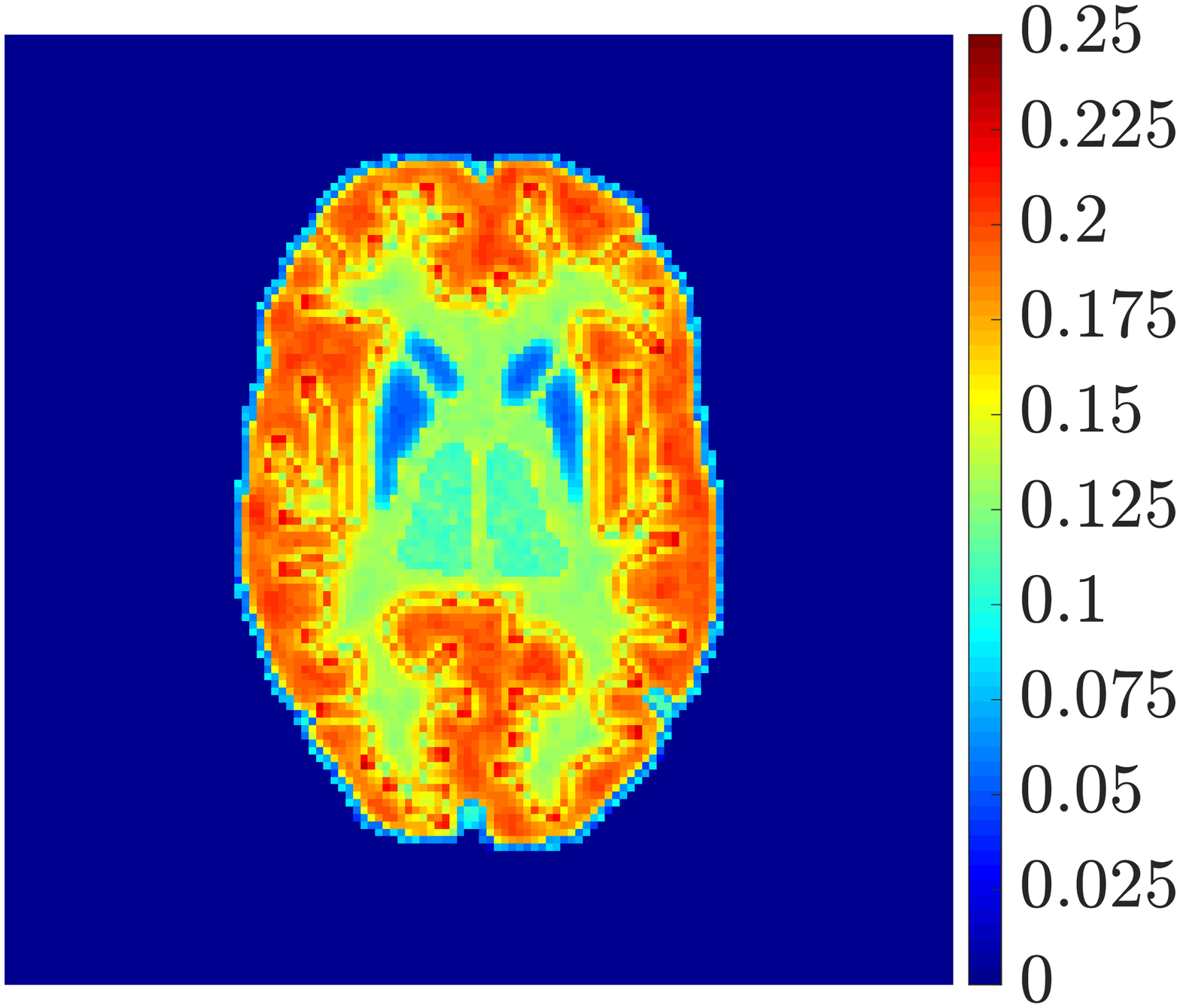}\quad    \hspace{0.55cm}
	\includegraphics[width=0.16\textwidth]{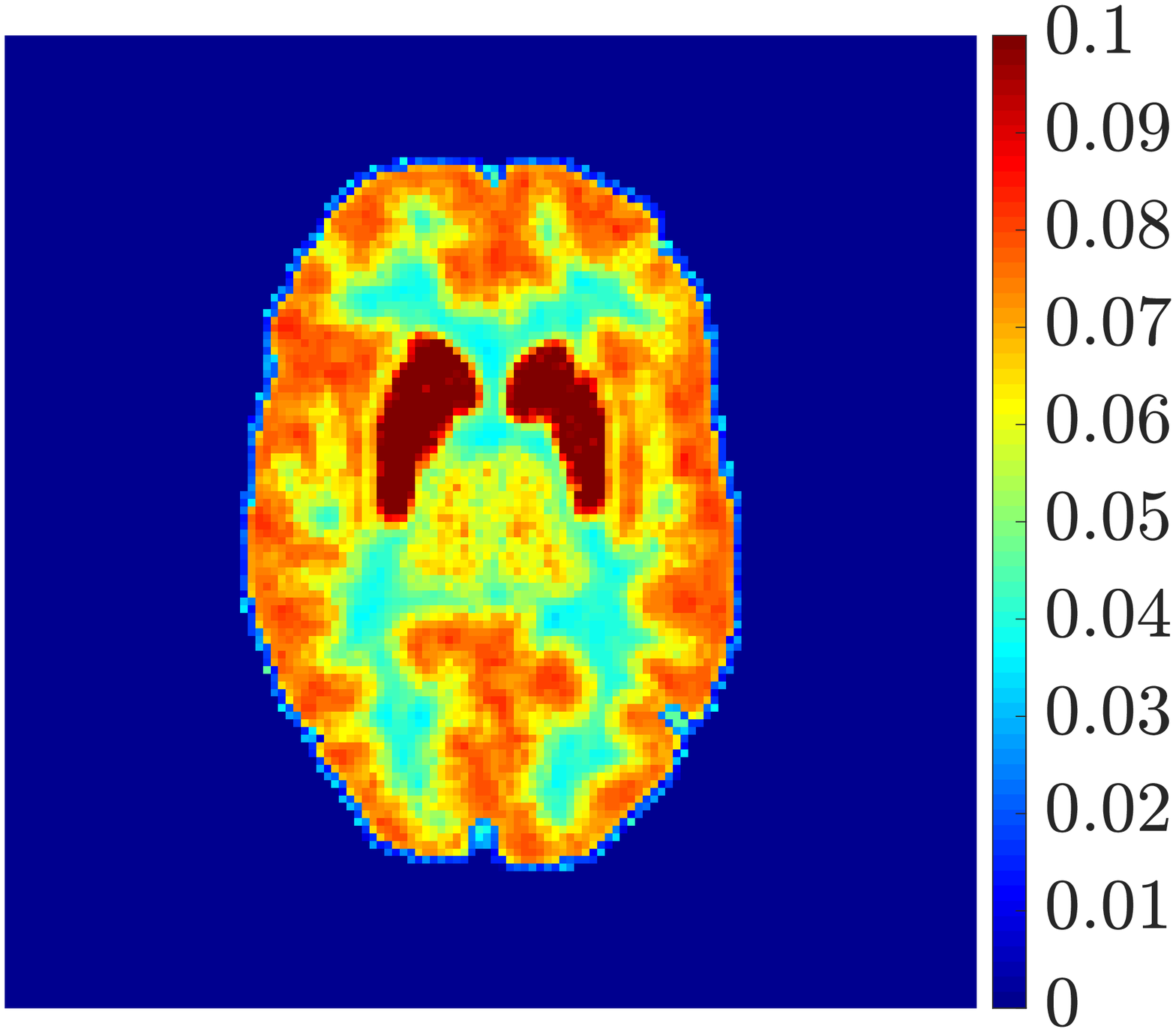}\quad    \hspace{0.55cm}
	\includegraphics[width=0.16\textwidth]{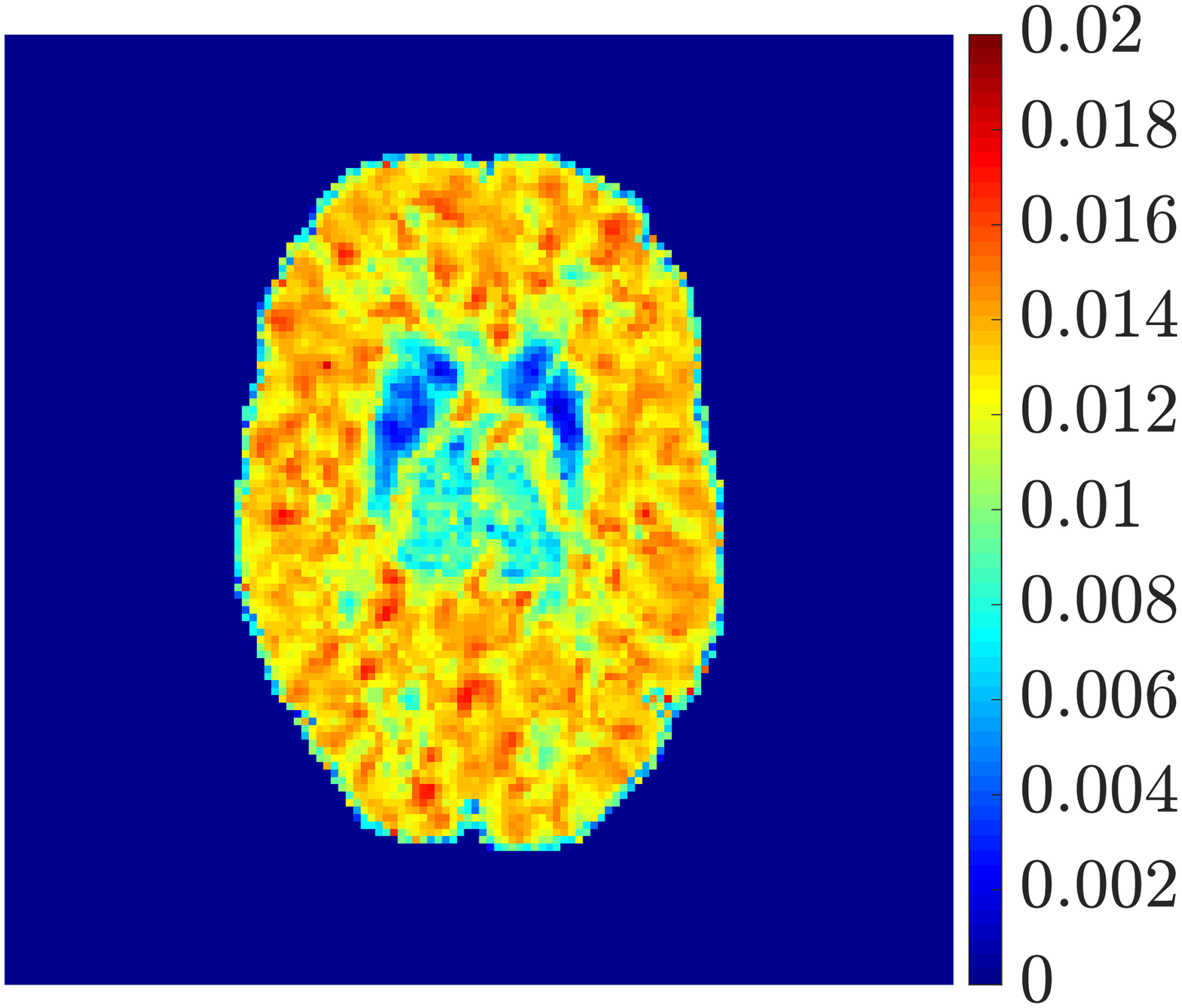}\\
	
	\includegraphics[width=0.16\textwidth]{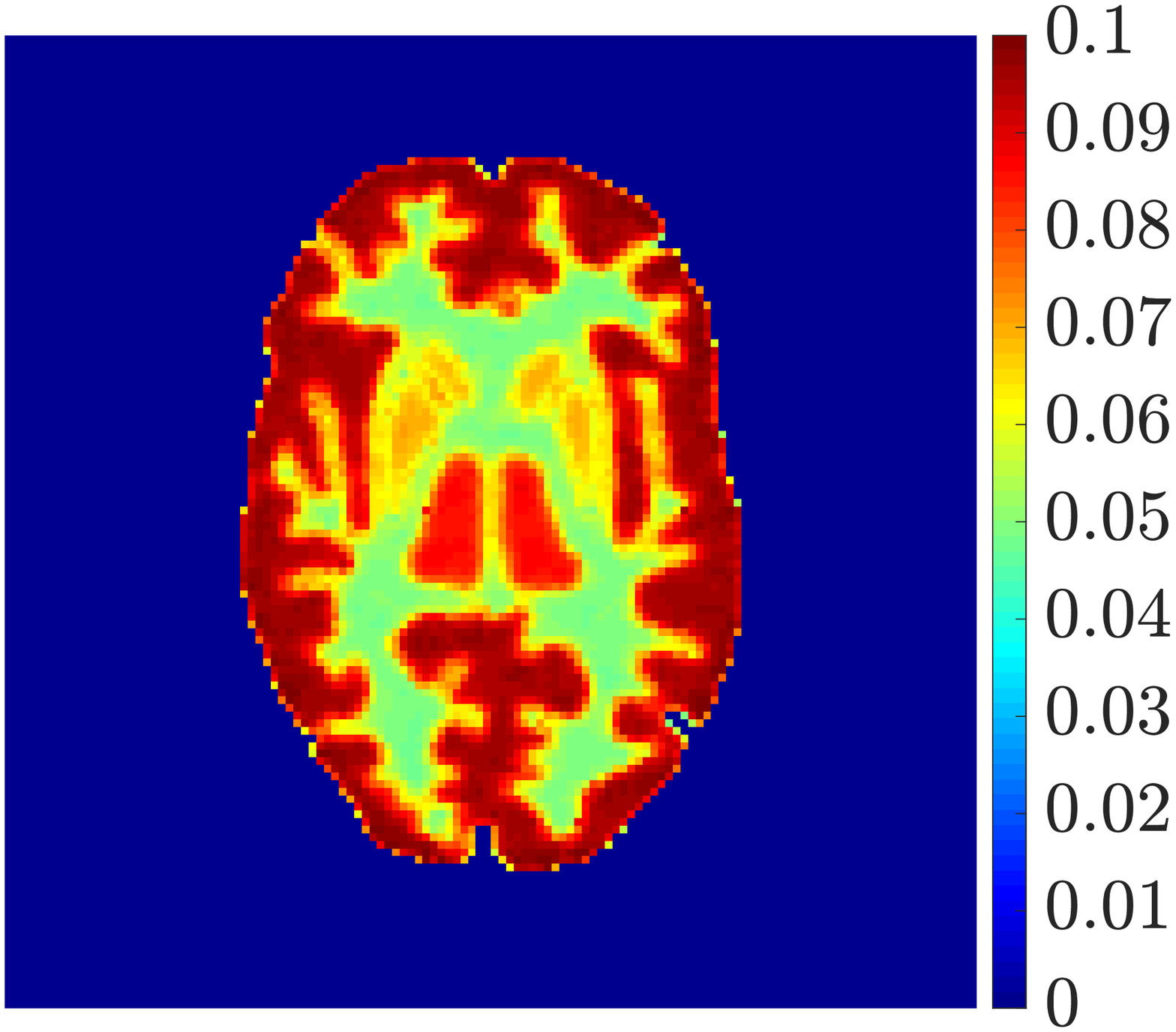}\quad    \hspace{0.55cm}
	\includegraphics[width=0.16\textwidth]{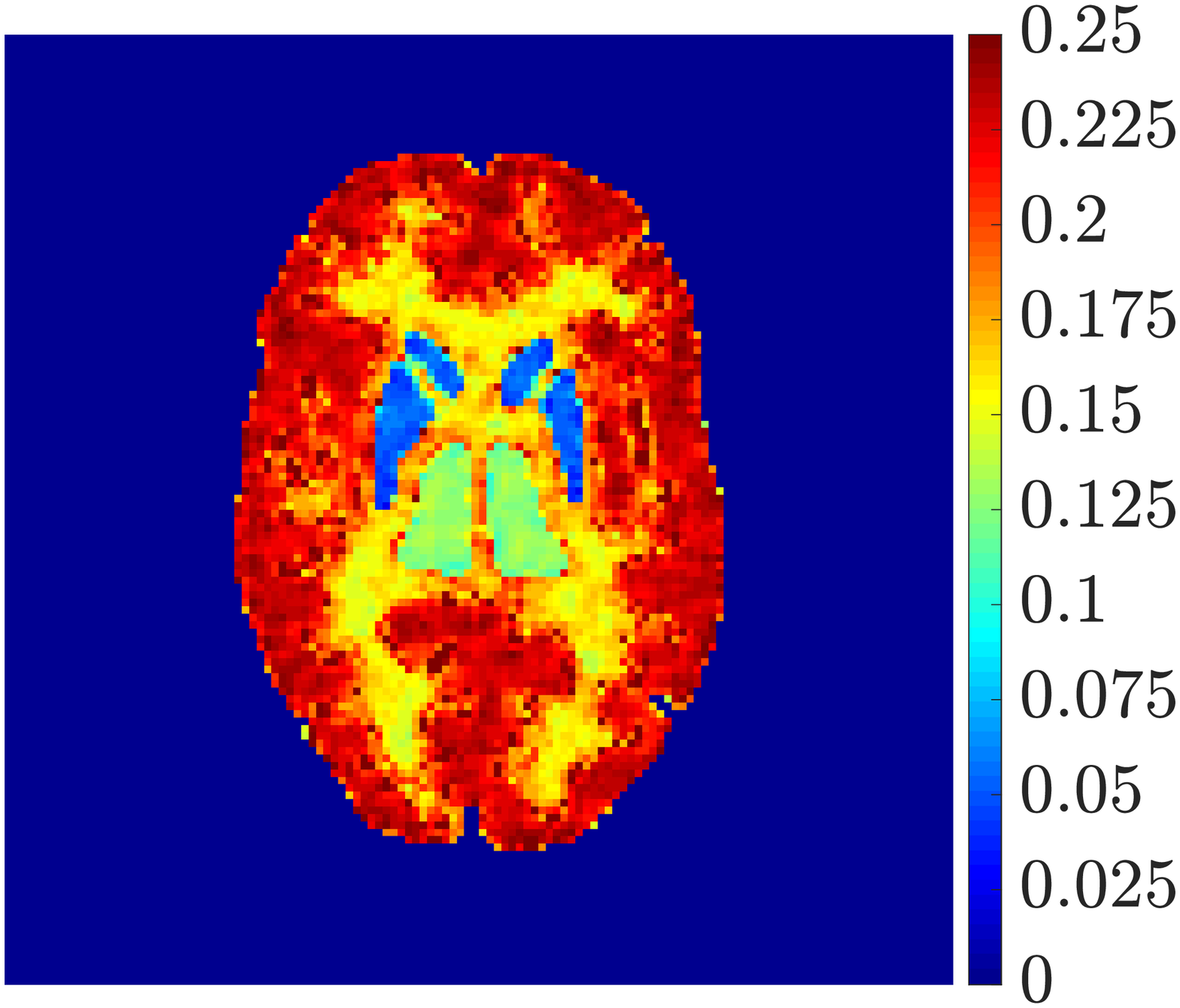}\quad    \hspace{0.55cm}
	\includegraphics[width=0.16\textwidth]{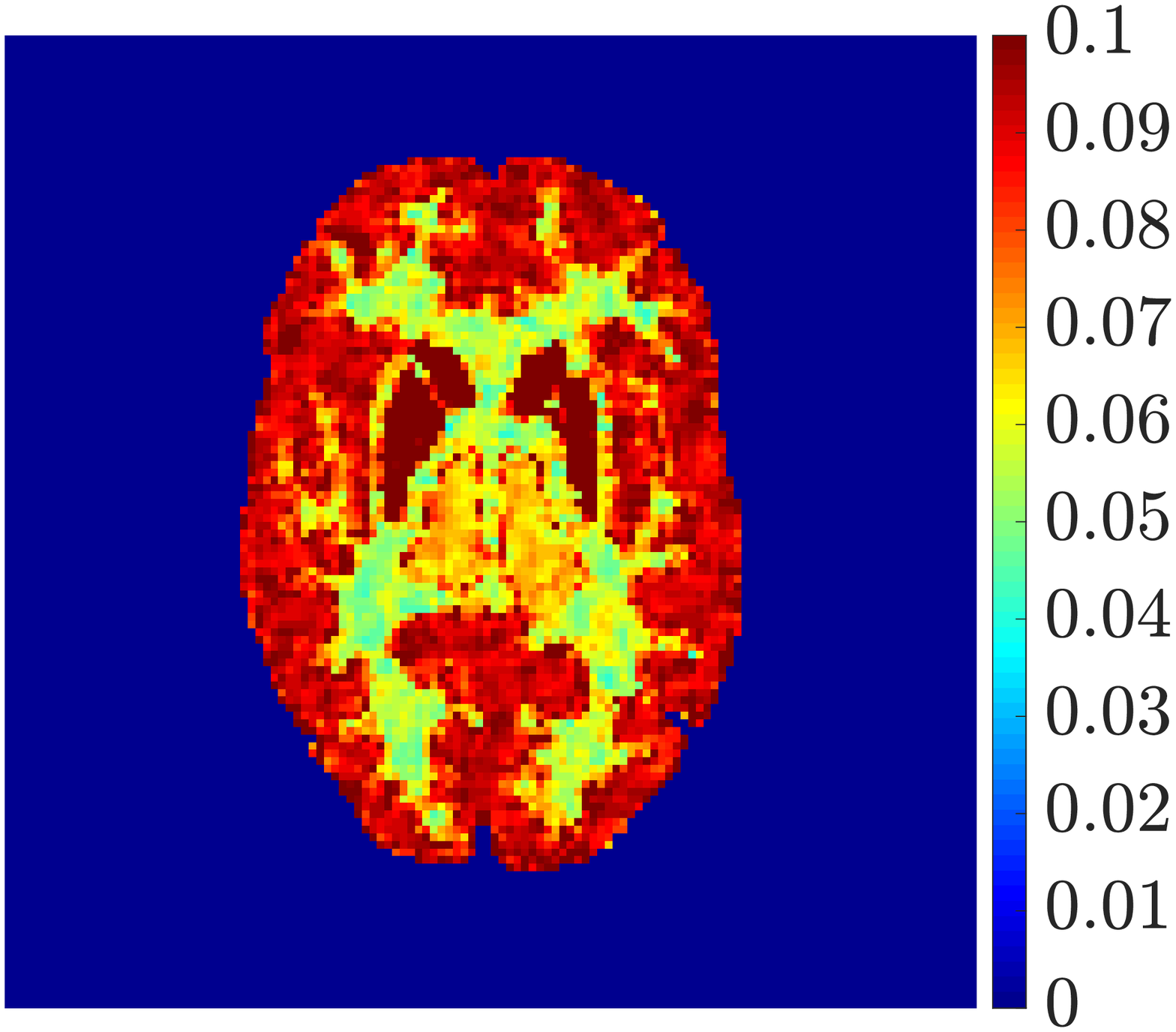}\quad    \hspace{0.55cm}
	\includegraphics[width=0.16\textwidth]{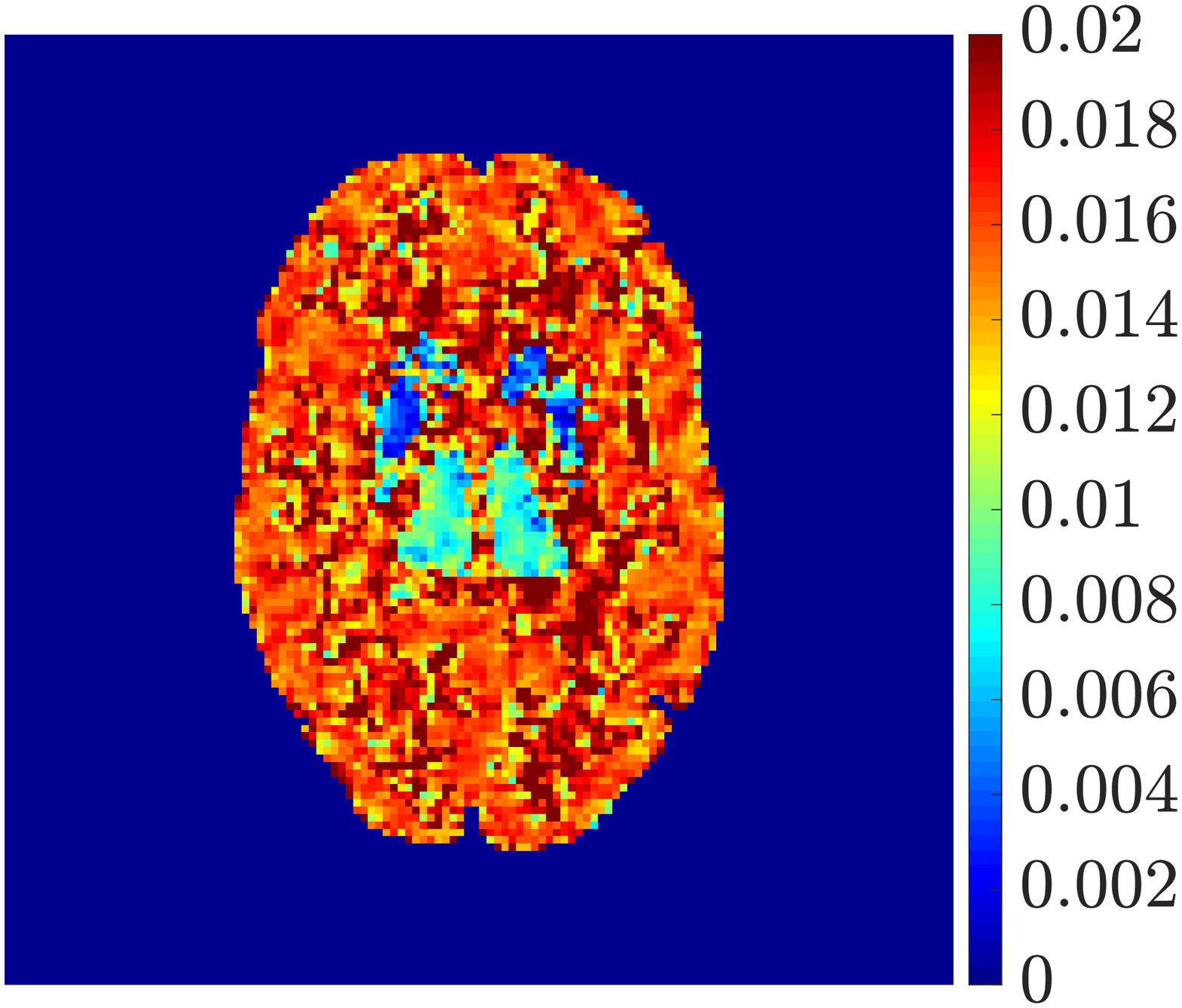}
		\caption{From left to right: mean parametric images corresponding to $k_1$, $k_2$, $k_3$, $k_4$, obtained by using reg-AS-TR (first row), reg-GN (second row), \emph{lsqcurvefit} (third row). Case $10\%$-noise IF.}
		\label{fig:k-rec-noise-10p}
	\end{figure}

\begin{figure}
		\centering
	    \includegraphics[width=0.16\textwidth]{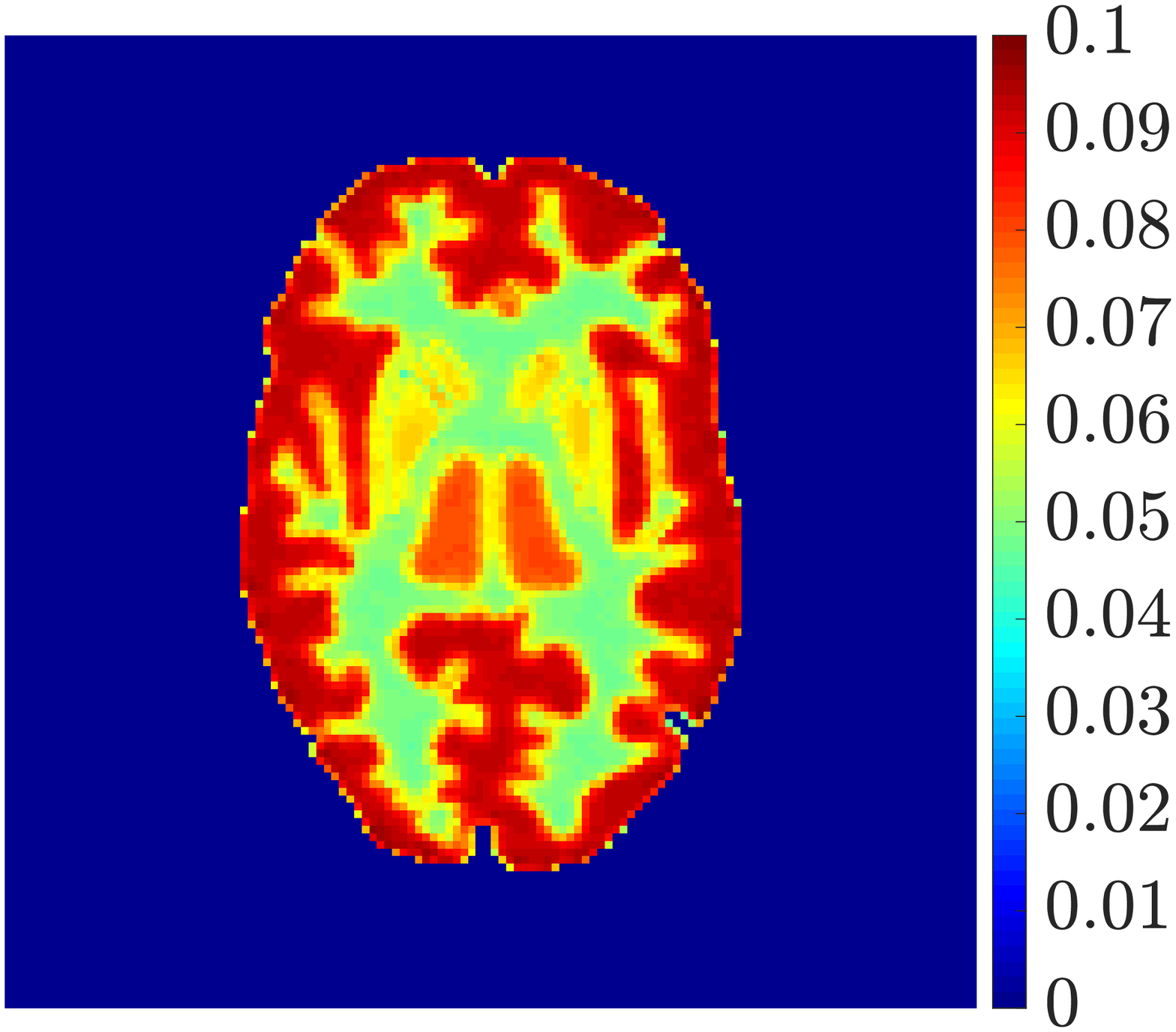}\quad    \hspace{0.55cm}
		\includegraphics[width=0.16\textwidth]{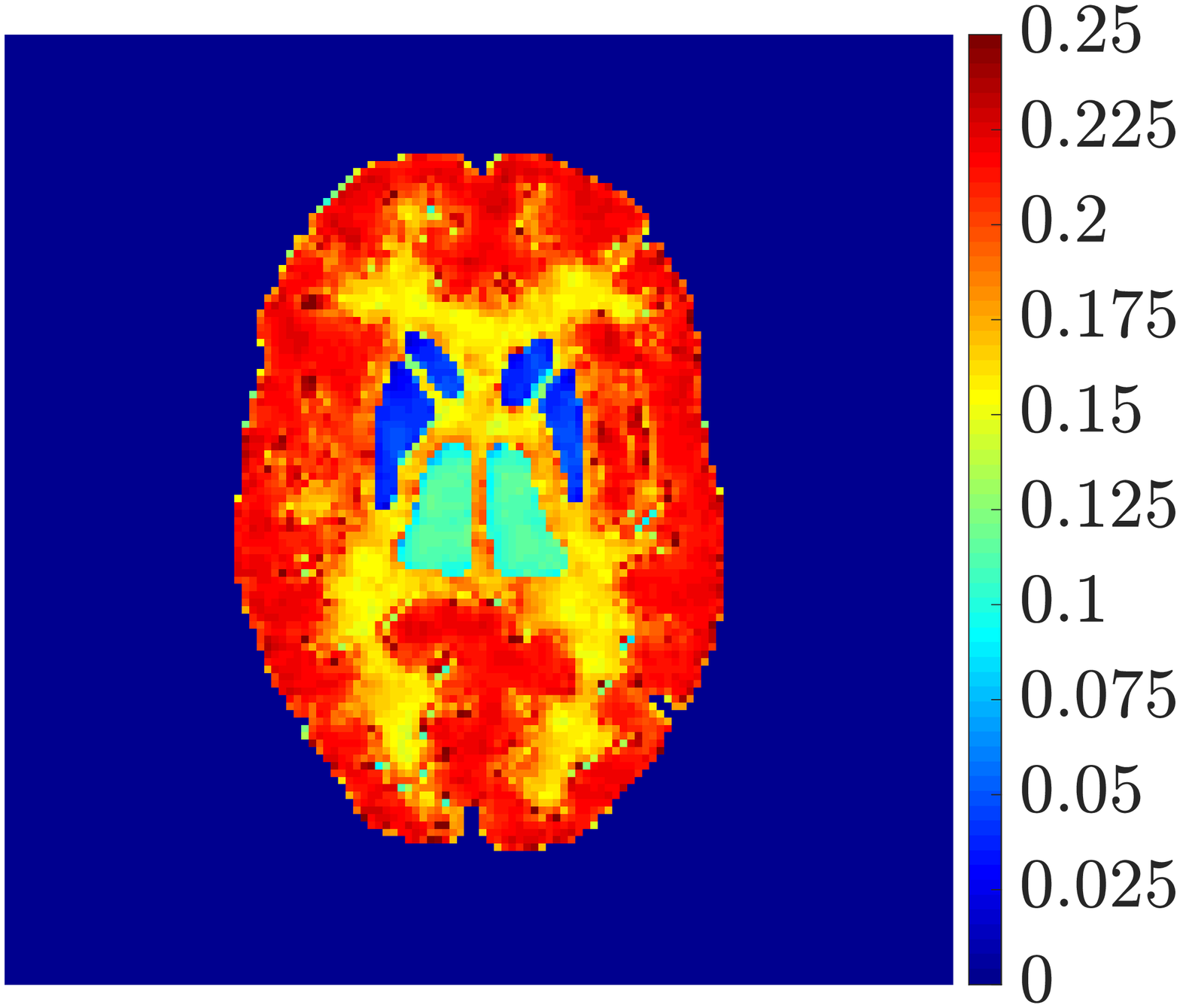}\quad   \hspace{0.55cm}
		\includegraphics[width=0.16\textwidth]{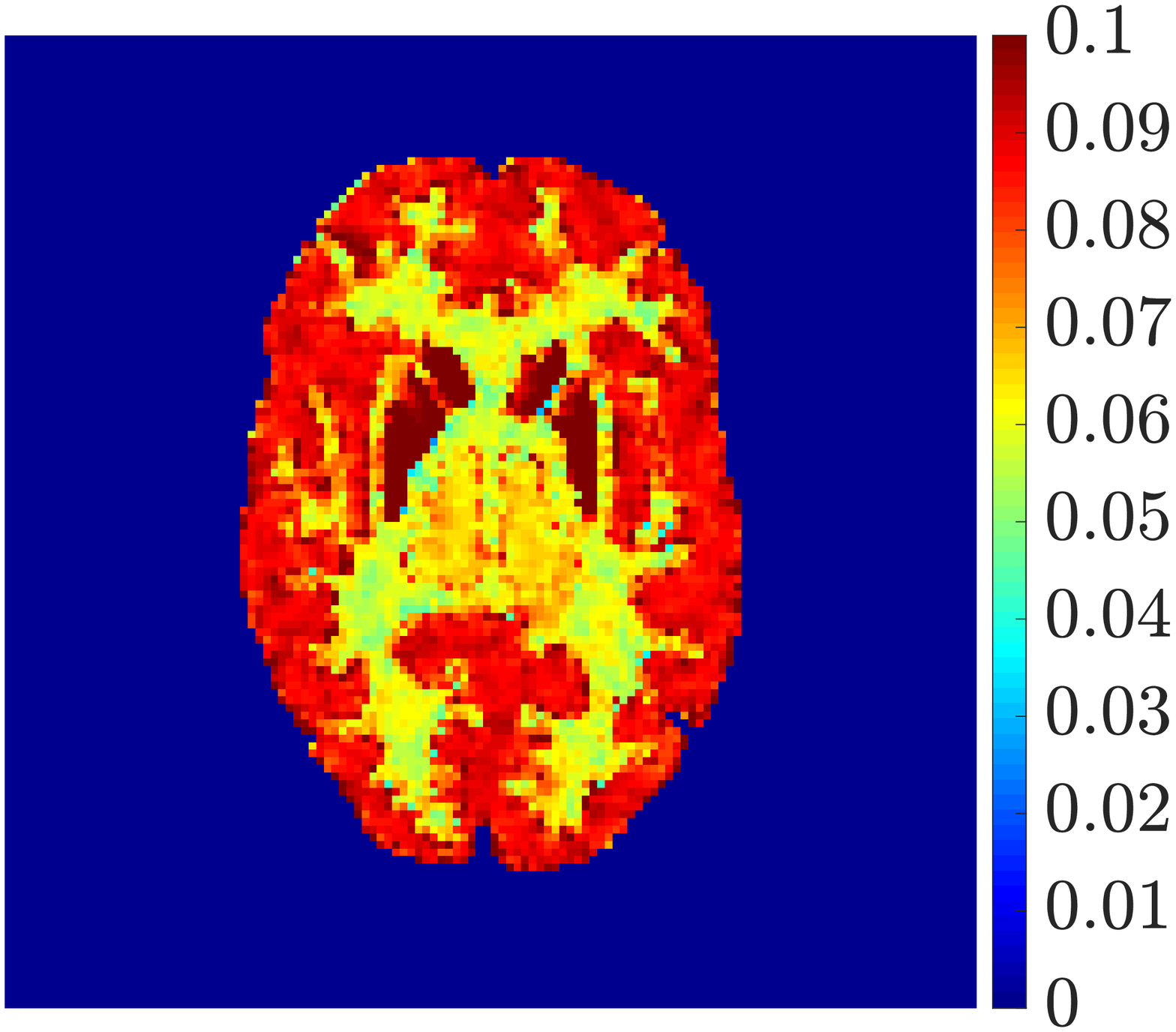}\quad     \hspace{0.55cm}
		\includegraphics[width=0.16\textwidth]{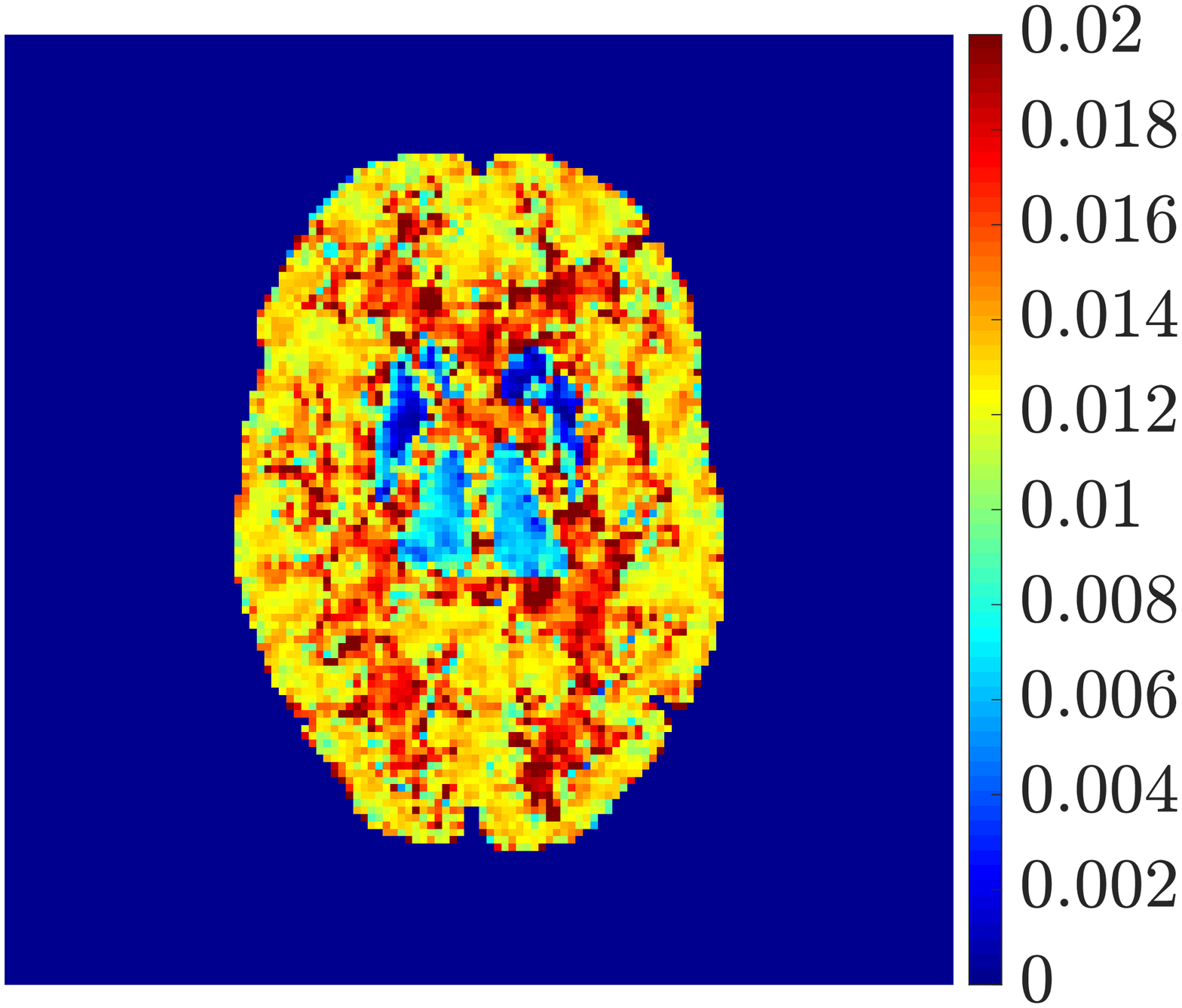}\\
		
		\includegraphics[width=0.16\textwidth]{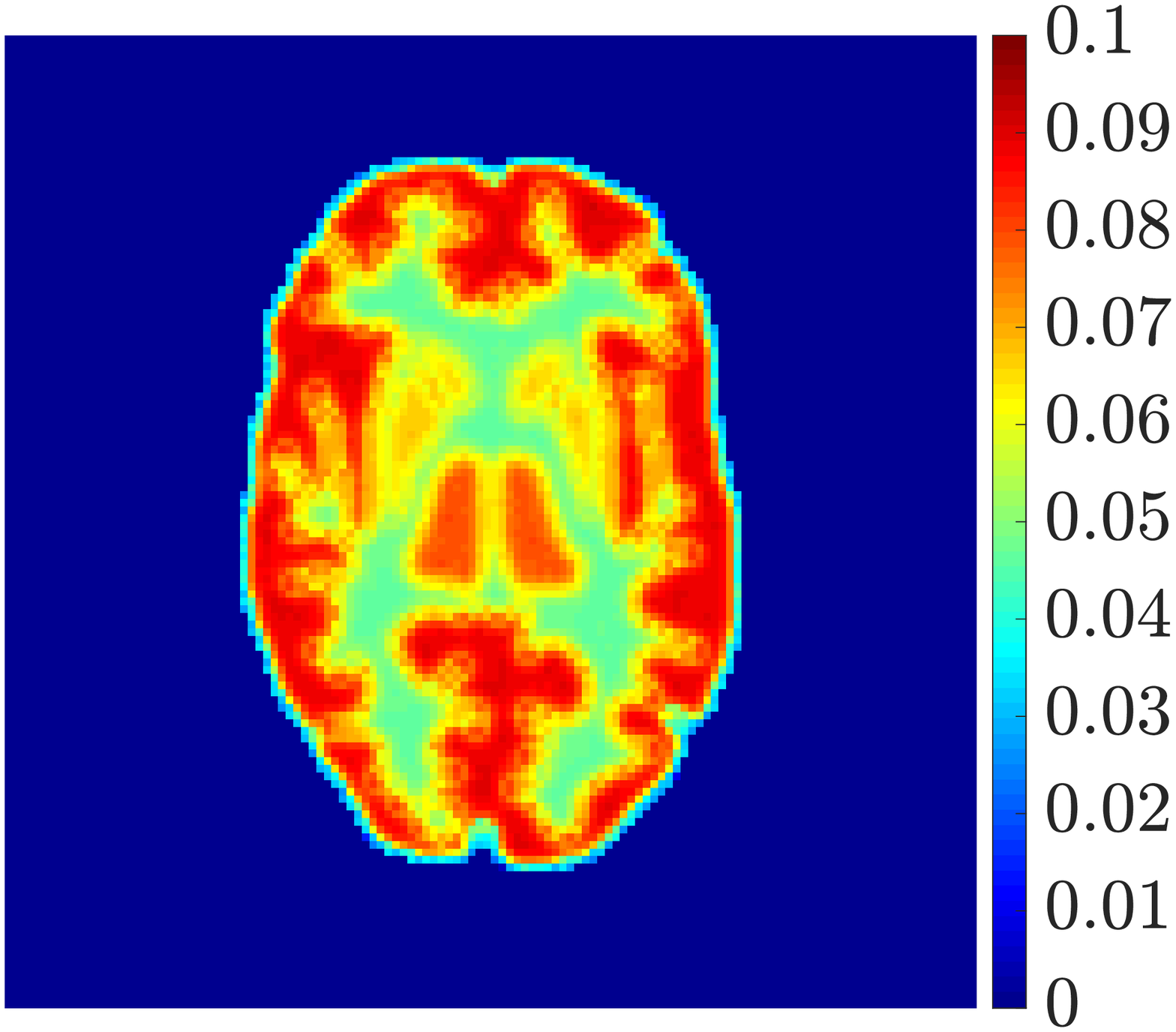}\quad    \hspace{0.55cm}
		\includegraphics[width=0.16\textwidth]{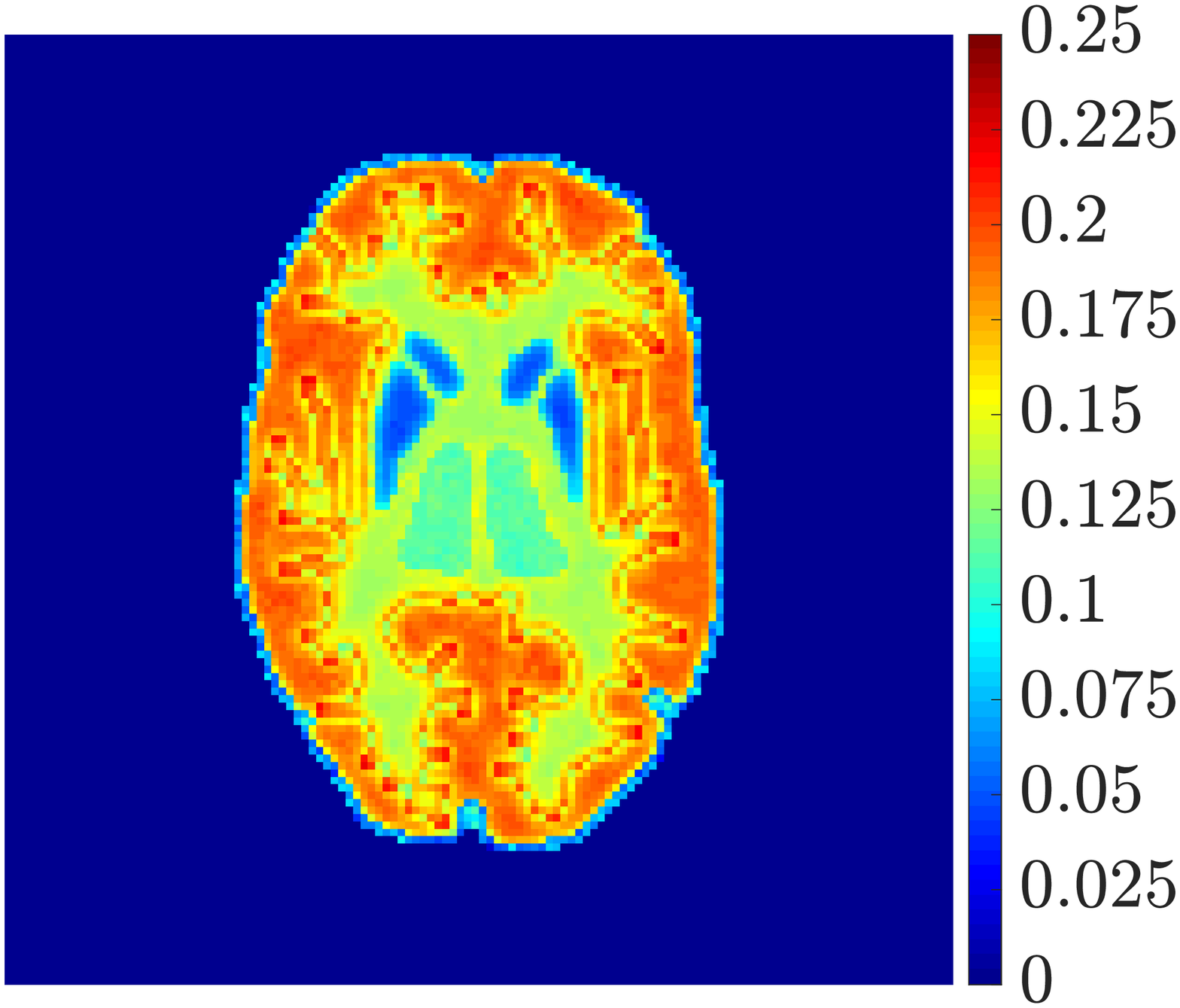}\quad    \hspace{0.55cm}
		\includegraphics[width=0.16\textwidth]{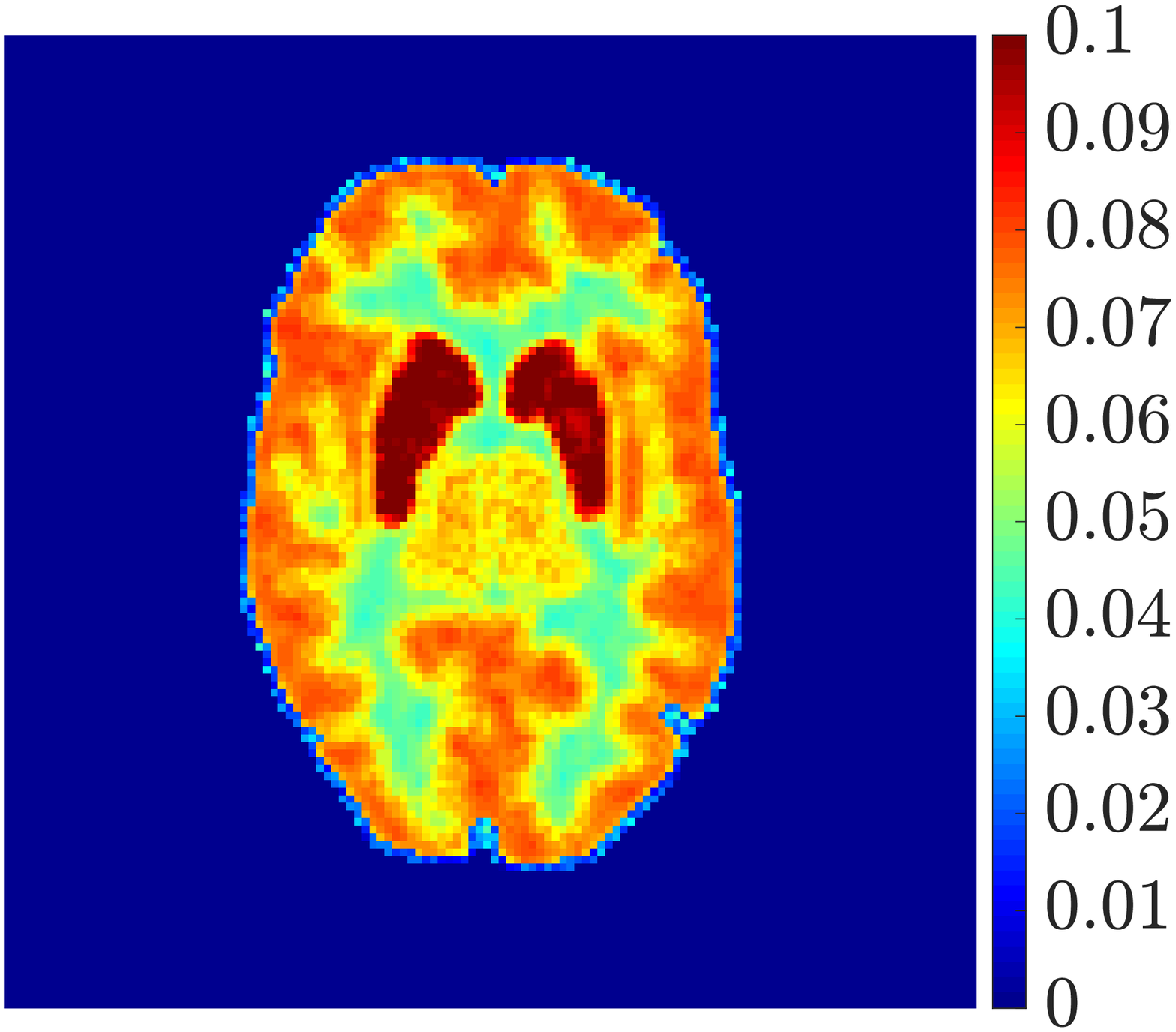}\quad    \hspace{0.55cm}
		\includegraphics[width=0.16\textwidth]{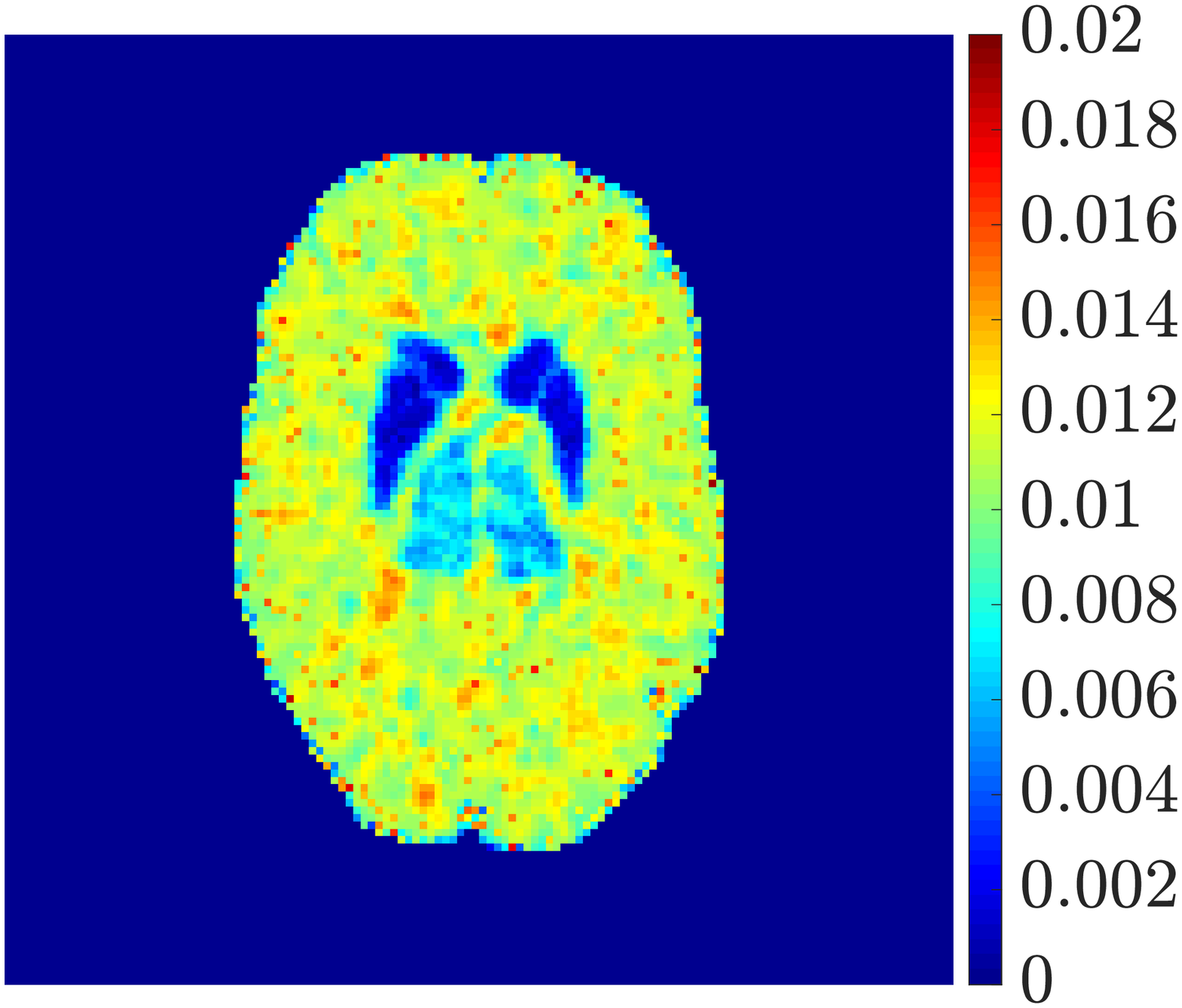}\\
		
		\includegraphics[width=0.16\textwidth]{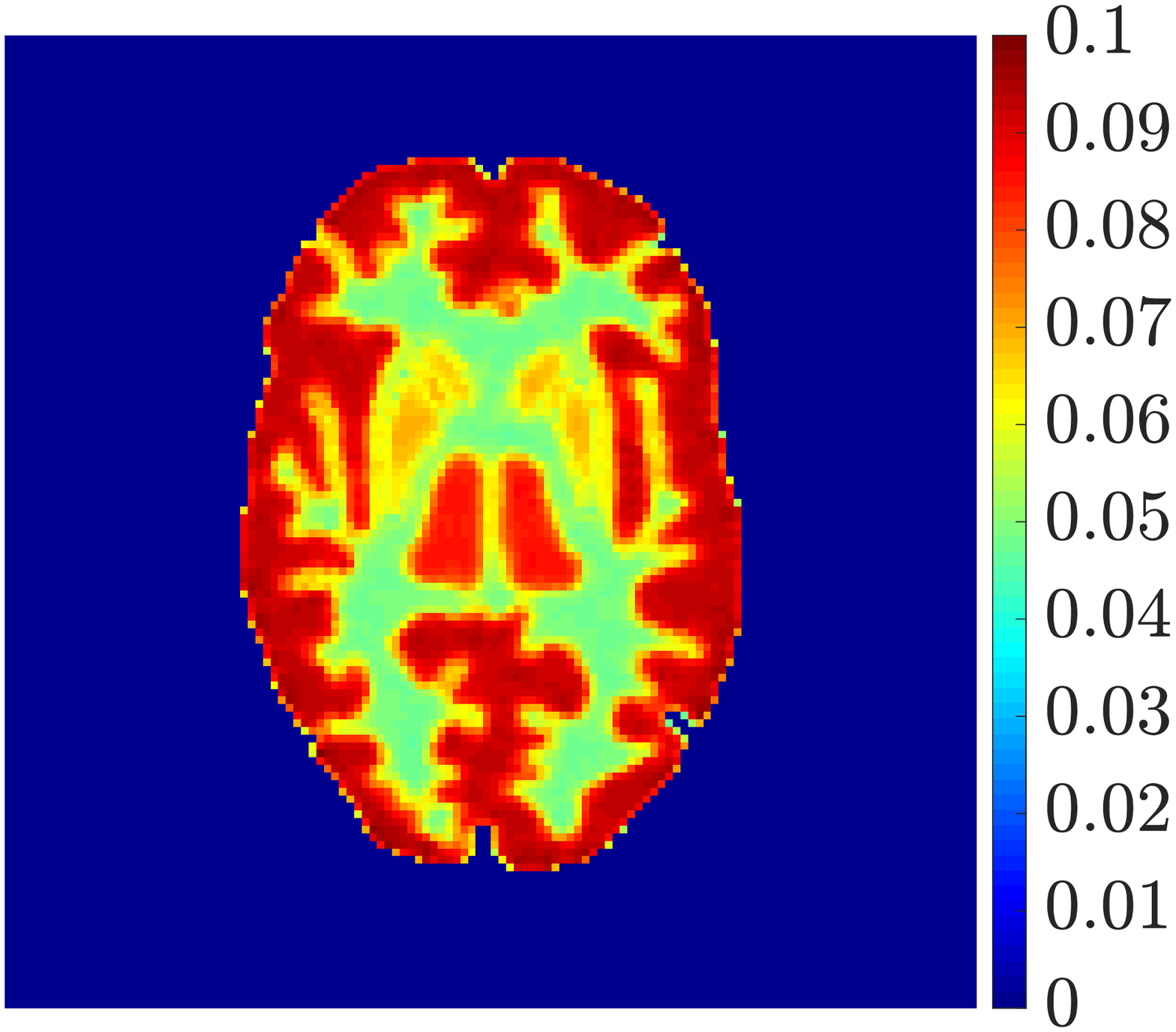}\quad   \hspace{0.55cm}
		\includegraphics[width=0.16\textwidth]{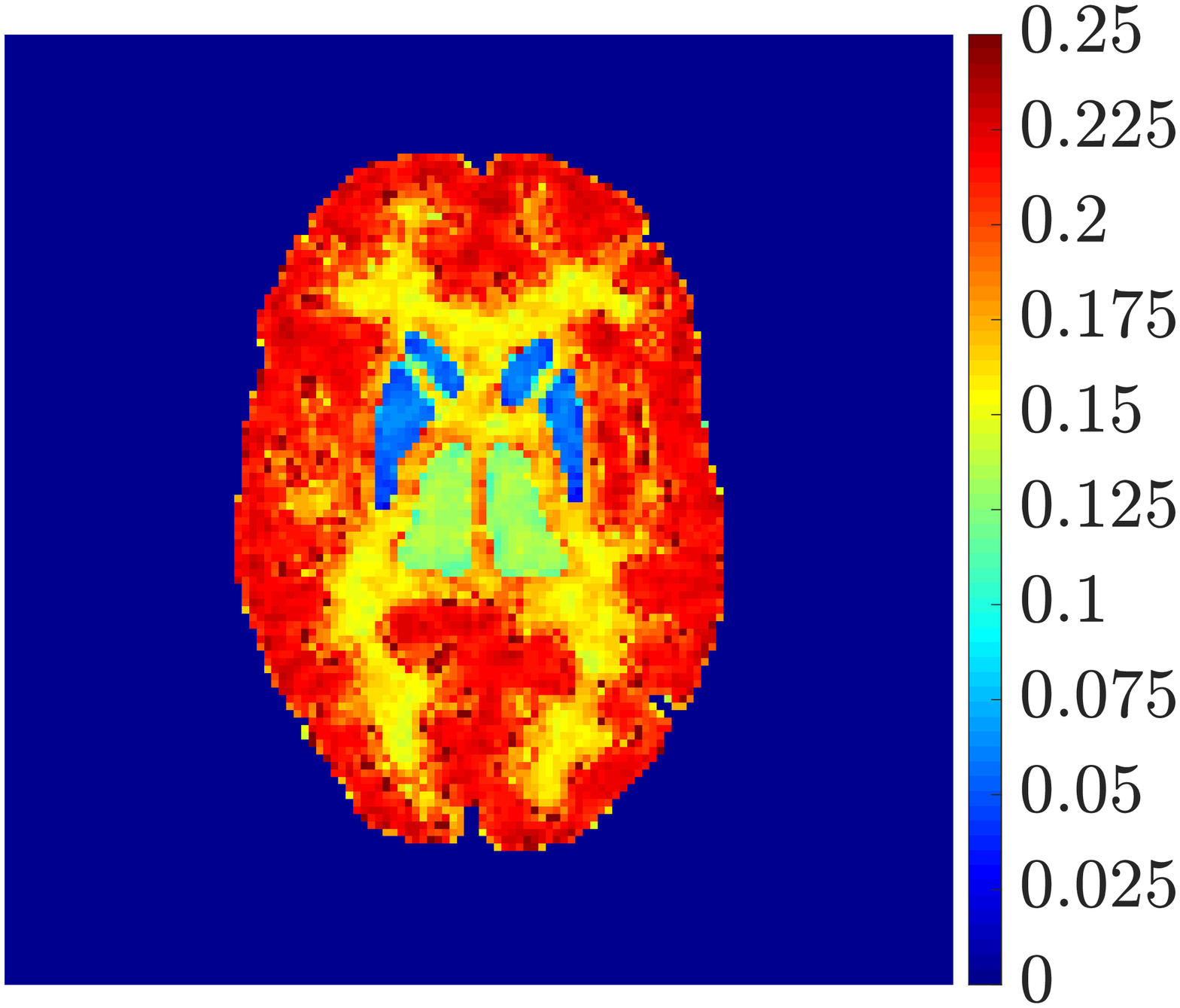}\quad   \hspace{0.55cm}
		\includegraphics[width=0.16\textwidth]{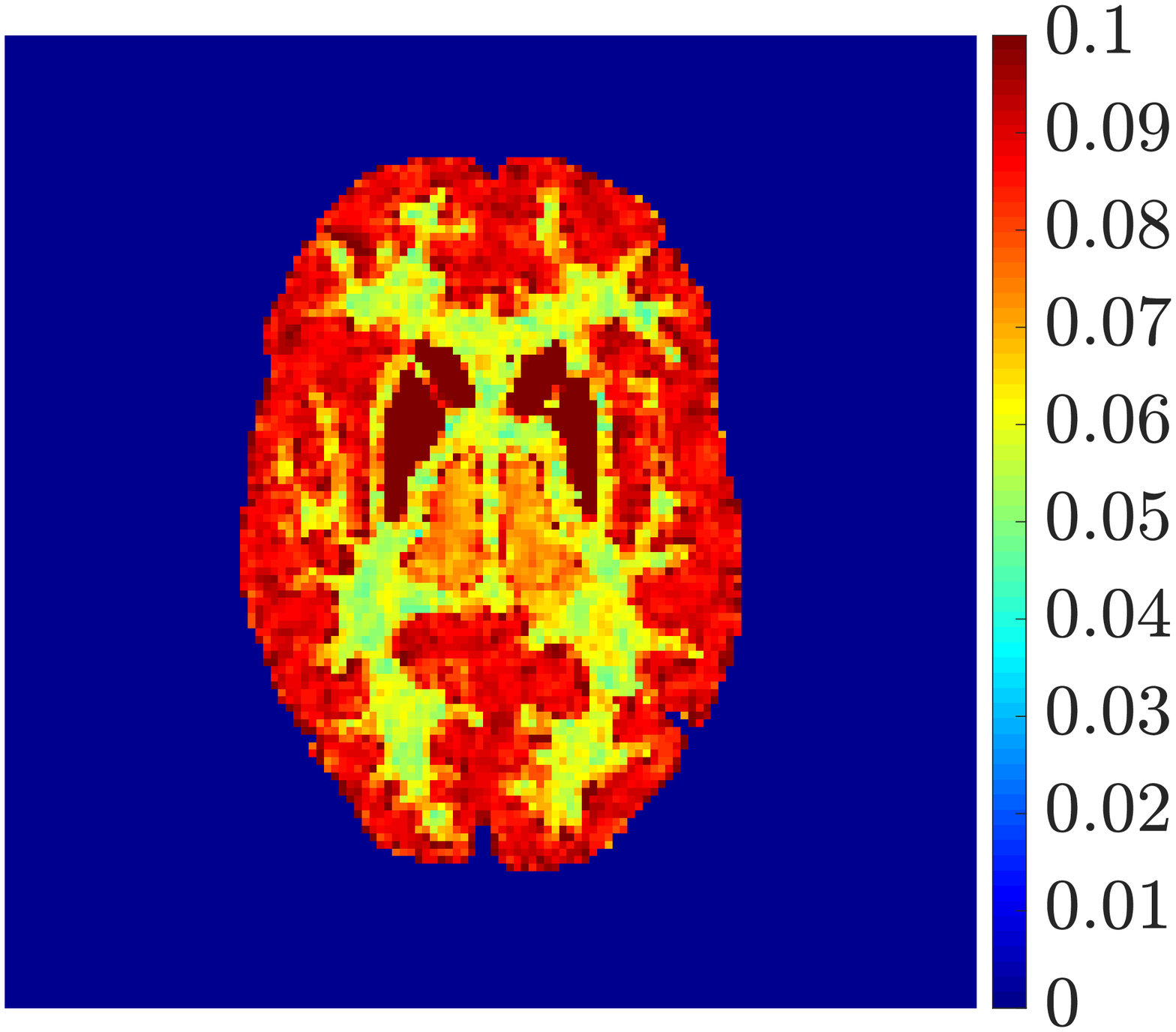}\quad    \hspace{0.55cm}
		\includegraphics[width=0.16\textwidth]{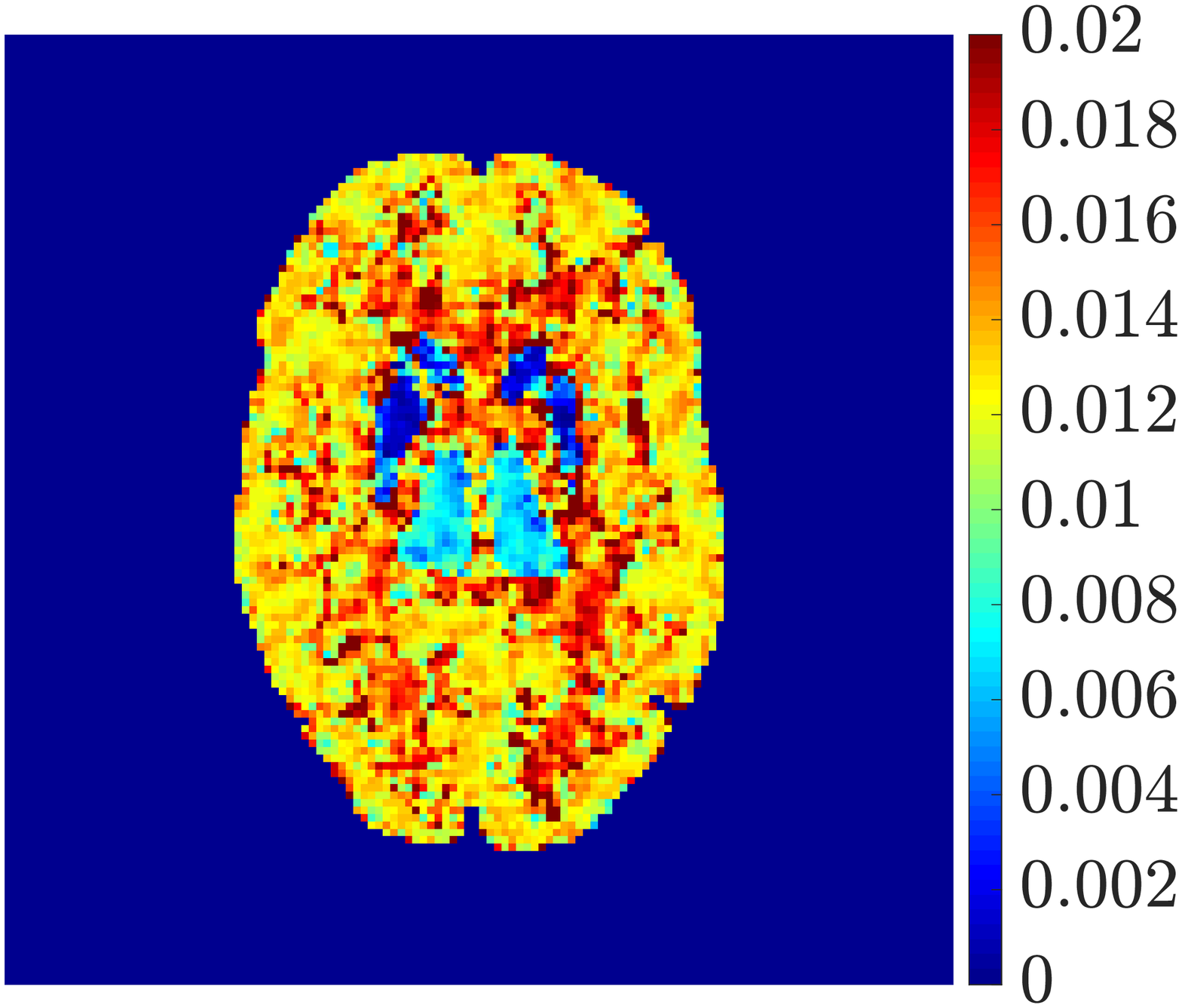}
		\caption{From left to right: mean parametric images corresponding to $k_1$, $k_2$, $k_3$, $k_4$, obtained by using reg-AS-TR (first row), reg-GN (second row), \emph{lsqcurvefit} (third row). Case $20\%$-noise IF.}
		\label{fig:k-rec-noise-20p}
	\end{figure}

Figure \ref{fig:k-rec-noise-free}, Figure \ref{fig:k-rec-noise-10p}, and Figure \ref{fig:k-rec-noise-20p} show the mean images computed over the ten reconstructions obtained by the methods reg-AS-TR, reg-GN, and by the Matlab routine \emph{lsqcurvefit} implementing a standard Trust-Region-Reflective least-squares algorithm \cite{Coleman,Coleman1994}. We used the noise-free IF and the perturbed IF with $10\%$ and $20\%$ of noise, respectively. 
 \figurename~\ref{fig:k-rec-all} contains mean and standard deviation values of the kinetic parameters computed over the ten reconstructions and over each one of the four homogeneous regions, for each one of the three noise levels on the IF.

Finally, \figurename~\ref{fig:data-rec} represents the last frame of the dynamic PET data reconstructed with the mean parametric values returned by reg-AS-TR, reg-GN, and \emph{lsqcurvefit}, with respect to the noise-free, 10$\%$-noise, and 20$\%$-noise IFs.

\section{Comments and conclusions}

\begin{figure}
	\centering
	\subfigure[$k_1$ \label{subfig:k1-rec}]
	{\includegraphics[width=0.35\textwidth]{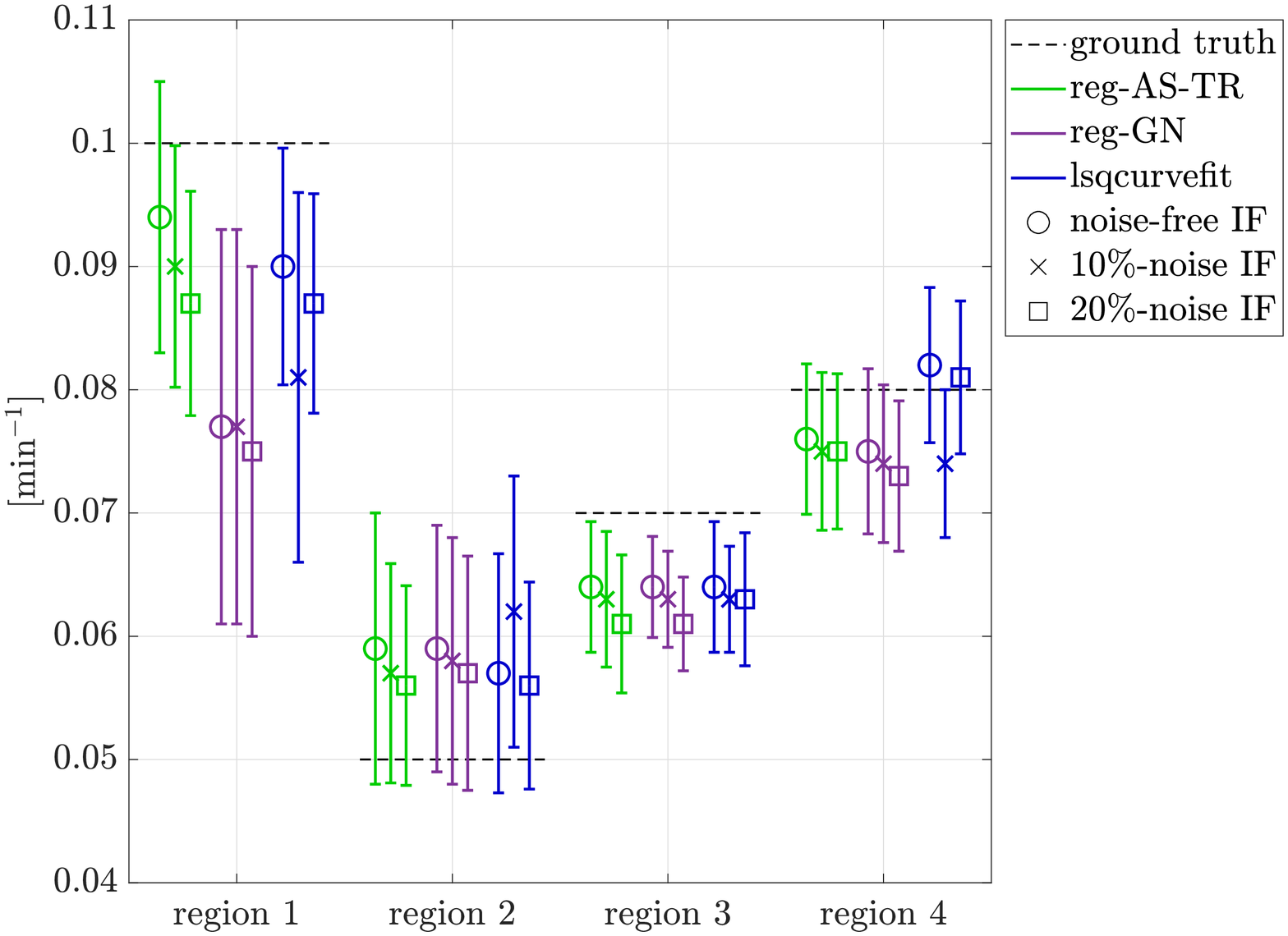}} \quad  \hspace{0.7cm}
	\subfigure[$k_2$ \label{subfig:k2-rec}]
	{\includegraphics[width=0.35\textwidth]{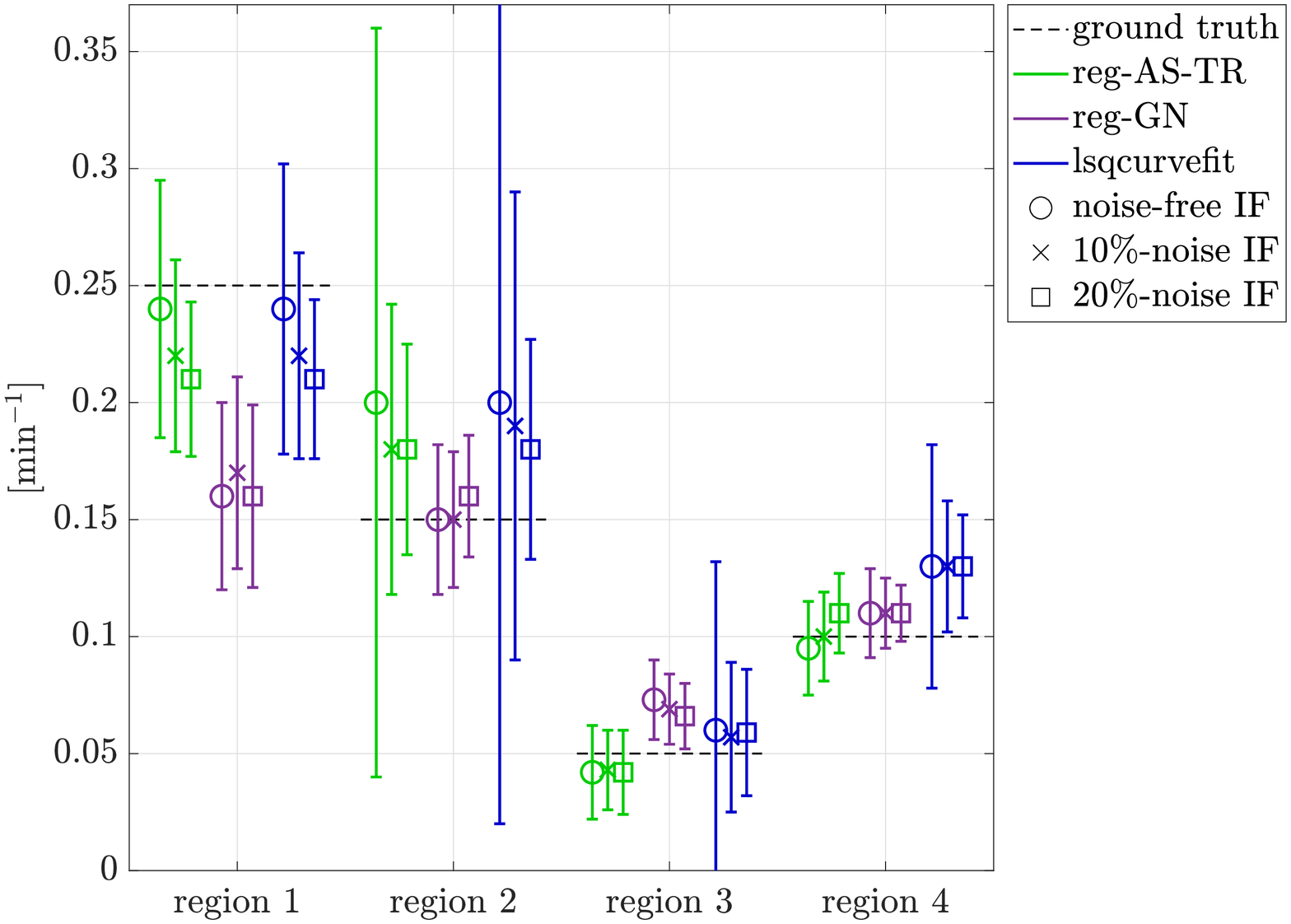}} \\
	\subfigure[$k_3$ \label{subfig:k3-rec}]
	{\includegraphics[width=0.35\textwidth]{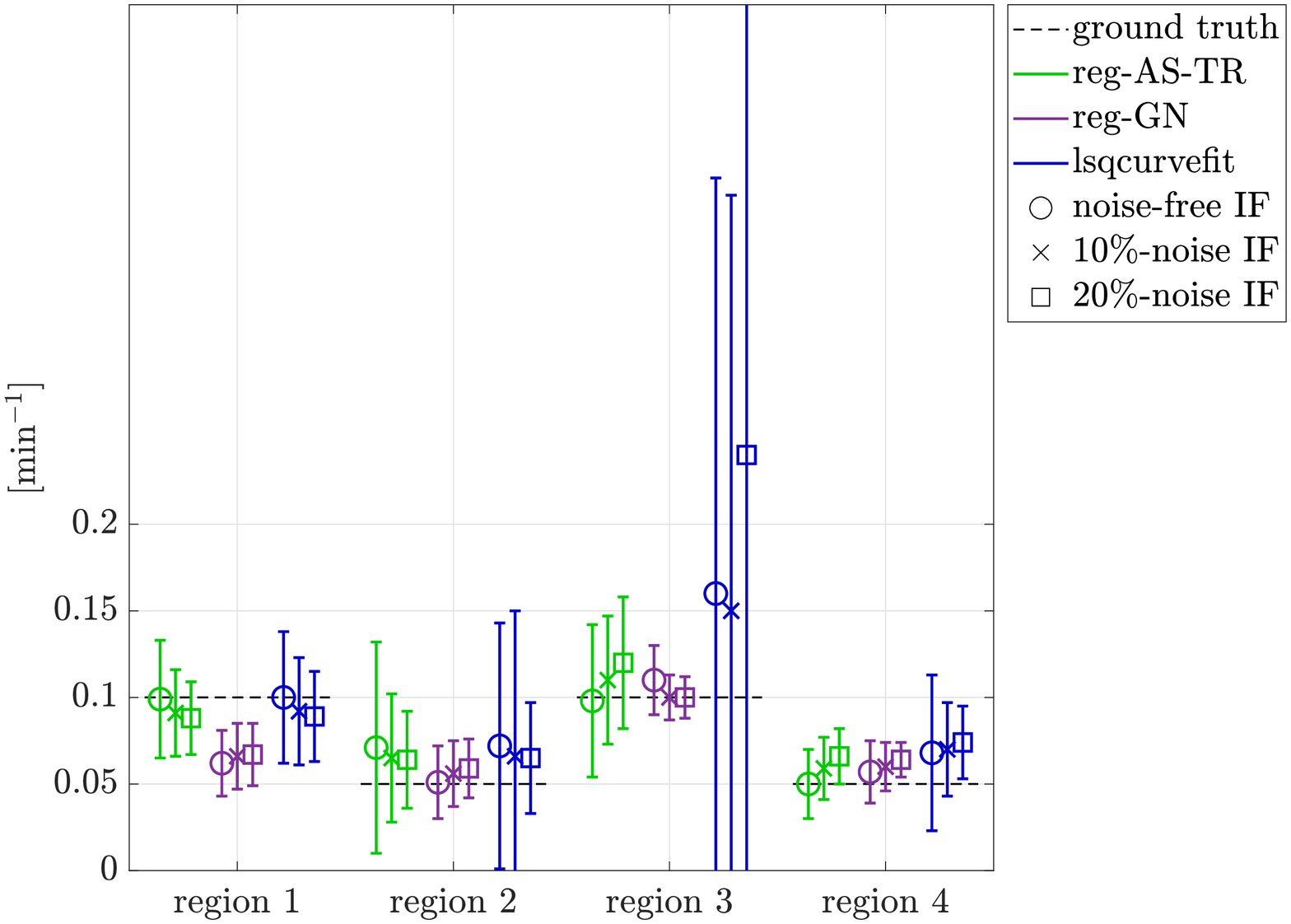}} \quad   \hspace{0.7cm}
	\subfigure[$k_4$ \label{subfig:k4-rec}]
	{\includegraphics[width=0.35\textwidth]{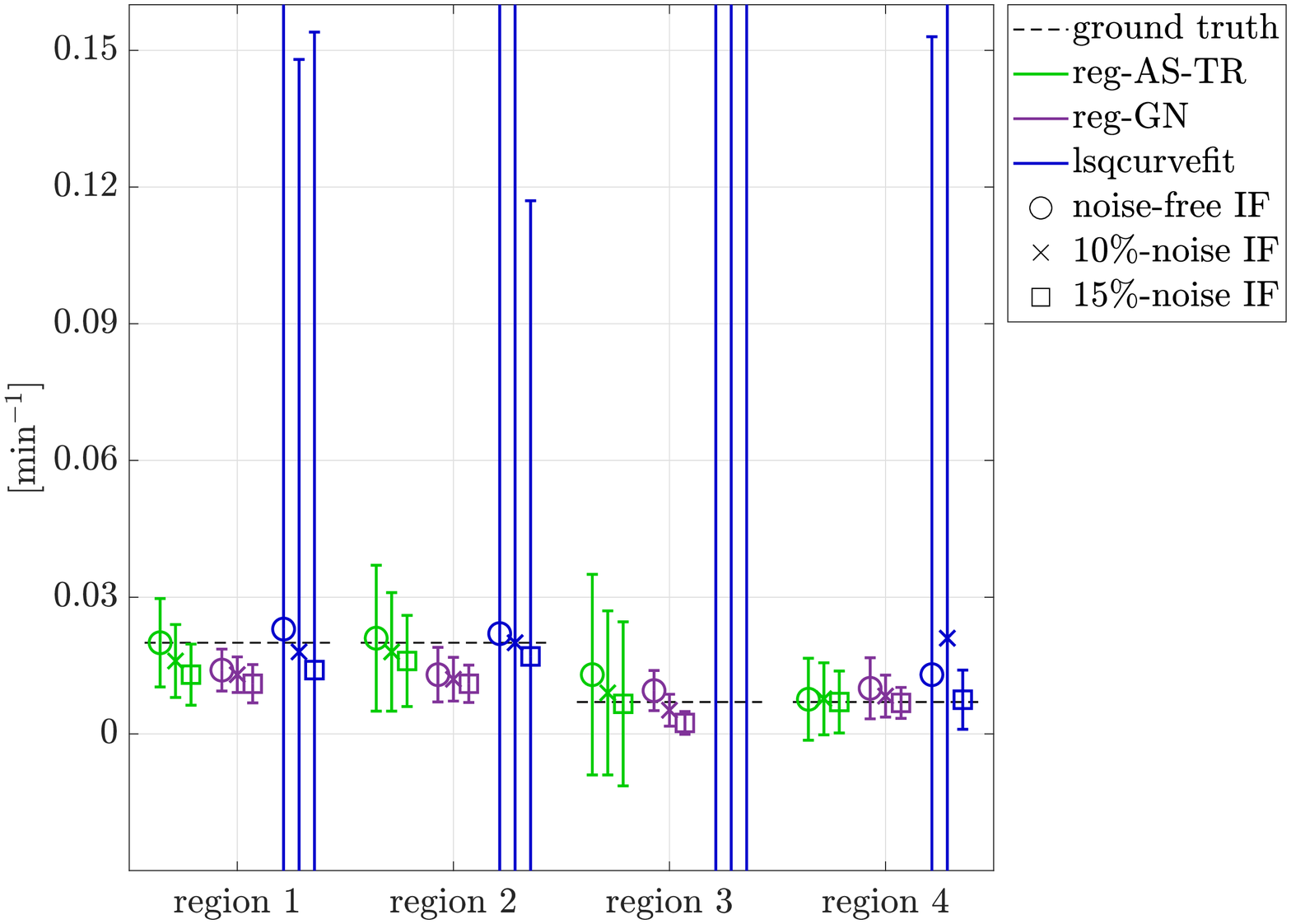}}
	\caption{Mean and standard deviation values of the kinetic parameters for the four homogeneous region as error bars over the reconstructions: reg-AS-TR (green bars), reg-GN (purple bars), \emph{lsqcurvefit} (blue bars); noise-free IF (circle), $10\%$-noise IF (cross), $20\%$-noise IF (square).}
	\label{fig:k-rec-all}
\end{figure}

\begin{figure}
	\centering
	\includegraphics[width=0.2\textwidth]{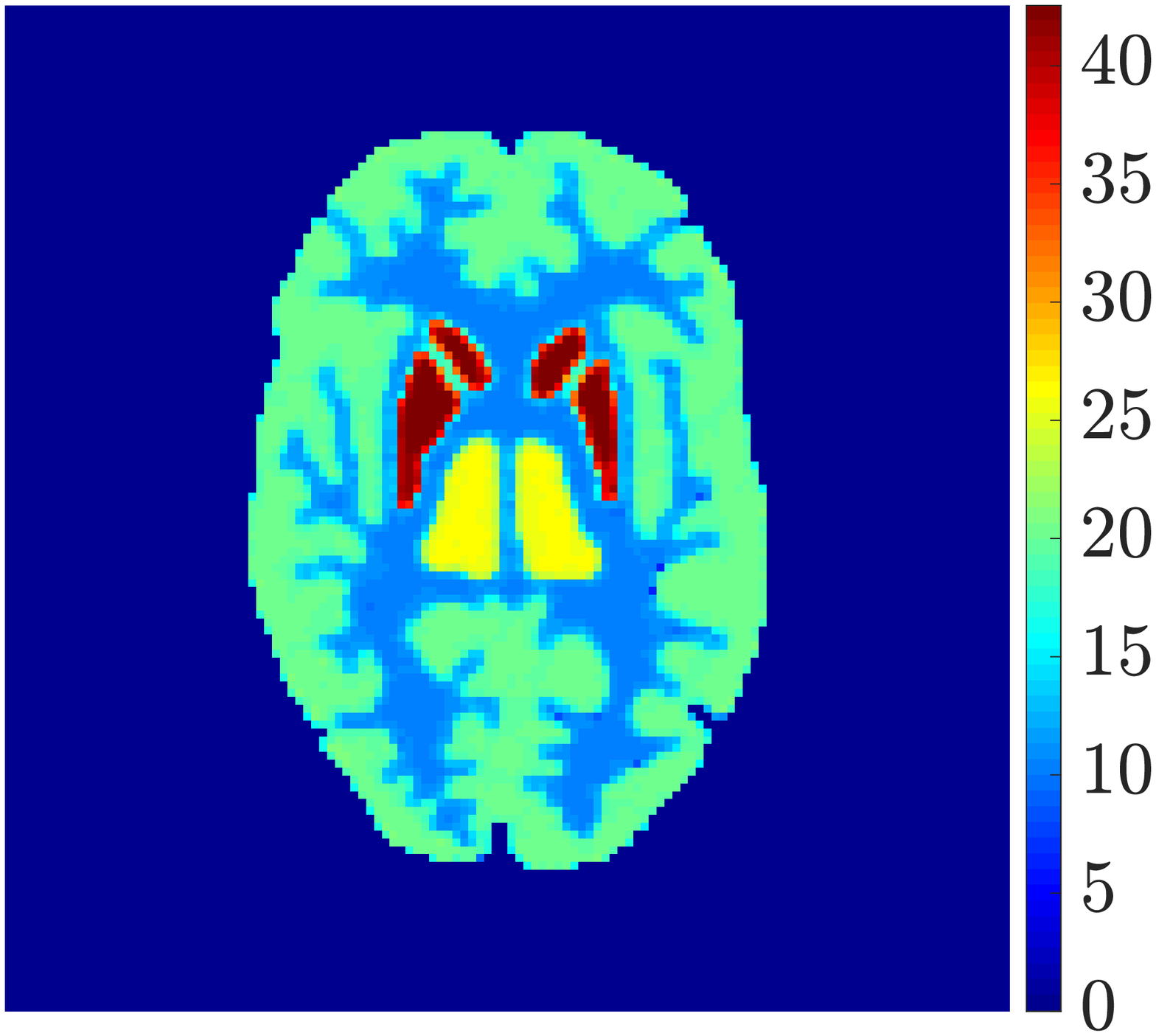}\quad   \hspace{1cm}
	\includegraphics[width=0.2\textwidth]{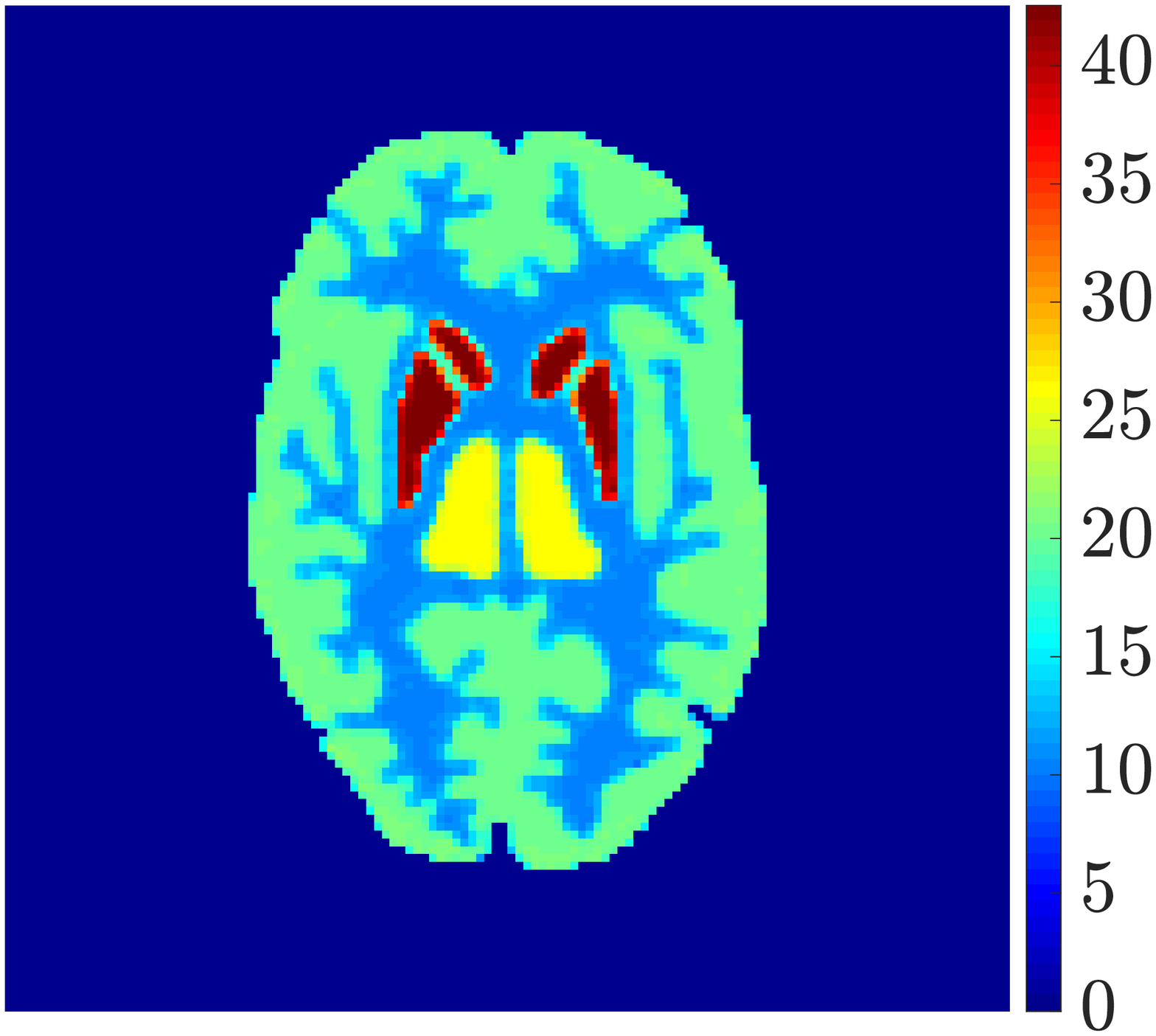}\quad  \hspace{1cm}
	\includegraphics[width=0.2\textwidth]{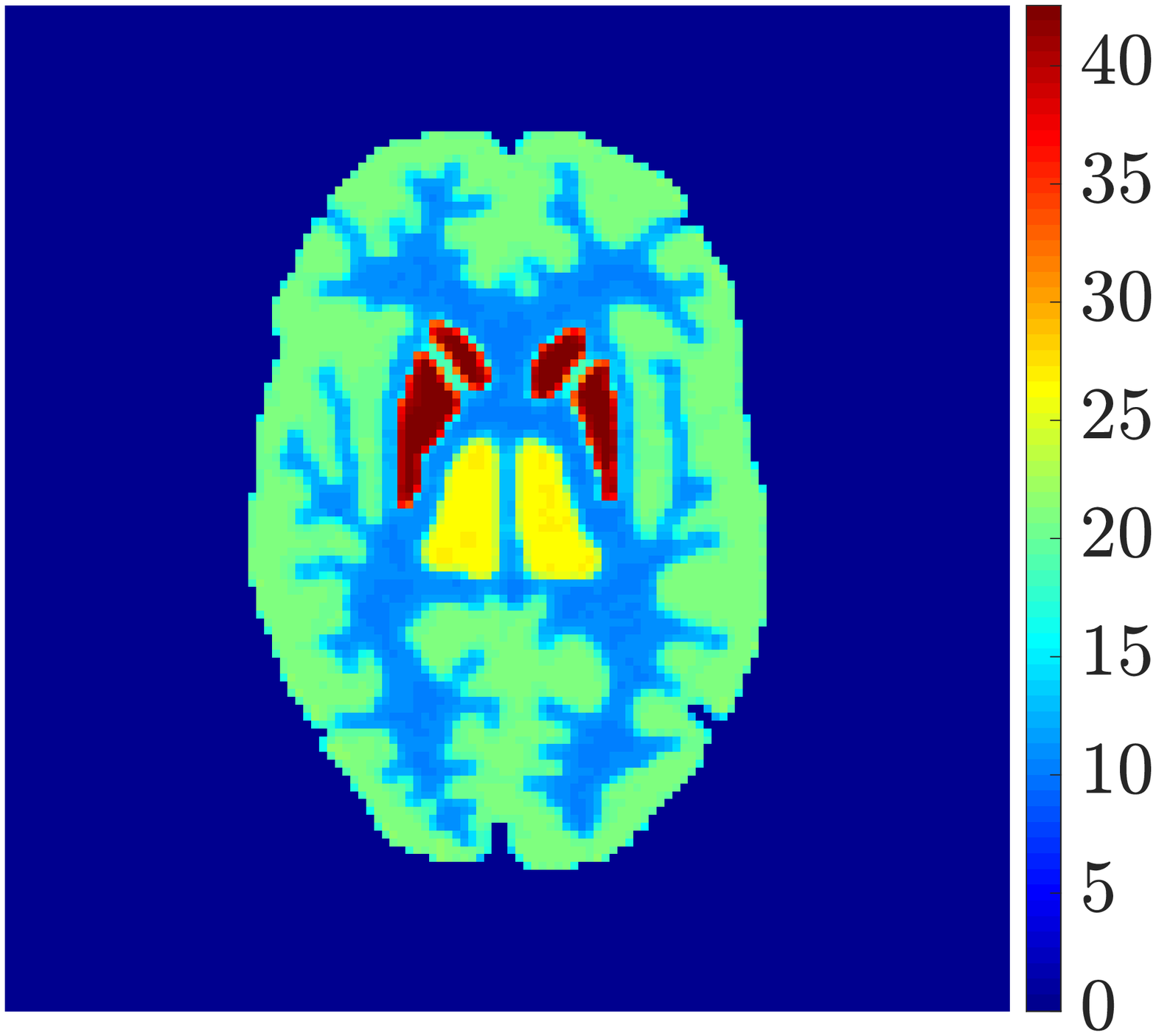}\\
	
	\includegraphics[width=0.2\textwidth]{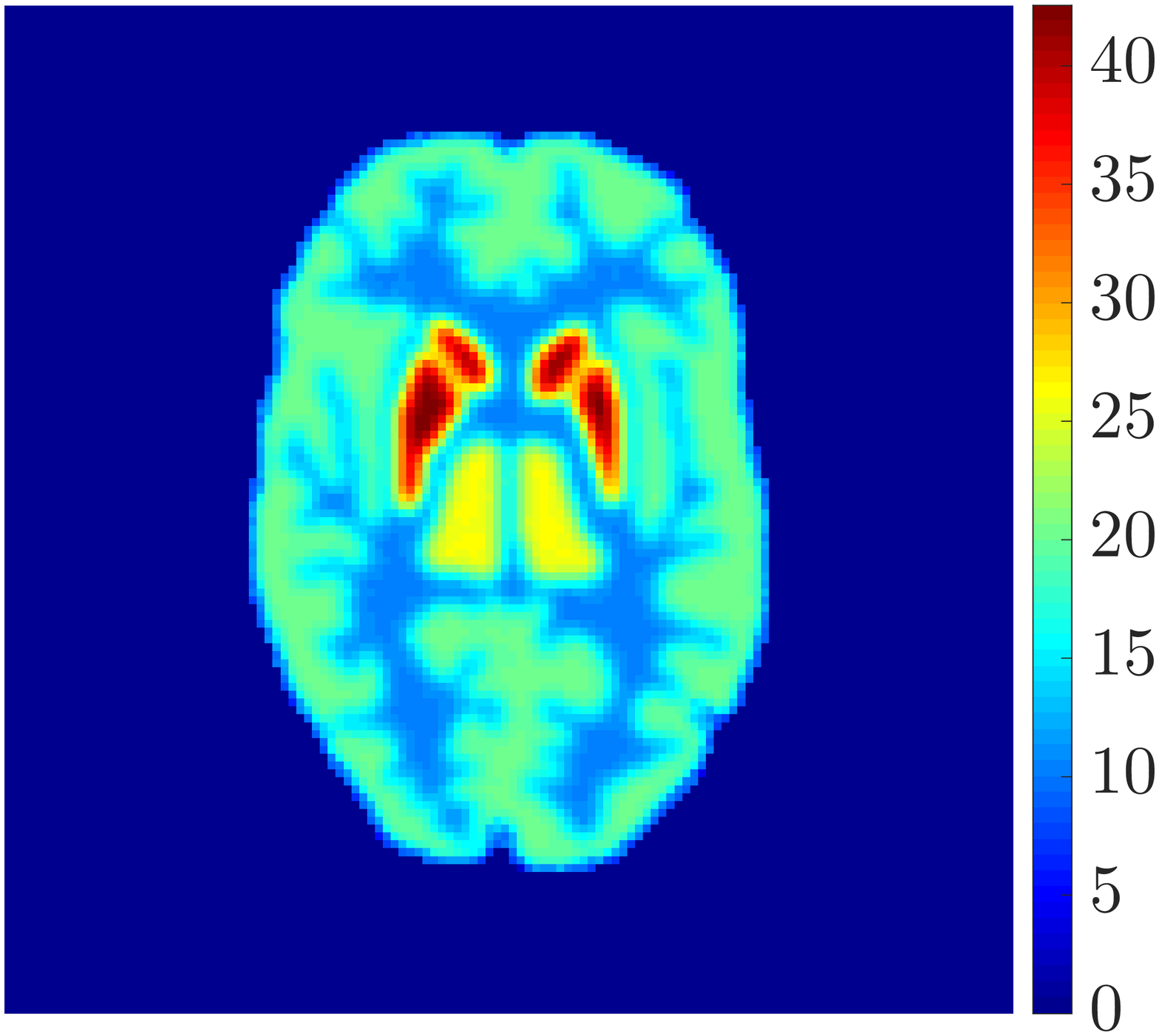}\quad  \hspace{1cm}
    \includegraphics[width=0.2\textwidth]{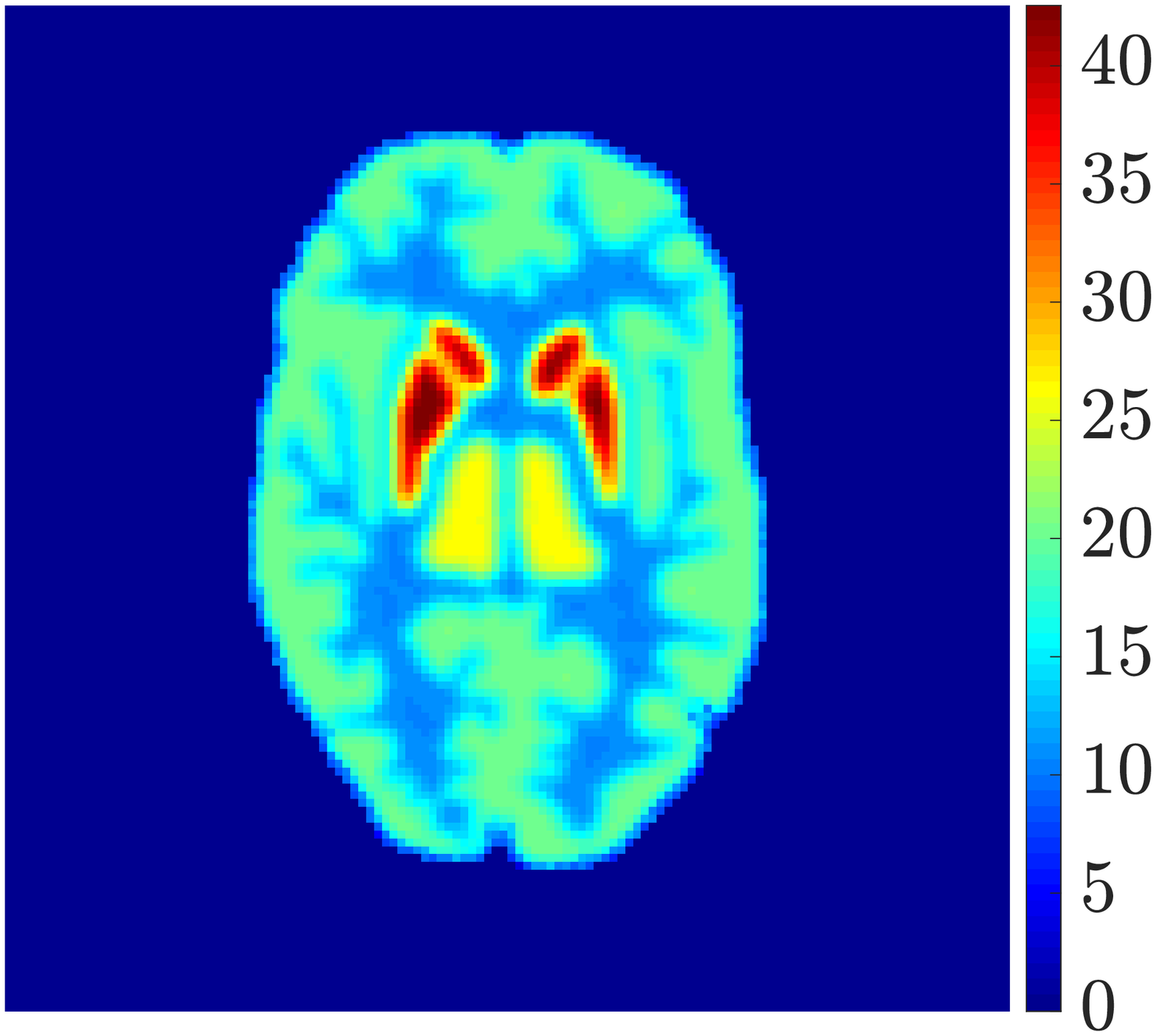}\quad  \hspace{1cm}
    \includegraphics[width=0.2\textwidth]{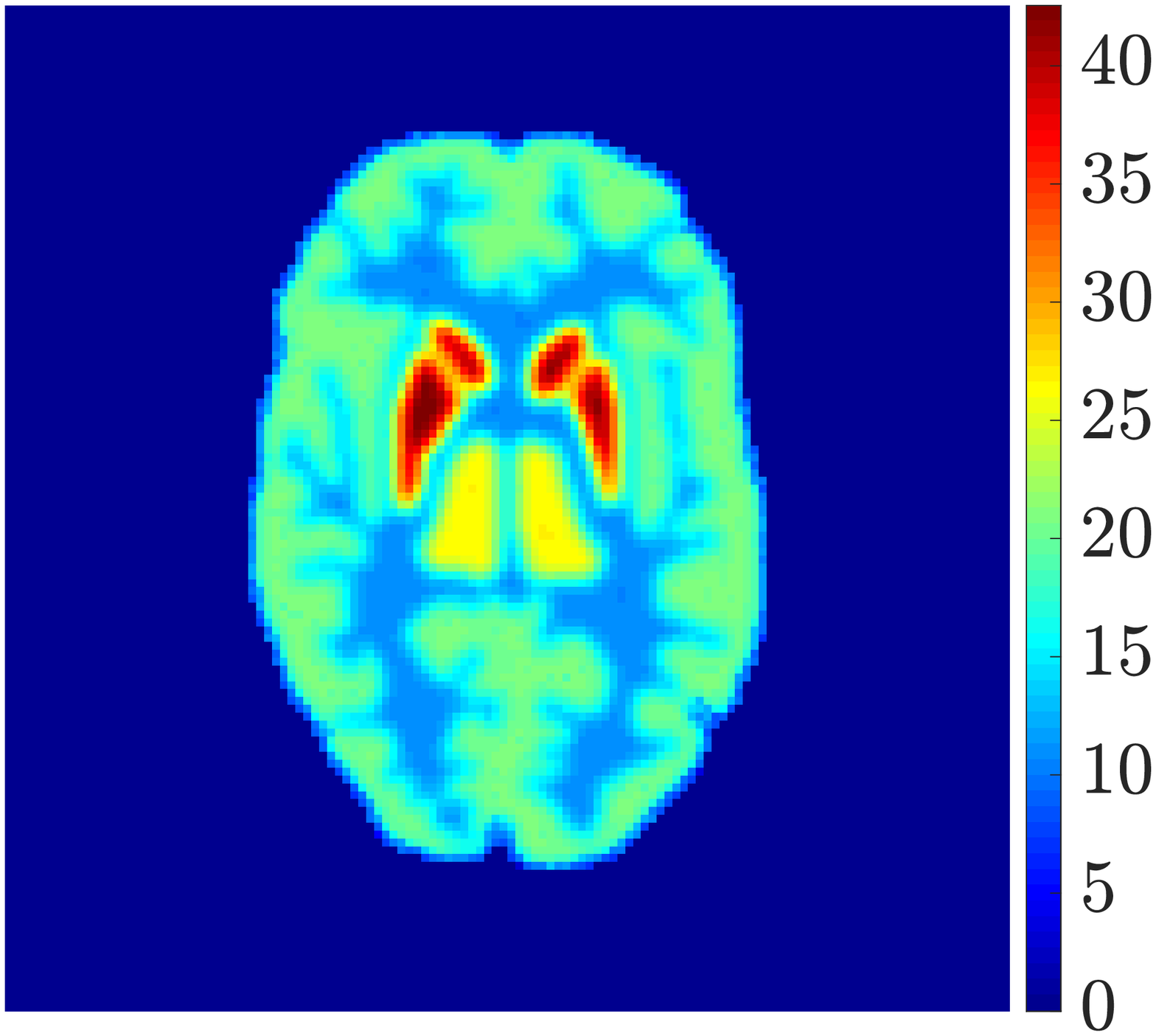}\\ 
	
	\includegraphics[width=0.2\textwidth]{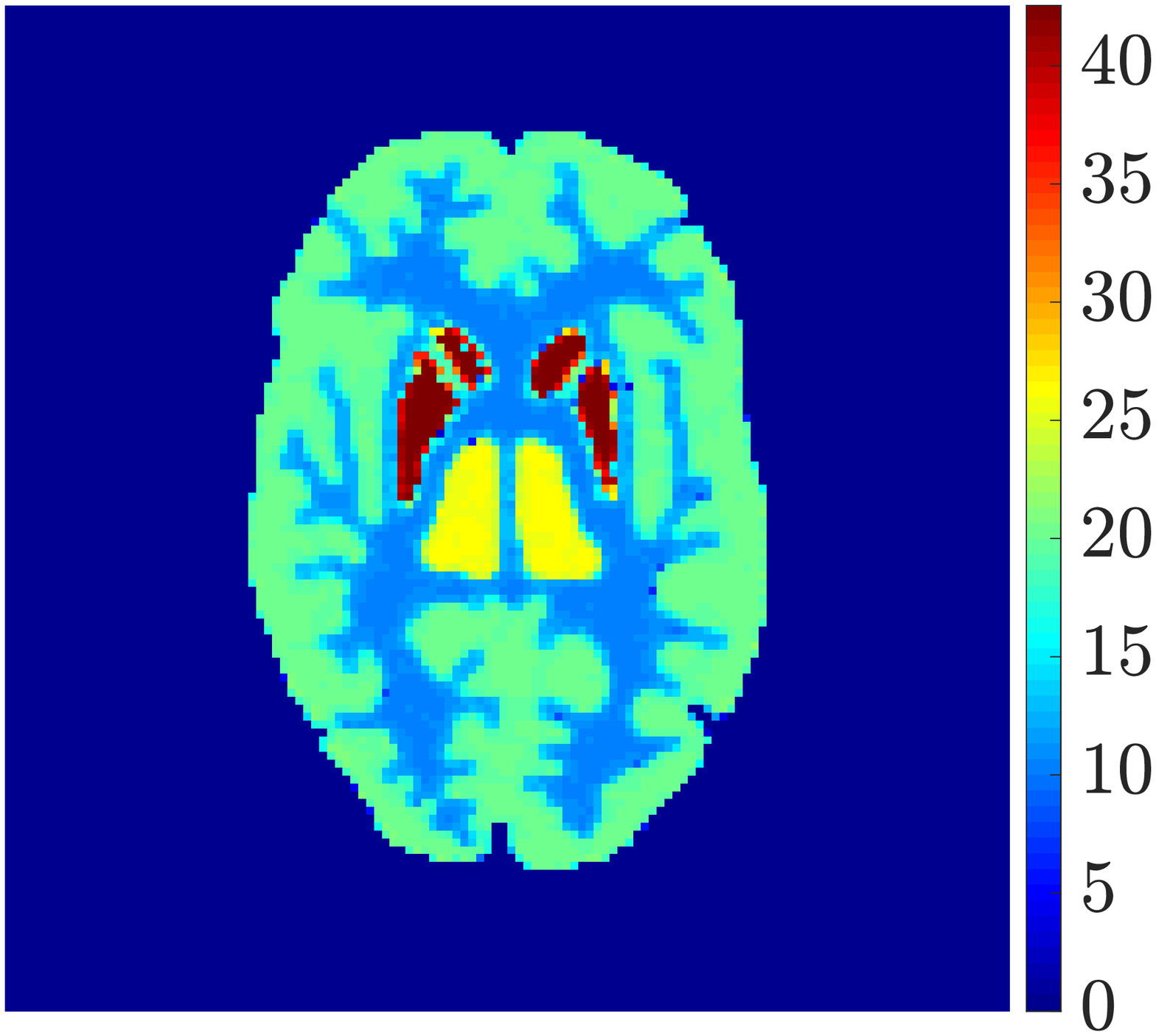}\quad   \hspace{1cm}
    \includegraphics[width=0.2\textwidth]{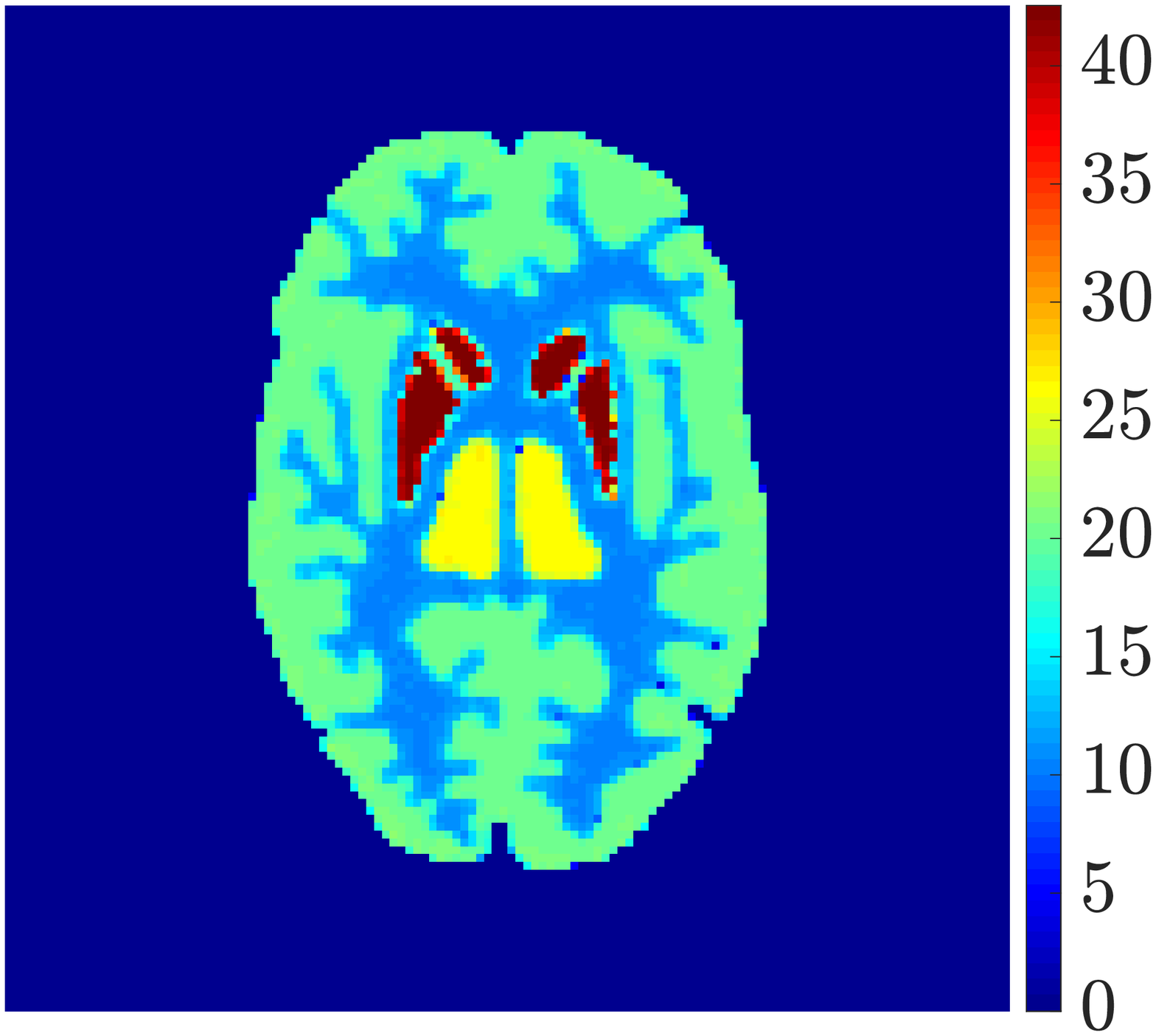}\quad   \hspace{1cm}
    \includegraphics[width=0.2\textwidth]{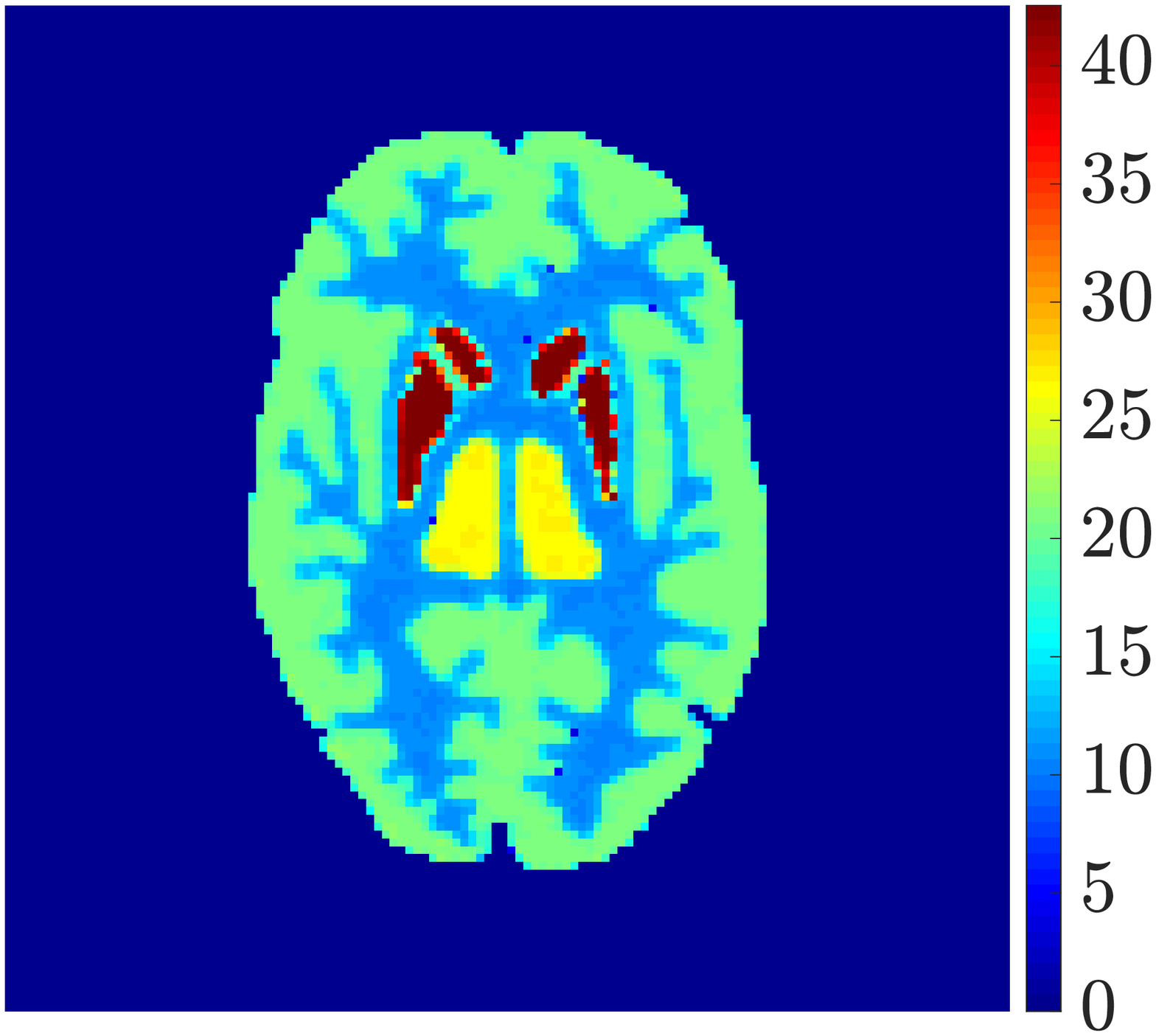}
  \caption{Last frame of the dynamic reconstruction obtained by using reg-AS-TR (first row), reg-GN (second row), and \emph{lsqcurvefit} (third row). From left to right: noise-free IF, $10\%$-noise IF, $20\%$-noise IF.}
  \label{fig:data-rec}
\end{figure}

In general, reg-AS-TR and \emph{lsqcurvefit} seem to provide similar mean reconstructions, although uncertainties associated to \emph{lsqcurvefit} are significantly bigger. On the other hand reg-GN seems to systematically underestimate the parameter values within region~$1$. Furthermore and as expected, for all methods the quality of the parametric reconstructions deteriorates with increasing noise levels; this is more clear from the $k_3$ and $k_4$ parametric images, probably due to the different sensitivities of the data with respect to the model parameters \cite{scetal19}. In reg-GN and \emph{lsqcurvefit} some artifacts can be observed at the edges of the homogeneous regions, especially around region $1$ and region $2$, whereas the effect of regularization in reg-AS-TR results in a reduced presence of artifacts while the structure of the regions is preserved. This general trend is confirmed by the error-bar plots of~\figurename~\ref{fig:k-rec-all}. Finally, the frames in \figurename~\ref{fig:data-rec} corresponding to reg-AS-TR show a significant improvement of the image quality with respect to what is provided by the other two approaches. 

The mean execution time for a single parametric reconstruction differs considerably between the reconstruction methods: reg-AS-TR requires about $20$ minutes, reg-GN needs a computational time in the range $75-120$ minutes with run time increasing with noise level on IF (as a consequence of the stopping criterion implemented) and Matlab \emph{lsqcurvefit} takes about $90$ minutes. Therefore reg-AS-TR seems to be the most efficient approach in terms of both computational time and reconstruction accuracy.

Next steps for this piece of research activity will be the validation of reg-AS-TR against several experimental datasets in the case of both humans' and small animals' dynamic PET images. Further, we are going to generalize reg-AS-TR to the case of more complex compartmental models like the ones for the assessment of FDG kinetics in liver \cite{gaetal15} and kidneys \cite{gaetal14}.

\newpage
	\section*{Appendix}
		\appendix
		\setcounter{section}{1}
	In this Appendix, we provide the proof of Proposition \ref{mono_proj}.
	\begin{proof}
		Let ${\ve p^j}^\dag  = \ve k^\dag - \ve k^j$ and $\ve k^{j+1}= \ve k^j + \ve {\bar p}^j$, with $\ve {\bar p}^j$ given in \eqref{proj_step}.
		If $\ve{\bar{p}}^j=\ve p^j$, we return back to the unconstrained case, for which Proposition 2.1 in \cite{Wang} holds.
		Let assume for some $j$ that
		$$\ve{\bar{p}}^j=t(\Pi(\ve k^j+\ve p^j)-\ve{k^j})\quad \quad 0<t<1.$$
		Let $\ve{\bar{k}}^j=\Pi(\ve{k^j}+\ve{p^j})$.  From the properties of the projection operator \cite[Proposition  2.1.3]{bertsekas},
we have
		\begin{equation}\label{prop-proj}
		(\ve{\bar{k}}^j-(\ve k^j+\ve p^j))^T (\ve k- \ve{\bar{k}}^j)\geq 0 \quad \forall\ \ve{k}\geq 0;
		\end{equation}
		in particular, Eq.~\eqref{prop-proj} holds for $\ve k^\dag\geq 0$ and, therefore, we obtain
		\begin{equation}
		(\ve{\bar{k}}^j-\ve{k^j})^T (\ve k^\dag- \ve{\bar{k}}^j)\geq {\ve p^j}^T (\ve k^\dag- \ve{\bar{k}}^j).\nonumber
		\end{equation}
		Setting $\ve{\hat{p}}^j=\ve{\bar{k}}^j-\ve k^j=\frac{1}{t} \ve{\bar{p}}^j$ and $\ve k^\dag- \ve{\bar{k}}^j= {\ve p^j}^\dag- \ve{\hat p}^j$, we can write the previous inequality as follows:
		\begin{equation}
		{\ve{\hat p}^j}^T {\ve p^j}^\dag \geq \|\ve{\hat p}^j\|^2 +  {\ve p^j}^T {\ve p^j}^\dag -  {\ve p^j}^T \ve{\hat{p}}^j. \nonumber
		\end{equation}
		From the identity $\| {\ve p^j}^\dag\|^2 -\|\ve{\bar{p}}^j - {\ve p^j}^\dag \|^2 =  2\ {\ve{\bar{p}}^j}^T{\ve p^j}^\dag -\|\ve{\bar{p}}^j\|^2$,
		the definition  $\ve{\bar{p}}^j=t\ve{\hat{p}}^j$, the previous inequality and $t\in (0,1)$,  we have
		\begin{eqnarray}
		\| {\ve p^j}^\dag\|^2 -\|\ve{\bar{p}}^j - {\ve p^j}^\dag \|^2 &=&  2\ {\ve{\bar{p}}^j}^T{\ve p^j}^\dag -\|\ve{\bar{p}}^j\|^2 =
		2t\  {\ve{\hat p}^j}^T {\ve p^j}^\dag - t^2\|\ve{\hat p}^j\|\geq \nonumber \\
		&\geq& 2t\ (  \|\ve{\hat p}^j\|^2 +  {\ve p^j}^T {\ve p^j}^\dag -  {\ve p^j}^T \ve{\hat{p}}^j)- t^2\|\ve{\hat p}^j\|  \nonumber\\
		&=& t\left( (2-t)  \|\ve{\hat p}^j\|^2 + 2{\ve p^j}^T {\ve p^j}^\dag -2 {\ve p^j}^T \ve{\hat{p}}^j  \right) \nonumber\\
		&>& t ( \|\ve{\hat p}^j\|^2 + 2{\ve p^j}^T {\ve p^j}^\dag + \|{\ve p^j}\|^2 - \|{\ve p^j}\|^2 -2 {\ve p^j}^T \ve{\hat{p}}^j)  \nonumber\\
		&=& t ( \|\ve{\hat p}^j - \ve p^j \|^2 +  2{\ve p^j}^T {\ve p^j}^\dag - \|{\ve p^j}\|^2  )  \nonumber \\
		&>&  t ( 2{\ve p^j}^T {\ve p^j}^\dag - \|{\ve p^j}\|^2)\nonumber \\
		&>&  t ( 2{\ve p^j}^T {\ve p^j}^\dag - 2\|{\ve p^j}\|^2 ). \label{ultima}
		\end{eqnarray}
		Now, we recall 
		that, in view of positive definiteness of the matrix $(\ve J(\ve k^j)^T \ve J(\ve k^j) +\alpha_j \ve I_n)$, the following matrix identities
		hold:
		{\footnotesize
			\begin{eqnarray*}
				&&  (\ve J(\ve k^j)^T \ve J(\ve k^j) +\alpha_j \ve I_n)^{-1} \ve J(\ve k^j)^T =\ve J(\ve k^j)^T (\ve J(\ve k^j) \ve J(\ve k^j)^T +\alpha_j \ve I_N  )^{-1} \\
				&&  \ve I_N- \ve J(\ve k^j) \ve J(\ve k^j)^T (\ve J(\ve k^j) \ve J(\ve k^j)^T +\alpha_j \ve I_N  )^{-1} = \alpha_j (\ve J(\ve k^j) \ve J(\ve k^j)^T +\alpha_j \ve I_N  )^{-1} 
			\end{eqnarray*}
		}
		As a consequence, setting $\ve r^j= \ve y^\delta -\ve F(\ve k^j)$, we can write
		\begin{eqnarray}
		&&  \ve r^j- \ve J(\ve k^j) \ve p^j = \alpha_j (\ve J(\ve k^j) \ve J(\ve k^j)^T +\alpha_j \ve I_N  )^{-1} \ve r^j. \label{tre}
		\end{eqnarray}
		In view of inequality \eqref{ultima}, the definition of $\ve p^j$ and the above identities, we have
		{\footnotesize
			\begin{eqnarray*}
				&& \| {\ve p^j}^\dag\|^2 -\|\ve{\bar{p}}^j - {\ve p^j}^\dag \|^2 \nonumber \\
				&&>t ( 2 {\ve r^j}^T(\ve J(\ve k^j) \ve J(\ve k^j)^T + \alpha_j \ve I_N)^{-1} \ve J(\ve k^j) {\ve p^j}^\dag
				- 2 ( \ve J(\ve k^j) {\ve p^j})^T (\ve J(\ve k^j) \ve J(\ve k^j)^T + \alpha_j \ve I_N)^{-1} \ve r^j) \nonumber\\
				&&= t ( 2 {\ve r^j}^T(\ve J(\ve k^j) \ve J(\ve k^j)^T + \alpha_j \ve I_N)^{-1} \ve J(\ve k^j) {\ve p^j}^\dag -2 {\ve{r}^j}^T (\ve J(\ve k^j)
				\ve J(\ve k^j)^T+\alpha_j \ve I_N)^{-1} \ve r^j +\nonumber\\
				&&+ 2\alpha_j {\ve r^j}^T  (\ve J(\ve k^j)  \ve J(\ve k^j) ^T + \alpha_j \ve I_N)^{-2}\ve r^j)\nonumber\\
				&&=t ( 2\alpha_j \| (\ve J(\ve k^j)  \ve J(\ve k^j) ^T + \alpha_j \ve I_N)^{-1} {\ve r^j}\|^2 
				-  2 {\ve r^j}^T (\ve J(\ve k^j)  \ve J(\ve k^j) ^T+\alpha_j \ve I_N)^{-1} (  \ve r^j -\ve J(\ve k^j)  {\ve p^j}^\dag ) )\nonumber\\
				&&>2t ( \alpha_j\|  (\ve J(\ve k^j)   \ve J(\ve k^j)  ^T + \alpha_j \ve I_N)^{-1} {\ve r^j}\|^2 
				-  \| (\ve J(\ve k^j) \ve J(\ve k^j)^T+\alpha_j \ve I_N)^{-1}\ve r^j\| \|\ve r^j -\ve J(\ve k^j) {\ve p^j}^\dag \| )\nonumber\\
				&&= 2t (\|\ve r^j- \ve J(\ve k^j)\ve p^j\|- \|\ve r^j -\ve J(\ve k^j) {\ve p^j}^\dag \|  ) \| (\ve J(\ve k^j) \ve J(\ve k^j)^T+\alpha_j \ve I_N)^{-1}\ve r^j\| \nonumber ,
			\end{eqnarray*}
		}
		where the last inequality follows from the Cauchy-Schwarz inequality.
		Then, the $q$-condition and the assumption \eqref{gamma_cond1} yields
		\begin{eqnarray*}
			\|{\ve p^j}^\dag\|^2 - \|\ve{\bar{p}}^j - {\ve p^j}^\dag \|^2  &>&\frac{2t(\gamma_\delta-1)q}{\gamma_\delta}\|\ve{y}^\delta-\ve F(\ve k^j)\|
			\|\ve v^j\|.\nonumber
		\end{eqnarray*}	
	\end{proof}


\bibliographystyle{plain}
\bibliography{mybiblio}

\end{document}